\documentclass[10pt]{article}

\usepackage{amsmath}
\usepackage{amsfonts}
\usepackage{amssymb}
\usepackage{framed}
\usepackage{enumerate}
\usepackage{mathrsfs}
\usepackage[hidelinks]{hyperref}
\usepackage{enumitem}
\usepackage{mathrsfs}
\usepackage{scrextend}
\usepackage{amsthm}
\usepackage{bbm}
\usepackage{thmtools}
\usepackage{thm-restate}
\usepackage{mathabx}
\usepackage[margin=1.5in]{geometry}
\usepackage{rotating}

\usepackage{xcolor}
\hypersetup{
    colorlinks,
    linkcolor={red!50!black},
    citecolor={blue!50!black},
    urlcolor={blue!80!black}
}

\usepackage{tikz-cd}

\newcommand{\res}{\textrm{res}}

\theoremstyle{plain}
\newtheorem{thm}{Theorem}[section]
\newtheorem{lem}[thm]{Lemma}
\newtheorem{prop}[thm]{Proposition}
\newtheorem{cor}[thm]{Corollary}

\theoremstyle{definition}
\newtheorem{defn}[thm]{Definition}

\theoremstyle{remark}

\newtheorem*{notn}{Notation}

\tikzset{
  symbol/.style={
    draw=none,
    every to/.append style={
      edge node={node [sloped, allow upside down, auto=false]{$#1$}}}
  }
}

\newcommand{\Spec}{\textrm{Spec} \hspace{0.15em} }

\newcommand\restr[2]{{
	\left.\kern-\nulldelimiterspace
	#1
	\vphantom{\big|}
	\right|_{#2}
	}}
\newcommand{\an}[1]{#1^{\textrm{an}}}

\newcommand{\ep}{\varepsilon}

\newcommand{\FSch}{\textrm{FSch}}
\newcommand{\Set}{\textrm{Set}}

\newcommand{\Hom}{\textrm{Hom}}

\newcommand{\ch}[1]{\widecheck{{#1}}}

\newcommand{\GL}{\textrm{GL}}

\usepackage{amsmath,calligra,mathrsfs}

\usepackage[cal=boondoxo]{mathalfa}

\DeclareMathOperator{\sheafhom}{\mathcal{H \kern -1pt o \kern -2pt m}}
\DeclareMathOperator{\sheafend}{\mathcal{E \kern -1pt n \kern -2pt d}}

\title{Sets of Special Subvarieties of Bounded Degree}
\author{David Urbanik}

\begin{document}

\maketitle

\begin{abstract}
Let $f : X \to S$ be a family of smooth projective algebraic varieties over a smooth connected quasi-projective base $S$, and let $\mathbb{V} = R^{2k} f_{*} \mathbb{Z}(k)$ be the integral variation of Hodge structure coming from degree $2k$ cohomology it induces. Associated to $\mathbb{V}$ one has the so-called \emph{Hodge locus} $\textrm{HL}(S) \subset S$, which is a countable union of ``special'' algebraic subvarieties of $S$ parametrizing those fibres of $\mathbb{V}$ possessing extra Hodge tensors (and so conjecturally, those fibres of $f$ possessing extra algebraic cycles). The special subvarieties belong to a larger class of so-called weakly special subvarieties, which are subvarieties of $S$ maximal for their algebraic monodromy groups. For each positive integer $d$, we give an algorithm to compute the set of all weakly special subvarieties $Z \subset S$ of degree at most $d$ (with the degree taken relative to a choice of projective compactification $S \subset \overline{S}$ and very ample line bundle $\mathcal{L}$ on $\overline{S}$). As a corollary of our algorithm we prove conjectures of Daw-Ren and Daw-Javanpeykar-K\"uhne on the finiteness of sets of special and weakly special subvarieties of bounded degree.
\end{abstract}

\tableofcontents

\section{Introduction}

Suppose that $f : X \to S$ is a smooth projective morphism of quasi-projective varieties, with $S$ smooth and connected. We may view $f$ as an algebraic family of projective varieties, where a point $s$ in the parameter space $S$ gives rise to the projective fibre $X_{s}$. A general problem which encompasses many questions of classical interest is to understand what conditions on $s$ are imposed by the requirement that $X_{s}$ (or its Cartesian powers $X_{s} \times \cdots \times X_{s}$) possess ``more'' algebraic cycles than a general fibre of $f$. 

Under the Hodge conjecture, this problem can be reformulated in the following manner. Let $\mathbb{V} = R^{2k} f^{\textrm{an}}_{*} \mathbb{Z}(k)$ be the local system on $S$ whose fibres $\mathbb{V}_{s}$ are naturally identified with $H^{2k}(X_{s}, \mathbb{Z}(k))$. Then $\mathcal{H} = \mathbb{V} \otimes \mathcal{O}_{\an{S}}$ is a vector bundle on $S$, whose fibres may be naturally identified with $H^{2k}(X_{s}, \mathbb{C})$. Classical Hodge theory endows $\mathcal{H}$ with a decreasing filtration $F^{\bullet}$ by vector subbundles, making the data of $(\mathbb{V}, F^{\bullet})$ into a so-called \emph{variation of Hodge structure}, and the Hodge conjecture then predicts that the classes in $\mathcal{H}_{s}$ induced by algebraic cycles are exactly those in $\mathbb{V}_{\mathbb{Q}, s} \cap F^{0}$. We may thus reformulate our question as: what conditions on $s$ are imposed by the requirement that the fibres $\mathbb{V}_{\mathbb{Q}, s}$ (or their tensor powers) have more classes lying in $F^{0}$ than at a general point of $S$?

In this setting, the locus in $S$ determined by this condition is called the Hodge locus, commonly denoted $\textrm{HL}(S) \subset S$. It is known by a result of \cite{CDK} to be a countable union of closed algebraic subvarieties of $S$. The individual irreducible components of $\textrm{HL}(S)$ are commonly referred to as (maximal proper) special subvarieties in view of the relationship between such components and the notion of special subvariety arising in the theory of Shimura varieties. More precisely, we may describe a collection of ``special'' subvarieties of $S$ associated to $\mathbb{V}$ in the following manner:

\begin{defn}
Given a (pure integral weight zero) Hodge structure $(V, F^{\bullet})$, we may define the \emph{Mumford-Tate group} of $(V, F^{\bullet})$ to be the stabilizer in $\textrm{GL}(V_{\mathbb{Q}})$ of all \emph{Hodge tensors} associated to $(V, F^{\bullet})$, i.e., the stabilizer of the rational tensors 
\[ \textrm{Hg}(V) = \left( \bigoplus_{i,j \geq 0} V^{\otimes i}_{\mathbb{Q}} \otimes (V^{*}_{\mathbb{Q}})^{\otimes j} \right) \cap F^{0} . \] 
Given a subvariety $Z \subset S$ we define the Mumford-Tate group of $Z$ to be the Mumford-Tate group at a generic point of $Z$, i.e., at a point $s \in Z$ which lies outside of $\textrm{HL}(Z)$.
\end{defn}

\begin{defn}
A \emph{special subvariety} $Y \subset S$ is a geometrically irreducible algebraic subvariety maximal under inclusion for its Mumford-Tate group (c.f. \cite[Def. 1.2]{closurepositivelocus}).
\end{defn}

The special subvarieties belong to a larger class of subvarieties of $S$ which are termed ``weakly special'', defined in an analogous manner as follows: 

\begin{defn}
Given a variety $Z \subset S$, its algebraic monodromy group is the identity component of the Zariski closure of the monodromy representation associated to $\restr{\mathbb{V}}{Z^{\textrm{nor}}}$, where $Z^{\textrm{nor}} \to Z$ is the normalization.
\end{defn}

\begin{defn}
\label{weakspdef}
A \emph{weakly special subvariety} $Z \subset S$ is a geometrically irreducible algebraic subvariety maximal under inclusion for its algebraic monodromy group (c.f. \cite[Cor. 3.14]{closurepositivelocus}).
\end{defn}

\noindent It is known that the algebraic monodromy group of a variety $Z$ is always a $\mathbb{Q}$-normal subgroup of its Mumford-Tate group by a result of Andr\'e-Deligne \cite{Andre1992}, and that any special subvariety of $S$ is in fact weakly special \cite[Def. 3.1]{closurepositivelocus}. Thus studying the set of weakly special subvarieties is a generalization of the problem of studying the set of special subvarieties.

Our central result, from which the rest of the results in this paper will follow, is the following:

\begin{defn}
Given a variation of integral Hodge structure $\mathbb{V}$, we say say that $\mathbb{V}$ is polarizable if there exists a map $Q : \mathbb{V} \otimes \mathbb{V} \to \mathbb{Z}$ of local systems which fibrewise induces a polarization on the Hodge structure $\mathbb{V}_{s}$ for each $s$.
\end{defn}

\begin{thm}
\label{mainthm}
Let $(\mathbb{V}, F^{\bullet})$ be a polarizable variation of Hodge structure on a smooth quasi-projective algebraic variety $S$, and fix a projective compactification $S \subset \overline{S}$ with ample line bundle $\mathcal{L}$ on $\overline{S}$. Given a subvariety $Z \subset S$, write $\deg Z$ for the degree of its closure $\overline{Z}$ in $\overline{S}$ with respect to $\mathcal{L}$. Then for any positive integer $d$, there exists a terminating algorithm which computes, as a constructible locus in the Hilbert scheme $\textrm{Hilb}(\overline{S})$, the set of all weakly special subvarieties of $S$ of degree at most $d$.
\end{thm}

\noindent An immediate question that is raised by the statement of \autoref{mainthm} is how one plans to computationally represent the data $(\mathbb{V}, F^{\bullet})$. In the geometric setting, for instance, this representation will be in terms of the algebraic de Rham cohomology and the algebraic incarnation of the Gauss-Manin connection. We will elabourate on the computational model in which we work in \autoref{compmodelsec}, but for now let us describe how we interpret the claim that the algorithm ``computes the set of all weakly special subvarieties''. As this is an uncountable set (for instance, when $\mathbb{V}$ induces a quasi-finite period map, every point of $S$ is itself a zero-dimensional weakly special subvariety), we will need to provide a finitary description of it. What we actually do is consider the Hilbert scheme $\textrm{Hilb}(\overline{S})$ and the union $C_{d}$ of the finitely many components parametrizing (geometrically irreducible, reduced) subvarieties of $\overline{S}$ of degree at most $d$. We then consider the locus $\mathcal{W}_{d} \subset C_{d}$ consisting of points $[\overline{Z}]$ such that $Z = S \cap \overline{Z}$ is weakly special. We prove that:

\begin{thm}
\label{mainthm2}
The locus $\mathcal{W}_{d}$ is a constructible algebraic subset of $C_{d}$.
\end{thm}

\noindent We then interpret the claim of \autoref{mainthm} as saying that there exists an algorithm that outputs algebraic equations defining the locus $\mathcal{W}_{d}$ inside $C_{d}$.

\subsection{Applications} 

\subsubsection{Computing Rigid Specials}

Ultimately, one of our principal goals is to understand special subvarieties, not necessarily weakly special ones, so let us take some time to describe how information about special subvarieties can be gleaned from an algebraic description of the locus $\mathcal{W}_{d}$. Certainly some information is lost: in the case of a quasi-finite period map each point of $S$ is weakly special, and so we learn nothing about special points of $S$ this way. On the other hand, positive-dimensional\footnote{Here we really mean positive-dimensional in the Hodge-theoretic sense of \cite[Def. 1.3]{closurepositivelocus}, which coincides with being positive-dimensional as a subvariety in the quasi-finite period map case.} special subvarieties need not lie inside families of weakly specials; those special subvarieties which don't are referred to as ``weakly non-factor'' in \cite{fieldsofdef}, and for instance any positive-dimensional component of the Hodge locus of the universal family of degree $m$ hypersurfaces in $\mathbb{P}^{n+1}$ is of this type.\footnote{This is due to the simplicity of derived subgroup of the Mumford-Tate group of the variation, this latter fact being a consequence of the result that the global algebraic monodromy group is maximal due to \cite{beauville1986groupe}.}

We will see that the weakly non-factor special subvarieties all correspond to isolated points of $\mathcal{W}_{d}$, leading to the following refinement of \autoref{mainthm}:

\begin{cor}
\label{mainthmspcor}
In the situation of \autoref{mainthm}, there exists an algorithm which computes a finite set $\mathcal{S} \subset \mathcal{W}_{d}$ containing all the weakly non-factor special subvarieties of degree at most $d$.
\end{cor}

\noindent The set $\mathcal{S}$ referred to by \autoref{mainthmspcor} can simply be taken to be the set of isolated points in $\mathcal{W}_{d}$. The reason why not all such points correspond to weakly non-factor specials is because families of weakly specials need not be flat, and so points in $\mathcal{W}_{d}$ can be isolated for reasons unrelated to Hodge theory. Identifying which of the points of $\mathcal{S}$ correspond to special subvarieties would require further analysis, which appears difficult to carry out in full generality. 

\subsubsection{Conjectures on (Weakly) Special Degrees}

A second application of our algorithm, and of \autoref{mainthm2} in particular, is to two conjectures on the degrees of special subvarieties by Daw-Ren and Daw-Javanpeykar-K\"uhne. Fix a point $s \in S$ outside of $\textrm{HL}(S)$, let $\textrm{MT}(S) \subset \textrm{GL}(\mathbb{V}_{\mathbb{Q}, s})$ be the Mumford-Tate group of $\mathbb{V}_{s}$, and denote by $\Gamma_{S} \subset \textrm{GL}(\mathbb{V}_{s})$ the image of $\pi_{1}(S, s)$. For each subvariety $Z \subset S$, the algebraic monodromy group $\mathbf{H}_{Z}$ of $Z$ defines a $\Gamma_{S}$-conjugacy class of semisimple subgroup of $\textrm{MT}(S)$ by parallel translation. A conjecture of Daw-Ren is then contained in the following statement:

\begin{cor}
\label{finitesubgp}
As $Z \subset S$ ranges over all geometrically irreducible subvarieties of $S$ of degree at most $d$, the algebraic monodromy groups $\mathbf{H}_{Z}$ define only finitely many $\Gamma_{S}$-conjugacy classes in $\textrm{MT}(S)$.
\end{cor}

This corollary is a substantial strengthening of \cite[Conj. 10.4]{daw2018applications}, where it is stated in the special case of variations associated to Shimura varieties, and where $Z$ ranges only over the special subvarieties. Actually, Daw and Ren do not use the language of variations of Hodge structures, preferring instead the language of Shimura datum; we will give a brief translation between the two pictures when we prove \autoref{finitesubgp}. In their case, the corollary is of interest to implement a certain strategy for proving the Zilber-Pink conjecture, generalizing the Pila-Zannier strategy for the Andr\'e-Oort conjecture, which has had substantial success in recent years (c.f. \cite{klingler2018bi}). In particular, Daw and Ren prove several strong consequences of the Zilber-Pink conjecture (Theorem 14.2 and Theorem 14.3 in \cite{daw2018applications}) conditional on their Conjecture 10.4 (or rather Conjecture 10.3, which Conjecture 10.4 implies) and a collection of conjectures relating to heights and Galois orbits associated with special points. Our result therefore removes an obstacle to the implementation of Daw and Ren's generalization of the Pila-Zannier strategy.

In a paper \cite{daw2020effective} by Daw-Javanpeykar-K\"uhne a related conjecture is made in Remark 3.8, which is a special case of the following statement:

\begin{cor}
\label{finweaknonfactor}
For any positive integer $d$, there are finitely many weakly non-factor special subvarieties of $S$ of degree at most $d$.
\end{cor}

\noindent This is in fact contained in what we have already said, as it amounts to the statement that the set of (closures in $\overline{S}$ of) weakly non-factor subvarieties in $C_{d}$ is an algebraic subset of dimension zero, hence finite as $C_{d}$ has finitely many components.

\subsection{Prior Work}

The closest work in spirit to ours appears to be a recent preprint \cite{binyamini2021effective} by Binyamini and Daw. The authors work in the Shimura variety setting, and are primarily concerned with applications to the ``geometric part'' of the Zilber-Pink conjecture. They give an effective version of the result in \cite{daw2018applications} that says that if $S = \Gamma \backslash D$ is a Shimura variety and $V \subset S$ is a subvariety, then the so-called weakly optimal subvarieties of $V$ belong to finitely many constructible algebraic families. To recover this sort of result in our language, one would equip $S$, after possibly passing to a finite covering, with a variation $\mathbb{V}$ of integral Hodge structure coming from the Shimura datum. The corresponding families of weakly optimal varieties would then be obtained from suitable families of weakly specials of $S$ by intersecting their fibres with $V$; one would however need the effective bounds on degrees of weakly optimal subvarieties in \cite{binyamini2021effective} to know which such families to consider.

\subsection{Acknowledgements}

The author thanks Jacob Tsimerman for numerous discussions related to this work.

\section{Preliminaries on Computational Models}
\label{compmodelsec}

\subsection{Elementary Operations}

The algorithms we describe in this paper use as elementary steps many of the standard ``operations'' of algebraic geometry; for instance, computing the image of a constructible set under an algebraic map, trivializing vector bundles, etc. We indicate which such elementary operations we will require here, and where the reader may find references for the existence of algorithms for performing such operations. We further note that specifying an algorithm to compute an object is the same as giving a constructive existence proof for that object, so we will often use the two notions interchangeably. 

In what follows all spaces we consider are algebraic varieties, not necessarily reduced or irreducible, over a field $K$ whose field operations are computable. An affine variety $V$ we will regard as a finitely-presented $K$-algebra $R$, and morphisms of affine varieties we model as morphisms of the associated $K$-algebras. To model a more general variety $V$, we will have a finite collection of affine varieties $\{ V_{i} \}_{i = 1}^{n}$ with $V_{i} = \Spec R_{i}$, and ideals $I_{ij} \subset R_{i} \otimes_{K} R_{j}$ representing the intersection of the diagonal $\Delta_{V} \subset V \times V$ with $V_{i} \times V_{j}$. We may regard this as giving a constructive version of the condition that $V$ be separated. A morphism $f : V \to W$ of varieties we represent as the data of:
\begin{itemize}
\item covers $\{ V_{i} \}_{i = 1}^{n}$ and $\{ W_{j} \}_{j = 1}^{m}$ such that for each $j$, the inverse image $f^{-1}(W_{j})$ may be written as a union $\bigcup_{i \in I_{j}} V_{i}$ for some subset $I_{j} \subset \{ 1, \hdots, n \}$;
\item for each $i \in I_{j}$ a map $f_{ij} : V_{i} \to W_{j}$.
\end{itemize}
Thus given two varieties $V$ and $W$, in order to represent a morphism $f : V \to W$ between them for the purposes of computation, one might first have to refine the cover of $V$. We also note that given a map $g : W \to X$ represented by maps $g_{jk}$ such that the cover of $W$ for $g$ agrees with the cover of $W$ for $f$, we may compute the composition by computing appropriate compositions $g_{jk} \circ f_{ij}$. 

Let us also note that given a model $\{ V_{i} \}_{i = 1}^{n}$ for $V$ as above, we obtain canonical models for the intersections $V_{i} \cap V_{j}$ as the quotients $T_{ij} = R_{i} \otimes_{K} R_{j} / I_{ij}$. We then regard a quasi-coherent sheaf $\mathcal{H}$ over $V$ as consisting of the following data:
\begin{itemize}
\item for each $V_{i} = \Spec R_{i}$ in the cover of $V$, a finitely-presented $R_{i}$-module $M_{i}$;
\item for all $1 \leq i, j \leq n$ a finitely-presented $T_{ij}$-module $N_{ij}$, along with a restriction morphism $\textrm{res}_{ij} : M_{i} \to N_{ij}$ of $R_{i}$-modules;
\item for all $1 \leq i, j \leq n$ morphisms of finitely-presented $T_{ij}$-modules $\xi_{ij} : N_{ij} \xrightarrow{\sim} N_{ji}$ satisfying the axioms of a gluing datum.
\end{itemize}
We will regard the data of these modules and morphisms of modules as a model of $\mathcal{H}$ \emph{relative to the cover} $V_{i}$ for $V$. Whenever we refer to algorithms involving ``vector bundles'' we mean to work with the associated quasi-coherent sheaf. 

A closed subvariety $F \subset \Spec R$ of an affine variety we will represent by an ideal $I \subset R$, which is itself represented in terms of a finite set of generators. An open subvariety $U \subset \Spec R$ of an affine variety we will represent by finitely many elements $f_{1}, \hdots, f_{k} \in R$ such that $U = \bigcup_{i = 1}^{k} D(f_{i})$, where the $D(f_{i}) = \Spec R_{f_{i}}$ are the distinguished affine opens. A constructible set $C$ in $\Spec R$ will then be represented by a finite list $F_{1}, U_{1}, \hdots, F_{m}, U_{m}$ of closed and open subvarieties such that $C = \bigcup_{i = 1}^{m} F_{i} \cap U_{i}$. In the case of closed, open, and constructible sets inside a general variety $V = \{ V_{i} \}_{i = 1}^{n}$, we will represent these by compatible closed, open and constructible sets inside each $V_{i}$. 

\begin{prop}
\label{algbasicsprop}
Let $K$ be a field admitting a computational model, let $X, Y$ be $K$-varieties, let $f : X \to Y$ be a map of $K$-varieties, let $C \subset X$ and $D \subset Y$ be $K$-constructible sets, and let $\mathcal{V}$ and $\mathcal{W}$ be $K$-vector bundles on $Y$. Then there exists algorithms for:
\begin{itemize}
\item[(i)] given a presentation of $Y$ in terms of affine opens $\{ Y_{i} \}_{i = 1}^{n}$ and affine open covers $ \{ U_{ij} \}_{j = 1}^{m(i)}$ for $Y_{i}$, computing a presentation of $Y$ in terms of the cover $\{ U_{ij} \}$;
\item[(ii)] computing the image $f(C) \subset Y$ and the inverse image $f^{-1}(D) \subset X$;
\item[(iii)] computing a cover $\{ E_{i} \}_{i = 1}^{q}$, refining the given cover $\{ Y_{i} \}_{i = 1}^{m}$ for $Y$, and a presentation of $\mathcal{V}$ by free modules with respect to the cover $\{ E_{i} \}_{i = 1}^{q}$;
\item[(iv)] computing the pull back bundle $f^{*} \mathcal{V}$, the ``tensor product'' bundle $\mathcal{V} \otimes \mathcal{W}$, and the bundle $\sheafhom(\mathcal{V}, \mathcal{W})$.
\end{itemize}
\end{prop}

\begin{proof}
~ 
\begin{itemize}
\item[(i)] The statement reduces to the following problem: given inclusions of affine varieties $\iota_{k} : U_{k} \hookrightarrow Y_{k}$ and $\iota_{\ell} : U_{\ell} \hookrightarrow Y_{\ell}$, represented as maps of finitely-generated algebras $g_{k} : R_{k} \to R'_{k}$ and $g_{\ell} : R_{\ell} \to R'_{\ell}$, compute a generating set for the ideal $I'_{k\ell} \subset R'_{k} \otimes_{K} R'_{\ell}$ corresponding to the inverse image of the diagonal under the map $\iota_{k} \times \iota_{\ell}$. We claim it suffices to take as generators for $I'_{k\ell}$ the images of the generators of $I_{k\ell} \subset R_{k} \otimes_{K} R_{\ell}$ under the map $g_{k} \otimes g_{\ell}$. Indeed, this is easily checked from the fact that if $g : R \to R'$ and $q : R \to R/I$ are ring maps, then the pushout $R' \otimes_{R} R/I$ is represented by $R'/I'$, where $I'$ is the inflation of $I$ in $R'$.

\item[(ii)] In both situations we may reduce to the affine case where $X = \Spec B$, $Y = \Spec A$ and $f$ is given by a map $a : A \to B$. For the case of the computation of $f(C)$ an algorithm can be found for instance in \cite{barakat2019algorithmic}, among other places. In the case of $f^{-1}(D)$ it suffices to consider the cases where $D$ is either open or closed. In the closed case $D$ is represented by an ideal $I \subset A$, and the ideal defining $f^{-1}(D)$ is the inflation $B \cdot a(I)$ for which generators are easily computed from a generating set of $I$. In the case where $D$ is open it is represented by finitely many elements $f_{1}, \hdots, f_{k} \in A$ such that $D = \bigcup_{i = 1}^{k} \Spec A_{f_{i}}$, and then $f^{-1}(D)$ is represented by the elements $a(f_{1}), \hdots, a(f_{k})$.

\item[(iii)] Let us first deal with the case where $Y = \Spec R$ is affine, and $\mathcal{V}$ is represented by a projective $R$-module $M$. Then the claim then reduces to finding finite sets of elements $f_{1}, \hdots, f_{k} \in R$ generating the unit ideal such that for each $j$, the module $M_{f_{j}}$ is free over the localization $R_{f_{j}}$, and computing a basis for $M_{f_{j}}$. This we can do using the constructive proof that any projective module over a Noetherian ring is locally free found in \cite[Chap. X, \S1]{constrcommalg}. 

\hspace{1em} Using (i) to refine the cover of $Y$ as necessary, the general case then reduces to the following situation. We have two affine open sets $\Spec R \subset Y$ and $\Spec R' \subset Y$, modules $M$ and $M'$ representing $\mathcal{V}$ on $\Spec R$ and $\Spec R'$, restriction maps $\textrm{res} : M \to N$ and $\textrm{res}' : M' \to N'$, a gluing map $\xi : N \to N'$, and two localizations $M_{f}$ and $M'_{f'}$ of $M$ and $M'$. We then have to compute the associated restriction morphisms $\textrm{res}_{f} : M_{f} \to N_{f}$ and $\textrm{res}_{f'} : M'_{f'} \to N'_{f'}$ as well as the gluing morphism $\xi_{f} : N_{f} \xrightarrow{\sim} N'_{f'}$. As the maps $\textrm{res}, \textrm{res}'$ and $\xi$ are given in terms of the finite presentations for the modules $M, M', N$ and $N'$, the maps $\textrm{res}_{f}, \textrm{res}_{f'}$ and $\xi_{f}$ may be simply computed by extending scalars along the maps $R \to R_{f}$ and $R' \to R'_{f'}$. 

\item[(iv)] By (iii) we may refine the cover $\{ Y_{i} \}_{i = 1}^{m}$ of $Y$ to a cover $\{ Y'_{j} \}_{j = 1}^{m'}$ such that $\mathcal{V}$ and $\mathcal{W}$ are presented using free modules with explicit bases with respect to $\{ Y'_{j} \}_{j = 1}^{m}$. We note that after refining the cover $\{ Y_{i} \}_{i = 1}^{m}$, we may use (ii) with $D = Y'_{j}$ for $j = 1, \hdots, m'$ and part (i) to compute a cover $\{ X'_{k} \}_{k = 1}^{n'}$ of $X$ such that the morphism $f$ may be presented with respect to $\{ X'_{k} \}_{k = 1}^{n'}$ and $\{ Y'_{j} \}_{j = 1}^{m}$.

\hspace{1em} We begin with $f^{*} \mathcal{V}$. Let $\{ Y_{i} \}_{i = 1}^{m}$ and $\{ X_{k} \}_{k=1}^{n}$ be the affine covers of $X$ and $Y$, denote the rings associated to $Y_{i}$ by $R_{i}$, and the rings associated to the intersections $Y_{i} \cap Y_{j}$ by $T_{ij}$. The data of $\mathcal{V}$ is then a collection of free $R_{i}$-modules $M_{i}$, free $T_{ij}$-modules $N_{ij}$ with restrictions $\textrm{res}_{ij} : M_{i} \to N_{ij}$, and gluing data $\xi_{ij} : N_{ij} \xrightarrow{\sim} N_{ji}$. The cover $\{ X_{k} \}_{k=1}^{n}$ corresponds to rings $S_{k}$ representing the $X_{k}$, and we may compute rings $U_{k\ell}$ representing $X_{k} \cap X_{\ell}$ using the diagonal ideals $I_{k\ell}$. Our assumption on the cover for $Y$ as it relates to the computational model for $f$ ensures that we may associate to each $k$ an index $f(k)$ such that a morphism $\alpha_{k} : R_{f(k)} \to S_{k}$ is part of the data for $f$. Given two such morphisms $\alpha_{k} : R_{f(k)} \to S_{k}$ and $\alpha_{\ell} : R_{f(\ell)} \to S_{\ell}$, if we denote by $I_{f(k)f(\ell)} \subset R_{f(k)} \otimes_{K} R_{f(\ell)}$ and $I_{k\ell} \subset S_{k} \otimes_{K} S_{\ell}$ the diagonal ideals, we then obtain for each pair $(k, \ell)$ a canonical morphism 
\[ \beta_{k\ell} : T_{f(k)f(\ell)} = (R_{f(k)} \otimes_{K} R_{f(\ell)})/I_{f(k)f(\ell)} \to (S_{k} \otimes_{K} S_{\ell})/I_{k\ell} = U_{k\ell} , \]
which is compatible with the maps $\alpha_{k}, \alpha_{\ell}$, and the natural restrictions.

\hspace{1em} We now construct $f^{*} \mathcal{V}$. We compute $M'_{k} = M_{f(k)} \otimes_{R_{f(k)}} S_{k}$ as the free module obtained from $M_{f(k)}$ by extending scalars, and similarly for the modules $N'_{k\ell} = N_{f(k)f(\ell)} \otimes_{T_{f(k)f(\ell)}} U_{k \ell}$. The restrictions $\textrm{res}'_{k\ell} : M'_{k} \to N'_{k\ell}$ are likewise given by scalar extension of $\res_{f(k)f(\ell)}$ using the maps $\alpha_{k}$ and $\beta_{k\ell}$, and the gluing datum $\xi'_{k\ell} : N'_{k \ell} \xrightarrow{\sim} N'_{\ell k}$ is computed via scalar extension of $\xi_{f(k)f(\ell)}$ along the map $\beta_{k\ell}$.

\hspace{1em} For the tensor product $\mathcal{V} \otimes \mathcal{W}$ we keep the notation from the preceding paragraph, and suppose that $\mathcal{W}$ is represented by free $R_{i}$-modules $P_{i}$. Then one easily computes a representation for $\mathcal{V} \otimes \mathcal{W}$ using the modules $M_{i} \otimes_{R_{i}} P_{i}$, using the obvious construction for the gluing data.

\hspace{1em} For the construction of $\sheafhom(\mathcal{V}, \mathcal{W})$ we may reduce to the case where $\mathcal{W} = \mathcal{O}_{Y}$, using the fact that $\sheafhom(\mathcal{V}, \mathcal{W}) \simeq \mathcal{W} \otimes \sheafhom(\mathcal{V}, \mathcal{O}_{Y})$. As each of the constituent modules $M_{i}$ are free with an explicit basis, we may take the dual basis to construct $\textrm{Hom}(M_{i}, R_{i})$. Gluing the modules $\textrm{Hom}(M_{i}, R_{i})$ is then easily done using the gluing data for $\mathcal{V}$.
\end{itemize}
\end{proof}

\subsection{Representing Variations Algebraically}

\subsubsection{Passing to Algebraic Data}

A variation of integral Hodge structure $(\mathbb{V}, F^{\bullet})$ on an smooth algebraic variety $S$ consists of an integral local system $\mathbb{V}$ on $\an{S}$ and a filtration $F^{\bullet}$ on the holomorphic vector bundle $\an{\mathcal{H}} = \mathbb{V} \otimes_{\mathbb{Z}} \mathcal{O}_{\an{S}}$. In general, this data is highly transcendental, owing to the fact that the image in the fibre $\an{\mathcal{H}}_{s}$ of the lattice $\mathbb{V}_{s}$ is a transcendental object with respect to any coordinate system on $\an{\mathcal{H}}$ coming from algebraic geometry. For this reason one cannot hope to represent the data $(\mathbb{V}, F^{\bullet})$ in a finitary manner for the purposes of computation, and in our algorithms we work entirely with the coarser data $(\mathbb{V}_{\mathbb{C}}, F^{\bullet})$, where $\mathbb{V}_{\mathbb{C}}$ is the complex local system associated to $\mathbb{V}$.

The advantage of passing to the coarser data $(\mathbb{V}_{\mathbb{C}}, F^{\bullet})$ is that the complex local system $\mathbb{V}_{\mathbb{C}}$ is equivalent, by the Riemann-Hilbert correspondence, to the data of an algebraic vector bundle $\mathcal{H}$ with a flat algebraic connection $\nabla : \mathcal{H} \to \Omega^{1}_{S} \otimes \mathcal{H}$, from which $\mathbb{V}_{\mathbb{C}}$ is recovered by taking flat sections of $\mathcal{H}$ with respect to $\nabla$. The filtration $F^{\bullet}$ on $\an{\mathcal{H}}$ is in fact algebraic, see \cite[pg. 235]{Schmid1973}.
\begin{defn}
Given a variation $(\mathbb{V}, F^{\bullet})$, we call the tuple $(\mathcal{H}, F^{\bullet}, \nabla)$ the \emph{associated algebraic data} of the variation. If each element of the tuple is defined over a subfield $K \subset \mathbb{C}$, we will also say the \emph{associated $K$-algebraic data}.
\end{defn}
The algorithms we describe will work entirely with the associated algebraic data $(\mathcal{H}, F^{\bullet}, \nabla)$. (Of course, for the purposes of proving the correctness of our computations, we assume that this data is associated to a polarizable integral variation of Hodge structure.) We therefore devote the following two subsections to explaining how to obtain computational models for the data $(\mathcal{H}, F^{\bullet}, \nabla)$ in the situation where the variation $\mathbb{V}$ comes from geometry, which is the primary situation of interest.

\subsubsection{Models in the Geometric Setting}

Suppose that $f : X \to S$ is a smooth projective\footnote{Here we mean in the sense of \cite{hartshorne}, i.e., there exists an embedding $X \hookrightarrow \mathbb{P}^{n}_{S}$ of $S$-schemes for some $n$.} morphism of $K$-varieties, with $S$ smooth, connected, and quasi-projective. Suppose that $(\mathbb{V}, F^{\textrm{an}, \bullet})$ is the associated variation of integral Hodge structure, with $\mathbb{V} = R^{k} f_{*} \mathbb{Z}(k)$ and $F^{\textrm{an}, \bullet}$ the Hodge filtration. There is a standard way of obtaining a $K$-algebraic model $(\mathcal{H}, F^{\bullet}, \nabla)$ for the data $(\mathbb{V}_{\mathbb{C}}, F^{\textrm{an}, \bullet})$, which we now describe. 

We take $\mathcal{H} = R^{k} f_{*} \Omega^{\bullet}_{X/S}$, where $\Omega^{\bullet}_{X/S}$ is the relative algebraic de Rham complex. To give the filtration on $\mathcal{H}$, we define the Hodge filtration $F^{i} \Omega^{\bullet}_{X/S} = \Omega^{\bullet \geq i}_{X/S}$, from which the filtration on $\mathcal{H}$ is obtained as the image
\[ F^{i} \mathcal{H} = \textrm{image} \, \left[ R^{k} f_{*} F^{i} \Omega^{\bullet}_{X/S} \to R^{k} f_{*} \Omega^{\bullet}_{X/S} \right] . \]
To construct the connection, we follow \cite{katz1968}. Define a filtration $L^{\bullet}$ on the complex $\Omega^{\bullet}_{X}$ via
\[ L^{i} \Omega^{\bullet}_{X} = \textrm{image} \left[ \Omega^{\bullet - i}_{X} \otimes_{\mathcal{O}_{X}} f^{*} (\Omega^{i}_{S}) \to \Omega^{\bullet}_{X} \right] . \]
We have an exact sequence of complexes
\begin{equation}
\label{exactseq}
0 \to f^{*} \Omega^{1}_{S} \otimes \Omega^{\bullet-1}_{X/S} \to \Omega^{\bullet}_{X} / L^{2} \Omega^{\bullet}_{X} \to \Omega^{\bullet}_{X/S} \to 0 .
\end{equation}
Then we have
\begin{thm}[Katz-Oda]
\label{GMconn}
The connecting homomorphism
\[ \mathcal{H} = R^{k} f_{*} \Omega^{\bullet}_{X/S} \xrightarrow{\delta} \Omega^{1}_{S} \otimes R^{k+1} f_{*} \Omega^{\bullet - 1}_{X/S} = \Omega^{1}_{S} \otimes \mathcal{H} , \]
is the Gauss-Manin connection $\nabla$.
\end{thm}

\medskip

The algorithms we describe in this paper will actually require something slightly weaker than the associated algebraic data $(\mathcal{H}, F^{\bullet}, \nabla)$; it will suffice to have an open cover $\{ S_{i} \}_{i = 1}^{n}$ of $S$ and models for the restricted data $\restr{(\mathcal{H}, F^{\bullet}, \nabla)}{S_{i}}$ for each $i$. (We will not require the gluing data.) For this reason, the rest of this section is devoted to the proof of the following:

\begin{thm}
\label{modelscanbecomputed}
Given a smooth projective $K$-algebraic family $f : X \to S$ over a smooth base $S$, consider the variation of Hodge structure with underlying local system $\mathbb{V} = R^{k} f_{*} \mathbb{Z}(k)$. Then there exists an algorithm for computing:
\begin{itemize}
\item[(i)] an affine open cover $\{ S_{i} \}_{i = 1}^{n}$ of $S$ with $S_{i} = \Spec R_{i}$;
\item[(ii)] for each $i$, the $K$-algebraic data $(\mathcal{H}_{i}, F^{\bullet}_{i}, \nabla_{i})$ associated to the map $f_{i} : X_{S_{i}} \to S_{i}$ obtained by base-change.
\end{itemize}
\end{thm}

Our approach will be to adapt a method of \cite{scheiblechner2012castelnuovo} to the relative setting in order to compute $\mathcal{H}$ and $F^{\bullet}$ over a suitable affine base, and then to describe how to compute the connecting homomorphism in \autoref{GMconn}. As necessary setup, we first recall how to use a \v{C}ech resolution to describe the cohomology sheaves of the complexes appearing in the exact sequence (\ref{exactseq}) in the case where $S$ is affine, following \cite{katz1968}. 

Fix a finite cover $\mathcal{U} = \{ U_{i} \}_{i = 1}^{k}$ of $X$. Given a complex $\mathcal{L}^{\bullet}$ of quasi-coherent sheaves on $X$, we define a double complex $C^{\bullet}(\mathcal{U}, \mathcal{L}^{\bullet})$, where $C^{q}(\mathcal{U}, \mathcal{L}^{p})$ is the set of $q$-cochains with values in $\mathcal{L}^{p}$; an element $\beta \in C^{q}(\mathcal{U}, \mathcal{L}^{p})$ is a map which assigns to a set of indices $i_{0} < i_{1} < \cdots < i_{q}$ an element of $\mathcal{L}^{p}(\mathcal{U}_{i_{0}} \cap \cdots \cap \mathcal{U}_{i_{q}})$, where $1 \leq i_{j} \leq k$ for $0 \leq j \leq q$. We then have two differentials, $d : C^{q}(\mathcal{U}, \mathcal{L}^{p}) \to C^{q}(\mathcal{U}, \mathcal{L}^{p+1})$ defined by
\[ (d\beta)(i_{0}, \hdots, i_{q}) = d(\beta(i_{0}, \hdots, i_{q})) , \]
and a differential $\delta : C^{q} (\mathcal{U}, \mathcal{L}^{p}) \to C^{q+1}(\mathcal{U}, \mathcal{L}^{p})$ defined by
\[ (\delta \beta)(i_{0}, \hdots, i_{q+1}) = (-1)^{p} \sum_{j = 0}^{q+1} (-1)^{j} \beta(i_{0}, \hdots, \widehat{i_{j}}, \hdots, i_{q+1}) . \]
With these definitions, the double complex $K^{\bullet}(\mathcal{L}^{\bullet})$ with $K^{n}(\mathcal{L}^{\bullet}) = \bigoplus_{p+q=n} C^{q}(\mathcal{U}, \mathcal{L}^{p})$ and differential $d+\delta$ computes the cohomology of $\mathcal{L}^{\bullet}$ for each $\mathcal{L}^{\bullet}$ appearing in the exact sequence (\ref{exactseq}).

Although the \v{C}ech complex is given by explicit generators, simply constructing the \v{C}ech complex is not enough to compute the cohomology sheaves $R^{k} f_{*} \Omega^{\bullet}_{X/S}$, even if $S = \Spec R$ is affine. The issue is that although all the usual abelian category constructions, including kernels, images and quotients, may be carried out constructively in the category of finitely-presented $R$-modules, the modules that appear in $K^{\bullet}(\Omega^{\bullet}_{X/S})$ are in general infinite dimensional over $R$. Thus our goal for computing $\mathcal{H}$ (and also $F^{\bullet}$) will be to construct a finite-dimensional subcomplex of $K^{\bullet}(\Omega^{\bullet}_{X/S})$ which computes the cohomology sheaves, for which we adapt the method of \cite{scheiblechner2012castelnuovo} to the relative setting. This requires the following:

\begin{lem}
\label{hyperplaneclaim}
Let $d$ be the relative dimension of the family $f$, and fix an embedding $\iota : X \hookrightarrow \mathbb{P}^{m}_{S}$ of $S$-schemes for some $m$. Then there exists an algorithm to compute:
\begin{itemize}
\item a model for $S$ in terms of affine open subsets $\{ S_{i} \}_{i = 1}^{k}$; and
\item relative hyperplanes $h_{ij} : H_{ij} \to S_{i}$ for $j = 0, \hdots, d$ and $i = 1, \hdots, k$, with $H_{ij} \subset \mathbb{P}^{m}_{S}$, such that for each $s \in S_{i}$ the divisor $X_{s} \cap (H_{i0,s} \cup \cdots \cup H_{id,s})$ has normal crossings in $X_{s}$, and $H_{ij}$ is of the form $H \times S \subset \mathbb{P}^{m} \times S = \mathbb{P}^{m}_{S}$ for some $K$-algebraic hyperplane $H \subset \mathbb{P}^{m}$.
\end{itemize}
\end{lem}

\begin{proof}
Applying \autoref{algbasicsprop}(i) it suffices to assume that $S$ is affine. As $K$ is a field admitting a computational model, the collection of possible $K$-algebraic relative hyperplanes $h : H \to S$ of the prescribed type may be constructively enumerated, as can $(d+1)$-tuples $\overline{h} = (h_{0}, \hdots, h_{d})$ of such hyperplanes. Given such a $d$-tuple with $h_{j} : H_{j} \to S$, the condition that $X_{s} \cap (H_{0,s} \cup \cdots \cup H_{d,s})$ have normal crossings in $X_{s}$ is an open algebraic condition on $S$, and an open subset $S_{\overline{h}} \subset S$ lying in the locus where it holds can be computed by adapting the procedure in \cite[\S7]{scheiblechner2012castelnuovo}; namely, for each $\Spec R \subset S$ in an affine cover of $S$ we compute the ideal $I$ defining the intersection appearing in \cite[Eq. 21]{scheiblechner2012castelnuovo} over a $R$ instead of over $\mathbb{C}$, and then the locus $S_{\overline{h}} \cap \Spec R$ will be the locus where the coefficients of the generators of $I$ vanish. If we start computing the loci $S_{\overline{h}}$ for all possible choices of $\overline{h}$ in sequence, one has at the $i$'th stage affine open subsets $S_{\overline{h}_{1}}, \hdots, S_{\overline{h}_{i}}$ each with an appropriate relative family of hyperplanes, and it suffices to show that we have $S = S_{\overline{h}_{1}} \cup \cdots \cup S_{\overline{h}_{i}}$ for some $i$. This reduces by \cite[Lem. 7.1]{scheiblechner2012castelnuovo} to showing that for each proper complex algebraic subvariety $X \subset \mathbb{P}^{m}$, there exists a set of $K$-algebraic hyperplanes $H_{0}, \hdots, H_{d} \subset \mathbb{P}^{m}$ such that for each subset $0 \leq i_{0} \leq \cdots \leq i_{q} \leq d$ of the indices the intersection $X \cap (H_{i_{0}} \cup \cdots \cup H_{i_{q}})$ is everywhere transverse. The moduli space parametrizing tuples of $(d+1)$ hyperplanes in $\mathbb{P}^{m}$ is $\mathbb{Q}$-algebraic and unirational, and the condition on $(H_{0}, \hdots, H_{d})$ is an open algebraic condition, so the result follows as any complex algebraic open set in a unirational variety over $\mathbb{Q}$ will contain a $K$-point for any subfield $K$ of $\mathbb{C}$.
\end{proof}

\medskip

Using \autoref{hyperplaneclaim}, we now see that to prove \autoref{modelscanbecomputed} it suffices to assume that $S$ is affine and we have relative hyperplanes $h_{j} : H_{j} \to S$ for $j = 0, \hdots, d$ of the prescribed type. Following \cite[\S5]{scheiblechner2012castelnuovo} we set $\mathcal{U} = \{ U_{i} \}_{i = 0}^{d}$, where $U_{i} = X \setminus (X \cap H_{i})$; we note that $U_{i}$ is an affine subset of $X$. The complex $C^{\bullet}(\mathcal{U}, \Omega^{\bullet}_{X/S})$ then admits, for each integer $t$, a natural subcomplex $D^{\bullet, \bullet}$ given by 
\[ D^{p,q} = \bigoplus_{i_{0} < \cdots < i_{q}} \Gamma(X, \Omega^{p}_{X/S}((t + p) (H_{i_{0}} \cup \cdots \cup H_{i_{q}}))) , \]
where we understand $\Omega^{p}_{X/S}(V)$ for a relative divisor $V$ to be the tensor product with the relative sheaf $\mathcal{O}_{X/S}(V)$. We then have the following generalization of the result appearing in \cite[Lem. 5.2]{scheiblechner2012castelnuovo}:

\begin{lem}
\label{subcompworks}
Suppose that $t \geq d (ed + 1) k$ where $\deg X \leq k$ and $\textrm{codim}_{\mathbb{P}^{m}_{S}} X \leq e$. Then the total complex associated to the subcomplex $D^{\bullet,\bullet}$ computes the algebraic de Rham cohomology sheaves $R^{k} f_{*} \Omega^{\bullet}_{X/S}$, and the total complex associated to the subcomplex $D^{\bullet \geq i, \bullet}$ computes the cohomology sheaves of $R^{k} f_{*} F^{i} \Omega^{\bullet}_{X/S}$.
\end{lem}

\begin{proof}
We first handle the case of $R^{\bullet} f_{*} \Omega^{\bullet}_{X/S}$. This is proven in the non-relative setting by Scheiblechner, so it suffices to reduce to that case.  Let us denote by $T^{\bullet}(\Omega^{\bullet}_{X/S})$ the total complex associated to $D^{\bullet,\bullet}$, so that $T^{\bullet}(\Omega^{\bullet}_{X/S}) \subset K^{\bullet}(\Omega^{\bullet}_{X/S})$. We then have for each $i$ a natural map $\tau : h^{i}(T^{\bullet}(\Omega^{\bullet}_{X/S})) \to h^{i}(K^{\bullet}(\Omega^{\bullet}_{X/S}))$ of $R$-modules, where $S = \Spec R$. Given any complex point $s : \Spec \mathbb{C} \to S$, we obtain a map $\tau_{s} : h^{i}(T^{\bullet}(\Omega^{\bullet}_{X_{s}/\mathbb{C}})) \to h^{i}(K^{\bullet}(\Omega^{\bullet}_{X_{s}/\mathbb{C}}))$ on the level of fibres, and the result of \cite[Lem. 5.2]{scheiblechner2012castelnuovo} then tells us that $\tau_{s}$ is an isomorphism. The $R$-module $h^{i}(K^{\bullet}(\Omega^{\bullet}_{X/S}))$ is locally-free, hence has constant fibre dimension. This then implies that the $R$-module $h^{i}(T^{\bullet}(\Omega^{\bullet}_{X/S}))$ also has constant fibre dimension. Since $R$ is assumed reduced, we find by \cite[\S II., Ex. 5.8]{hartshorne} that $h^{i}(T^{\bullet}(\Omega^{\bullet}_{X/S}))$ is locally free. A map of locally free modules which is an isomorphism on fibres is an isomorphism, hence the result.

In the filtered case, we may argue as above to reduce to showing that $D^{\bullet \geq i, \bullet}$ computes $R^{k} f_{*} F^{i} \Omega^{\bullet}_{X/S}$ in the case where $S = \Spec \mathbb{C}$. For this we inspect the proof of \cite[Lem. 5.2]{scheiblechner2012castelnuovo}: he considers the natural spectral sequence $E_{r}(D^{\bullet,\bullet}) \to E_{r}(C^{\bullet}(\mathcal{U}, \Omega^{\bullet}_{X/S}))$ associated to the Hodge filtration (i.e., where the first page is obtained by taking cohomology with respect to the differential $d$ and not $\delta$), and shows that it is an isomorphism at the first stage. Necessarily this implies (see the proof of \cite[5.2.12]{weibel}) that the map $\textrm{Tot}(D^{\bullet,\bullet}) \to \textrm{Tot}(C^{\bullet}(\mathcal{U}, \Omega^{\bullet}_{X/S}))$ of total complexes induces an isomorphism between each filtered piece in each cohomological degree, and so in particular the map $\textrm{Tot}(D^{\bullet \geq i, \bullet}) \to \textrm{Tot}(C^{\bullet}(\mathcal{U}, \Omega^{\bullet \geq i}_{X/S}))$ is a quasi-isomorphism for each $i$.
\end{proof}

\medskip

To complete the proof of \autoref{modelscanbecomputed}, it now suffices to explain how to compute the connection $\nabla$ for those modules $\mathcal{H}$ arising from the subcomplex $D^{\bullet, \bullet}$. Note that to ``compute'' $\nabla$ it suffices to show we can evaluate the map $\nabla : \mathcal{H} \to \Omega^{1}_{S} \otimes \mathcal{H}$ (that is, the connecting homomorphism in \autoref{GMconn}) on a vector $v \in \mathcal{H}$, since we can then evaluate $\nabla$ on a generating set for $\mathcal{H}$ and use $R$-linearity. To complete the proof, we therefore need to show two things:
\begin{itemize}
\item[(A)] consider the derived version
\begin{equation}
\label{derexactseq}
 0 \to K^{\bullet}(f^{*} \Omega^{1}_{S} \otimes \Omega^{\bullet-1}_{X/S}) \xrightarrow{\rho} K^{\bullet}(\Omega^{\bullet}_{X} / L^{2} \Omega^{\bullet}_{X}) \xrightarrow{\sigma} K^{\bullet}(\Omega^{\bullet}_{X/S}) \to 0 , 
\end{equation}
of the exact sequence (\ref{exactseq}). Then given an element $\beta \in K^{\bullet}(\Omega^{\bullet}_{X/S})$ which lies in the kernel of the differential, it is possible to compute a lift $\widetilde{\beta}$ such that $\beta = \sigma(\widetilde{\beta})$, and compute $\rho^{-1}(d\widetilde{\beta})$;
\item[(B)] with the above setup, if $\beta \in T^{\bullet}(\Omega^{\bullet}_{X/S}) = \textrm{Tot}(D^{\bullet,\bullet})$ we have in addition that $\rho^{-1}(d \widetilde{\beta})$ lies in the total complex associated to 
\[ f^{*} \Omega^{1}_{S} \otimes D^{\bullet, \bullet} \subset C^{\bullet}(\mathcal{U}, f^{*} \Omega^{1}_{S} \otimes \Omega^{\bullet}_{X/S}) . \]
\end{itemize}

To carry out Step (A) one could try to proceed in a manner analogous to \cite{katz1968}, who construct an explicit map of complexes $K^{\bullet}(\Omega^{\bullet}_{X/S}) \to K^{\bullet}(f^{*} \Omega^{1}_{S} \otimes \Omega^{\bullet}_{X/S})$ which induces the connecting homomorphism on cohomology. However because we are merely interested in the existence of an algorithm and not its efficiency, we can make do with the following brutal method:
\begin{lem}
\label{dumblem}
Let $\mathcal{A}$ and $\mathcal{B}$ be any two computable sets, let $\kappa : \mathcal{A} \to \mathcal{B}$ be a map between them, and supposed we are given $b \in \textrm{im}(\kappa)$. Then there exists an algorithm to compute a preimage $a \in \mathcal{A}$ of $b$.
\end{lem}
\begin{proof}
Enumerate all the elements of $\mathcal{A}$ and try them one by one until we succeed.
\end{proof}

To complete Step (A), it therefore suffices to show:
\begin{lem}
There exists a computational representation of the complexes appearing in (\ref{derexactseq}) such that the maps $\rho$ and $\sigma$ as well as the inclusions
\[ \textrm{Tot}(D^{\bullet,\bullet}) \hookrightarrow K^{\bullet}(\Omega^{\bullet}_{X/S}) \hspace{2em} \textrm{and} \hspace{2em} \textrm{Tot}(f^{*} \Omega^{1}_{S} \otimes D^{\bullet,\bullet}) \hookrightarrow K^{\bullet}(f^{*} \Omega^{1}_{S} \otimes \Omega^{\bullet}_{X/S}) \]
are all computable.
\end{lem}

\begin{proof}
It suffices to show that the required statements are true on the level of the associated double complexes. This, in turn, reduces to showing that the vector bundles in the exact sequence (\ref{exactseq}) and the associated morphisms are computable. This follows easily from the general results on vector bundle computations shown in \autoref{algbasicsprop} above, as well as the fact that all the usual abelian category constructions may be carried out constructively in any category of finitely presented modules over a finitely presented $K$-algebra, see \cite[\S IV, 4]{constrcommalg}.
\end{proof}

Step (B), and hence the proof, is then completed by the following:

\begin{lem}
If $\beta \in T^{\bullet}(\Omega^{\bullet}_{X/S})$, then the element $\rho^{-1}(d\widetilde{\beta}) \in K^{\bullet}(\Omega^{1}_{S} \otimes \Omega^{\bullet - 1}_{X/S})$ lies inside the total complex associated to the subcomplex 
\[ f^{*} \Omega^{1}_{S} \otimes D^{\bullet,\bullet} \subset C^{\bullet}(\mathcal{U}, f^{*} \Omega^{1}_{S} \otimes \Omega^{\bullet}_{X/S}) . \]
\end{lem}

\begin{proof}
The element $\beta$ is a collection $\beta^{p,q}$ of \v{C}ech cocycles, where $\beta^{p,q}$ associates to $(i_{0}, \hdots, i_{q})$ a relative meromorphic form on $X/S$ with at worst relative poles of order $t+p$ along the relative divisor $H_{i_{0}} \cup \cdots \cup H_{i_{q}}$. The lift $\widetilde{\beta}$ is the same data, but where each component is instead considered as a form on $X$, and modulo the equivalence relation coming from $L^{\bullet}$. If we apply the differential, the poles of $d \widetilde{\beta}$ will have order at worst $t+p+1$ along $H_{i_{0}} \cup \cdots \cup H_{i_{q}}$, and this will furthermore continue to be true for the inverse image $\rho^{-1}(d \widetilde{\beta})$. But this exactly means that $\rho^{-1}(d\widetilde{\beta})$ in an element of the total complex associated to $f^{*} \Omega^{1}_{S} \otimes D^{\bullet,\bullet}$.
\end{proof}

\subsubsection{The Hyperelliptic Case} 

Let us illustrate how the computation of the associated algebraic data $(\mathcal{H}, F^{\bullet}, \nabla)$ can be done in practice with an example, unrelated to the rest of the paper. Let $S = \mathbb{A}^{2g+1} = \Spec \mathbb{Q}[e_{1}, \hdots, e_{2g+1}]$ be an affine space of dimension $2g+1$. We consider the relative hyperelliptic curve $C \subset \mathbb{P}^{2}_{S}$ which on the relative affine open subset $\mathbb{A}^{2}_{S} \subset \mathbb{P}^{2}_{S}$ with coordinates $x, y$ is given by\footnote{The factor of $4$ which appears here is included for compatibility with [Enolski and Richter], whose calculation we will refer to shortly.}
\[ y^{2} = R(x) = 4 \prod_{i = 1}^{2g+1}(x - e_{i}) = \sum_{i = 0}^{2g+1} \lambda_{i} x^{i} . \]
We write $C^{\circ}$ for the intersection $C \cap \mathbb{A}^{2}_{S}$, and $f : C \to S$ and $f^{\circ} : C^{\circ} \to S$ the families obtained by restricting the projection $\mathbb{P}^{2}_{S} \to S$. We write $\mathbb{V} = R^{1} f_{*} \mathbb{Z}(1)$ for the associated variation of Hodge structure.

We now compute the $\mathbb{Q}$-algebraic data associated to $\mathbb{V}$ by utilizing the relationship between $f$ and $f^{\circ}$. Denote by $\infty : S \to C$ the natural section so that the image of $\infty$ is the complement of $C^{\circ}$. Then we may view $f$ as a relative projective compactification of $f^{\circ}$, and $\infty$ as a complementary relative divisor which (trivially) has normal crossings. The sheaf $\mathbb{V}^{\circ} = R^{1} f^{\circ}_{*} \mathbb{Z}$ is then a local system admitting the structure of a mixed variation of Hodge structure. In fact, since $H^{1}(C_{e}, \mathbb{Z}) \simeq H^{1}(C^{\circ}_{e}, \mathbb{Z})$ canonically for each $e \in S(\mathbb{C})$, the natural base-change map
\[ \mathbb{V} = R^{1} f_{*} \mathbb{Z} \to R^{1} f^{\circ}_{*} \iota^{*} \mathbb{Z} \simeq R^{1} f^{\circ}_{*} \mathbb{Z} , \]
is an isomorphism of mixed (and hence pure) variations of Hodge structure, where $\iota : C^{\circ} \hookrightarrow C$ is the natural inclusion. We may therefore replace $\mathbb{V}$ with $\mathbb{V}^{\circ}$ and compute the $\mathbb{Q}$-algebraic data associated to $\mathbb{V}^{\circ}$.

We begin with the algebraic Hodge bundle $\mathcal{H} = R^{1} f^{\circ}_{*} \Omega^{\bullet}_{C^{\circ}/S}$. Since $C^{\circ}$ is affine over $S$, this may be computed directly from the algebraic de Rham complex $\Omega^{\bullet}_{C^{\circ}/S}$, and we therefore have $\mathcal{H} = \Omega^{1}_{C^{\circ}/S}/d\mathcal{O}_{C^{\circ}/S}$. In particular, $\mathcal{H}$ is a trivial vector bundle of rank $2g$. The filtered piece $F^{1} \an{\mathcal{H}} \subset \an{\mathcal{H}}$ corresponds to the image under the map $F^{1} (R^{1} f_{*} \mathbb{Z} \otimes \mathcal{O}_{\an{S}}) \to R^{1} f^{\circ}_{*} \mathbb{Z} \otimes \mathcal{O}_{\an{S}}$, and so consists of exactly those (families over $S$ of) forms which admit extensions to global holomorphic forms on $C$. Using the usual description of the space of holomorphic forms on hyperelliptic curves, we may find $g$ independent sections $u_{1}, \hdots, u_{g}$ of $\mathcal{H}$ of the form $u_{i} = x^{i-1} dx / y = \mathcal{U}_{i}(x) dx / y$ which together span $F^{1} \mathcal{H}$. We extend this to a basis of global sections, following \cite[\S2]{enolski}, by defining
\[ u_{g+i} =  \left(\sum_{k = i}^{2g+1 - i} (k + 1 - i) \lambda_{k+1+i} x^{k} \right) dx / 4 y = \mathcal{R}_{i}(x) dx / 4 y , \hspace{3em} 1 \leq i \leq g . \]

We now are left with the task of computing $\nabla$ in terms of the frame $u_{1}, \hdots, u_{2g}$. The tangent bundle $TS$ carries a natural trivialization by the sections $\frac{\partial}{\partial e_{1}}, \hdots, \frac{\partial}{\partial e_{2g+1}}$, and it suffices to compute the differential operators $\nabla_{\partial_{i}}$ as sections of the endomorphism bundle $\sheafend(\mathcal{H})$. Suppose that $e \in S$ is a point, and let $B \subset S$ be a ball around $e$ such that we have an integral basis $\gamma_{1}, \hdots, \gamma_{2g} \in \mathbb{V}^{\circ}(B)$. We may view these basis elements as being dual to a basis of topological cycles $c_{1}, \hdots, c_{2g}$. Let $p_{jk} = \int_{c_{k}} u_{j}$ be the periods of the $u_{j}$ with respect to these cycles, denote the associated matrix by $P$, and let $B = [b_{ij}]$ equal $P^{-1}$. We may compute that
\begin{align*}
\nabla_{\partial_{i}} u_{j} &= \nabla_{\partial_{i}}  \left( \sum_{k = 1}^{2g} p_{jk} \gamma_{k} \right) \\
&= \sum_{k = 1}^{2g} \partial_{i} p_{jk} \gamma_{k} \\
&= \sum_{k = 1}^{2g} \partial_{i} p_{jk} \left( \sum_{\ell = 1}^{2g} b_{k\ell} u_{\ell} \right) \\
\nabla_{\partial_{i}} u_{j} &=  \sum_{\ell = 1}^{2g} \left(\sum_{k = 1}^{2g} \partial_{i} p_{jk} b_{k\ell} \right) u_{\ell} .
\end{align*}
The calculation shows that in the basis given by the $u_{i}$, the operator $\nabla_{\partial_{i}}$ is equal to $(\partial_{i} P \cdot P^{-1})^{t}$. To determine this matrix we now quote \cite[Thm. 4.1]{enolski}:

\begin{prop}[Enolski-Richter]
For every $e \in U$, there exists a neighbourhood $B$ of $e$ and a basis $\gamma_{1}, \hdots, \gamma_{2g}$ for $\mathbb{V}^{\circ}(B)$ such that 
\begin{align*}
(\partial_{\ell} P \cdot P^{-1})^{t} = \begin{pmatrix} \boldsymbol{\alpha}_{\ell}^{t} & \boldsymbol{\gamma}_{\ell} \\ \boldsymbol{\beta}_{\ell} & - \boldsymbol{\alpha} \end{pmatrix}
\end{align*}
where
\begin{align*}
\boldsymbol{\alpha}_{\ell} &= \frac{-1}{2} \left( \frac{1}{R'(e_{\ell})} \boldsymbol{\mathcal{U}}(e_{\ell}) \boldsymbol{\mathcal{R}}^{t}(e_{\ell}) - M_{\ell} \right) , \\
\boldsymbol{\beta}_{\ell} &= -2 \left( \frac{1}{R'(e_{\ell})} \boldsymbol{\mathcal{U}}(e_{\ell}) \boldsymbol{\mathcal{U}}^{t}(e_{\ell}) \right) , \\
\boldsymbol{\gamma}_{\ell} &= \frac{1}{8} \left( \frac{1}{R'(e_{\ell})} \boldsymbol{\mathcal{R}}(e_{\ell}) \boldsymbol{\mathcal{R}}^{t}(e_{\ell}) - N_{\ell} \right) ,
\end{align*}
with
\[ M_{\ell} = \begin{pmatrix} 0 & 0 & 0 & \hdots & 0 & 0 \\
1 & 0 & 0 & \hdots & 0 & 0 \\
e_{\ell} & 1 & 0 & \hdots & 0 & 0 \\
e_{\ell}^{2} & e_{\ell} & 1 & \hdots & 0 & 0 \\
\vdots & \vdots & \ddots & \ddots & \ddots & \vdots \\
e_{\ell}^{g-2} & e_{\ell}^{g-3} & \hdots & e_{\ell} & 1 & 0 \end{pmatrix} , \]
and
$N_{\ell} = e_{\ell}(M_{\ell} Q_{\ell} + Q_{\ell} M_{\ell}^{t}) + Q_{\ell}$, where $Q_{\ell}$ is the diagonal matrix with $(Q_{\ell})_{k,k} = \mathcal{R}_{k}(e_{\ell})/\mathcal{U}_{k+1}(e_{\ell})$. 
\qed
\end{prop}

\section{The Constructive Jet Correspondence}
\label{jetcorrespsec}

\subsection{Statement of the Correspondence}

\subsubsection{Sketch of The Idea}

In this section we develop our main computational tool for computing sets of weakly special subvarieties. Let us explain the basic idea. Fix a set of Hodge numbers $h^{p,0}, h^{p-1,1}, \hdots, h^{0,q}$ with $h^{p,0} + h^{p-1,1} + \cdots + h^{0,q} = m$, and let $\ch{L}$ be the variety of Hodge flags on $\mathbb{Z}^{m}$, not necessarily polarized. There is a natural quotient map $q : \GL_{m} \to \ch{L}$ sending a matrix to the flag represented by its column vectors, where we define $F^{i} \mathbb{C}^{m}$ to be the span of the first $\dim F^{i}$ columns. Given a variation of Hodge structure $(\mathbb{V}, F^{\bullet})$ on $S$ with the same Hodge numbers and a simply-connected ball $B \subset S$, one may construct a period map $\psi : B \to \ch{L}$ as follows. First, pick a filtration-compatible frame $v^{1}, \hdots, v^{m}$ for $\mathcal{H} = \mathbb{V} \otimes \mathcal{O}_{\an{S}}$ on $B$, i.e., a frame such that each filtered piece $F^{i}$ is spanned by some initial segment $v^{1}, \hdots, v^{k_{i}}$. Picking a basis $b^{1}, \hdots, b^{m}$ for $\mathbb{V}(B)$, we denote by $A(s) = [a_{ij}(s)]$ the change-of-basis matrix corresponding to the bases $v^{1}_{s}, \hdots, v^{m}_{s}$ to the basis $b^{1}_{s}, \hdots, b^{m}_{s}$ at $s \in B$. The map $q \circ A : B \to \ch{L}$ is then a period map for $(\mathbb{V}, F^{\bullet})$ on $B$. More generally, we will define:

\begin{defn}
\label{locpermapdef}
If in the above construction, the basis $b^{1}, \hdots, b^{m}$ is chosen more generally to be a basis of $\mathbb{V}_{\mathbb{C}}(B)$ (instead of an integral basis of $\mathbb{V}(B)$) then we will call the resulting map $\psi = q \circ A$ is called a \emph{local period map}.
\end{defn}

The point is as follows. Period maps defined using the integral lattice are hard to model computationally, as the integral lattice $\mathbb{V}$ is a transcendental object. But we will see that collection of all \emph{local period maps} can be described entirely in terms of the associated $K$-algebraic data $(\mathcal{H}, F^{\bullet}, \nabla)$; in particular, we will find that local period maps are exactly the solutions to certain $K$-algebraic systems of differential equations determined by this data. Our goal will then be to construct the loci $\mathcal{W}_{d} \subset C_{d}$ (recall the setup of \autoref{mainthm2}) of weakly special subvarieties in three steps:
\begin{itemize}
\item[(i)] we will construct finitely many subvarieties $\ch{V}_{1}, \hdots, \ch{V}_{k} \subset \ch{L}$;
\item[(ii)] we will show that a variety $Z \subset S$ is weakly special if and only if for some $\ch{V}_{j}$, for any germ $(Z, s)$ of $Z$ and any local period map $\psi = q \circ A$ defined at $s$, the image $\psi(Z, s)$ lands inside some $\GL_{m}$-translate of $\ch{V}_{j}$, and $Z$ is a maximal geometrically irreducible variety satisfying this property;
\item[(iii)] we will explain how to compute, as a constructible $K$-algebraic subset, the locus of $[\overline{Z}] \in C_{d}$ such that $Z = \overline{Z} \cap S$ satisfies the condition described in (ii).
\end{itemize}

For the remainder of \autoref{jetcorrespsec} we will develop the machinery necessary to work with local period maps computationally so that we may carry out this strategy in the coming sections.

\subsubsection{Jets and Canonical Jet Torsors}

The construction we will describe will use \emph{jet spaces}. These are schemes associated to an algebraic variety $S$ and integers $d, r \geq 0$ that parametrize infinitesimal maps to $S$ from a $d$-dimensional disk of order $r$. We let $A^{d}_{r} = K[t_{1}, \hdots, t_{d}]/(t_{1}, \hdots, t_{d})^{r+1}$, and $\mathbb{D}^{d}_{r} = \Spec A^{d}_{r}$. We denote by $\FSch_{K}$ the category of finite-type schemes over the field $K$.

\begin{defn}
\label{jetspacedef}
Let $S$ be an object in $\FSch_{K}$. We define $J^{d}_{r} S$ to be the variety representing the contravariant functor $\FSch_{K} \to \Set$ given by 
\[ T \mapsto \Hom_{K}(T \times_{K} \mathbb{D}^{d}_{r}, S), \hspace{1.5em} [T \to T'] \mapsto [\Hom_{K}(T' \times_{K} \mathbb{D}^{d}_{r}, S) \to \Hom_{K}(T \times_{K} \mathbb{D}^{d}_{r}, S)] , \]
where the natural map $\Hom_{K}(T' \times_{K} \mathbb{D}^{d}_{r}, S) \to \Hom_{K}(T \times_{K} \mathbb{D}^{d}_{r}, S)$ obtained by pulling back along $T \times_{K} \mathbb{D}^{d}_{r} \to T' \times_{K} \mathbb{D}^{d}_{r}$. We note that all maps here are of $K$-schemes.
\end{defn}
\noindent There is an analogous construction when $S$ is a complex analytic space, compatible with analytification  whenever $K \subset \mathbb{C}$ is a subfield; for details we refer to \cite[\S2]{periodimages}. We note that $J^{d}_{r}$ is a functor for each $d, r \geq 0$. We will typically view points $j \in (J^{d}_{r} S)(T)$ as maps $T \times_{K} \mathbb{D}^{d}_{r} \to S$ from the infinitesimal disk $T \times_{K} \mathbb{D}^{d}_{r}$ over $T$, and from this point of view the map $J^{d}_{r} g : J^{d}_{r} S \to J^{d}_{r} S'$ induced by a map $g : S \to S'$ sends $j$ to $g \circ j$.

By functoriality, the action of $\GL_{m}$ on $\ch{L}$ induces an action of $\GL_{m}$ on $J^{d}_{r} \ch{L}$. The main results of this section is are then the following:

\begin{thm}
\label{constrjetcorrespthm}
Given a tuple $(\mathcal{H}, F^{\bullet}, \nabla)$ of associated $K$-algebraic data on the smooth $K$-variety $S$, and for each $d, r \geq 0$, there exists a canonical map $\eta^{d}_{r} : J^{d}_{r} S \to \GL_{m} \backslash J^{d}_{r} \ch{L}$ of $K$-algebraic stacks with the following property: given any point $j \in (J^{d}_{r} S)(\mathbb{C})$ and any local period map $\psi : B \to \ch{L}$ such that $B$ contains the image of $j : \mathbb{D}^{d}_{r} \to S$, we have $\eta^{d}_{r}(j) = \psi \circ j$ modulo $\GL_{m}$.

Moreover, there exists an algorithm to compute the $\GL_{m}$-torsor $p^{d}_{r} : \mathcal{T}^{d}_{r} \to J^{d}_{r} S$ and the $\GL_{m}$-invariant map $\alpha^{d}_{r} : \mathcal{T}^{d}_{r} \to J^{d}_{r} \ch{L}$ which defines $\eta^{d}_{r}$.
\end{thm}

\vspace{0em}

\begin{cor}
\label{compdiffconst}
Suppose that $\mathbb{V}$ is a variation of Hodge structure on a smooth $K$-variety $S = \bigcup_{i = 1}^{n} S_{i}$ with each $S_{i}$ an affine open subset, and we have computational models $(\mathcal{H}_{i}, F^{\bullet}_{i}, \nabla_{i})$ for the $K$-algebraic data associated to $\restr{\mathbb{V}}{S_{i}}$. Then given a subset $\mathcal{K}$ of the points of $\GL_{m} \backslash J^{d}_{r} \ch{L}$ corresponding to a constructible $\GL_{m}$-invariant subset $\mathcal{R} \subset J^{d}_{r} \ch{L}$, there exists an algorithm to compute the constructible $K$-algebraic subset $(\eta^{r}_{d})^{-1}(\mathcal{K}) \subset J^{d}_{r} S$.
\end{cor}

\begin{proof}
As the map $\eta^{d}_{r}$ described in \autoref{constrjetcorrespthm} is functorial with respect to restrictions of $\mathbb{V}$, we may show the required claim on each $S_{i}$ separately. The required preimage is then just $p^{d}_{r}((\alpha^{d}_{r})^{-1}(\mathcal{R}))$. 
\end{proof}

A version of \autoref{constrjetcorrespthm} without the computability claim and with the pair $(\ch{L}, \GL_{m})$ replaced by the pair $(\ch{D}, \textrm{Aut}(\mathbb{Z}^{m}, Q))$, where $\ch{D}$ is the flag variety of polarized Hodge flags on the polarized lattice $(\mathbb{Z}^{m}, Q)$, appears in \cite[Thm 1.11]{periodimages}. The difference here is that we do not assume an algebraic model for the polarization as \cite{periodimages} does, meaning the notion of local period map we work with is more general (the period isomorphisms need not preserve polarizations), and correspondingly we need to take a quotient by a larger symmetry group. Nevertheless, we will essentially follow the construction in \cite[\S3]{periodimages} with only minor modifications.

\subsection{Construction of Canonical Jet Torsors}

Let us explain essential idea underlying the construction in \autoref{constrjetcorrespthm}. Suppose $j \in J^{d}_{r} S$ is a jet projecting onto $s \in S$, and $\psi : B \to \ch{L}$ is a local period map defined on a small ball $B$ containing $s$. Then in algebraic coordinates on $\ch{L}$, the coordinates of the jet $\psi \circ j$ are algebraic functions of the algebraic coordinates of $j$ and the values at $s$ of the derivatives of $\psi$. What we will see is that the values of the derivatives of $\psi$ at $s$ are in fact algebraic functions of an initial condition for a system of differential equations, and the space of initial conditions for $\psi$ at $s$ is a $\GL_{m}$-torsor. Constructing the map $\alpha^{d}_{r} : \mathcal{T}^{d}_{r} \to J^{d}_{r} \ch{L}$ referred to in \autoref{constrjetcorrespthm} then essentially amounts to the observation that $\psi \circ j$ is an algebraic function of $j$ and the initial condition defining $\psi$.

\subsubsection{Jet Evaluation Maps}

In this section we explain that the main result of \cite[\S3.1]{periodimages} is constructive:

\begin{prop}
\label{jetevalprop}
Suppose that $U$ is a smooth $K$-variety, $c_{ij}$ are global sections of $\Omega^{1}_{U}$ for $1 \leq i, j \leq m$, and consider the system of differential equations\footnote{We regard the equality as taking of (\ref{diffeq}) as taking place inside $\Omega^{1}_{\an{U}}(B)$ for some analytic neighbourhood $B \subset U$.}
\begin{equation}
\label{diffeq}
df_{jk} = \sum_{i = 1}^{m} f_{ik} c_{ij}
\end{equation}
for $1 \leq j, k \leq m$. Regard a solution to (\ref{diffeq}) as a matrix-valued map $f = [f_{jk}]$, and suppose that for each point $s \in U$ and initial condition $f(s) = M \in \GL_{m}$ there exists a unique $f : B \to \GL_{m}$ solving (\ref{diffeq}) on some neighbourhood $B$ of $s$. Then for each $d, r \geq 0$ there exists a $K$-algebraic map $\beta : J^{d}_{r} U \times \GL_{m} \to J^{d}_{r} \GL_{m}$ with the following property: if $f : B \to \GL_{m}$ solves the differential system (\ref{diffeq}), then $\beta(\sigma, f(\pi(\sigma))) = f \circ \sigma$ for all $\sigma \in J^{d}_{r} B$ where $\pi : J^{d}_{r} U \to U$ is the natural projection.

Moreover, there exists an algorithm to compute $\beta$ given the sections $c_{ij}$.
\end{prop}

\begin{proof}
The existence of this map and the verification that it satisfies the required property is handled by \cite[Prop. 3.1]{periodimages}. Let us recall enough of the proof in \cite{periodimages} to verify that $\beta$ is actually computable.

It suffices to construct $\beta$ locally on $U$, so we may assume $U$ is affine of the form $\Spec A$. As $U$ is smooth, the module $\Omega^{1}_{U}$ is projective, but in general need not be free. However it is possible to constructively find an affine open cover $\{ U_{i} \}_{i=1}^{n}$ of $U$, with $U_{i} = \Spec A_{f_{i}}$ for $f_{i} \in A$, such that $\Omega^{1}_{U_{i}}$ is free over $A_{f_{i}}$; see \cite[Chap. X, \S1]{constrcommalg}. Using this, we may further restrict to assume that $\Omega^{1}_{U}$ is free over $A$, and choose trivializing sections $dz_{1}, \hdots, dz_{n}$ of $\Omega^{1}_{U}$ for functions $z_{i} \in A$. We then obtain a map $g: U \to \mathbb{A}^{n} = \Spec K[x_{1}, \hdots, x_{n}]$ given by $x_{i} \mapsto z_{i}$.

Let us recall the construction of $\beta$ given in the proof in \cite{periodimages}. It is defined as $\beta = r \circ \zeta$, where $\zeta : J^{d}_{r} U \times \GL_{m} \to \Spec R$ is a map to a certain affine space with coordinate ring $R$, and the map $r$ is a map $r : \Spec R \to J^{d}_{r} \mathbb{M}$, where $\mathbb{M}$ is the space of all $m \times m$ matrices, not necessarily invertible. The map $r$ is independent of the input data, so we describe it first. The jet space $J^{d}_{r} \mathbb{A}^{n}$ is an affine space, with natural coordinates coming from the coefficients of the maps $j : \mathbb{D}^{d}_{r} \to \mathbb{A}^{n}$ from the formal $d$-dimensional disk of order $r$; in particular, one can view such a map as a collection of $n$ power series in $d$ variates, each truncated at order $r$. We denote these coordinates on $J^{d}_{r} \mathbb{A}^{n}$ by $a_{p, i}$; here the index $1 \leq i \leq n$ corresponds to the component of the map $j$, and the index $p \in \mathcal{P}^{d}_{r}$, where $\mathcal{P}^{d}_{r}$ is the set of partitions of integers from $0$ to $r$ with $d$ terms, corresponds to the term in the power series of which $a_{p, i}$ is a coefficient. The constant terms $a_{\varnothing, 1}, \hdots, a_{\varnothing, n}$ may be identified with the coordinates $x_{1}, \hdots, x_{n}$ on $\mathbb{A}^{n}$. The ring $R$ is then the $\mathbb{Q}$-algebra freely generated by the formal symbols $a_{p, i}$ for $1 \leq i \leq n$ and $p \in \mathcal{P}^{d}_{r}$ and the symbols $(\partial_{i_{1}} \cdots \partial_{i_{q}} f_{jk})(a_{\varnothing, 1}, \hdots, a_{\varnothing, n})$, where $1 \leq j, k \leq m$ and $i_{1}, \hdots, i_{q}$ is a sequence of length ranging beween $0$ and $r$ with $1 \leq i_{\ell} \leq n$ for each $i_{\ell}$. 

Less formally, the affine space $\Spec R$ has as coordinates those of the space $J^{d}_{r} \mathbb{A}^{n}$, and in addition the values of all partial derivatives up to order $r$ of an undetermined matrix-valued function $f = [f_{jk}]$ defined on a neighbourhood of $\mathbb{A}^{n}$ and with values in $\mathbb{M}$. Denote by $\pi : J^{d}_{r} \mathbb{A}^{n} \to \mathbb{A}^{n}$ the natural projection. The map $r$ is then defined by the property that if $f : B \to \mathbb{M}$ is any analytic map on a neighbourhood of $B \subset \mathbb{A}^{n}$ and $\varphi_{f} : J^{d}_{r} B \to \Spec R$ is defined by
\begin{equation}
\label{ffacdef}
\sigma \mapsto (\sigma, (\partial_{i_{1}} \cdots \partial_{i_{q}} f_{jk})(\pi(\sigma))) ,
\end{equation}
then the map $J^{d}_{r} f : J^{d}_{r} B \to J^{d}_{r} \mathbb{M}$ is given by $r \circ \varphi_{f}$. Note that in (\ref{ffacdef}) it is understood that the indices $j$ and $k$ range from $1$ to $m$ and $i_{1}, \hdots, i_{q}$ ranges over all appropriate sequences, as described above. After removing all the notational bookkeeping, the basic idea is simply that for any analytic function $f$ and jet $\sigma \in J^{d}_{r} \mathbb{A}^{n}$, the coordinates of $f \circ \sigma$ are $\mathbb{Q}$-algebraic functions of the coordinates of $\sigma$ and the derivatives of $f$, as can be seen for instance by writing out $\sigma$ as a truncated power series and applying the multivariate chain rule. These observations are explained in more detail in \cite[\S2.2]{periodimages}. In particular the map $r$ is universal, independent of the input data, and is defined using general $\mathbb{Q}$-algebraic polynomials obtained from differentiation rules that are easily computed in advance of running the algorithm.

Next we consider the map $\zeta$. This map is obtained from the differential forms $c_{ij}$ given as input. First, as we have trivialized $\Omega^{1}_{A}$ by a basis $dz_{1}, \hdots, dz_{n}$, we obtain for each $(i, j)$ functions $c_{ij,\ell} \in A$ with $1 \leq i \leq n$ which are the coefficients of $c_{ij}$ in this basis. The basis $dz_{1}, \hdots, dz_{n}$ induces a dual basis $\frac{\partial}{\partial z_{1}}, \hdots, \frac{\partial}{\partial z_{n}}$ for the tangent sheaf, each element of which may be identified with a computable differential operator $\frac{\partial}{\partial z_{i}} : A \to A$. Equation (\ref{diffeq}) above then becomes
\begin{equation}
\label{diffeq2}
\frac{\partial f_{jk}}{\partial z_{\ell}} = \sum_{i} f_{ik} c_{ij,\ell} , \hspace{4em} 1 \leq j, k \leq n . 
\end{equation}
We now compute recursively a series of polynomials $\xi_{(i_{1}, \hdots, i_{q}), jk}$, where $(i_{1}, \hdots, i_{q})$ is as above, which satisfy the property that $\xi_{(i_{1}, \hdots, i_{q}), jk} \in A[f_{tu}, 1 \leq t, u \leq m]$ (i.e., they are polynomials in $f_{tu}$ with coefficients in $A$) and we have
\begin{equation}
\label{deriveq}
\frac{\partial f_{jk}}{\partial z_{i_{1}} \cdots \partial z_{i_{q}}} = \xi_{(i_{1}, \hdots, i_{q}), jk}([f_{tu}]) ,
\end{equation}
for all choices of $(i_{1}, \hdots, i_{q})$ and $1 \leq j, k \leq m$. For the case where $q = 1$ we may simply use Equation (\ref{diffeq2}). In the cases for $q > 1$ we differentiate the relevant equation in the $q-1$ case and then use (\ref{diffeq2}) to substitute the first-order derivatives for polynomials in $\{ f_{tu} : 1 \leq t, u \leq m \}$ with coefficients in $A$. 

Given all this, the map $\zeta : J^{d}_{r} U \times \GL_{m} \to \Spec R$ is defined as
\[ (\sigma, [a_{tu}]) \mapsto (g \circ \sigma, \xi_{(i_{1}, \hdots, i_{q}), jk}([a_{tu}], \pi(\sigma))) , \]
where we range over $1 \leq j, k \leq m$ and the sequences $(i_{1}, \hdots, i_{q})$ above. This gives an explicit computable presentation for the map $\zeta$, and as $\beta = r \circ \zeta$ this completes the proof.
\end{proof}

\subsubsection{The Torsor $\mathcal{T}^{d}_{r}$}
\label{torsorconssec}

Let us now construct the torsor $\mathcal{T}^{d}_{r}$ and the map $\alpha^{d}_{r} : \mathcal{T}^{d}_{r} \to J^{d}_{r} \ch{L}$. Using \autoref{algbasicsprop}(iii) we fix a finite open cover $\{ U_{i} \}_{i = 1}^{n}$ of $S$ such that each $U_{i}$ carries a filtration-compatible frame $v^{1}_{i}, \hdots, v^{m}_{i}$ for the bundle $\restr{\mathcal{H}}{U_{i}}$. We denote by $C^{ij} = [c^{ij}_{k\ell}]$ the transition functions for $\mathcal{H}$ on the intersection $U_{i} \cap U_{j}$ with respect to the frames $v^{1}_{i}, \hdots, v^{m}_{i}$ and $v^{1}_{j}, \hdots, v^{m}_{j}$. We then construct the $\GL_{m}$-torsor $\mathcal{T}^{d}_{r}$ as follows. Over the open subset $J^{d}_{r} U_{i}$ of $J^{d}_{r} S$, we model it by $J^{d}_{r} U_{i} \times \GL_{m}$ with the obvious projection. The gluing data $J^{d}_{r} (U_{i} \cap U_{j}) \times \GL_{m} \xrightarrow{\sim} J^{d}_{r} (U_{i} \cap U_{j}) \times \GL_{m}$ over the open set $J^{d}_{r} (U_{i} \cap U_{j})$ is then given by $(\sigma, A) \mapsto (\sigma, A \cdot C^{ji})$. Note that the torsor $p^{d}_{r} : \mathcal{T}^{d}_{r} \to J^{d}_{r} S$ is nothing but the base-change of the frame bundle of $\mathcal{H}$ to $J^{d}_{r} S$.

Each filtration-compatible frame $v^{1}, \hdots, v^{m}$ on a $K$-algebraic open subset $U \subset S$ induces a system of differential equations on $U$ of the form in (\ref{diffeq}) in the following way. Write $c_{ij} \in \Omega^{1}_{U}$ for the unique forms giving the equality $\nabla v^{i} = \sum_{j = 1}^{m} c_{ij} \otimes v^{j}$. Suppose we wish to find analytic functions $f_{ij}$ on an open ball $B \subset U$ such that $b^{k} = \sum_{i = 1}^{m} f_{ik} v^{i}$ is a flat frame. We note that
\begin{align*}
\nabla b^{k} &= \nabla \left(\sum_{i = 1}^{m} f_{ik} v^{i} \right) \\
&= \sum_{j = 1}^{m} df_{jk} \otimes v^{j} + \sum_{i = 1}^{m} f_{ik} \left( \sum_{j = 1}^{m} c_{ij} \otimes v^{j} \right) \\
&= \sum_{j = 1}^{m} \left( df_{jk} + \sum_{i = 1}^{m} f_{ik} c_{ij} \right) \otimes v^{j} ,
\end{align*}
thus the condition that $b^{k}$ give a flat frame defines a system of differential equations of the form in (\ref{diffeq}) after absorbing a sign into the $c_{ij}$. 

Note that we can compute such a differential system for each frame $v^{1}_{i}, \hdots, v^{m}_{i}$ with respect to the cover $\{ U_{i} \}_{i = 1}^{n}$ after possibly refining the cover so that it is compatible with our computational model for $\nabla$. We therefore obtain maps $\beta_{i} : J^{d}_{r} U_{i} \times \GL_{m} \to J^{d}_{r} \GL_{m}$ for each differential system from \autoref{jetevalprop}, and we may define $\alpha^{d}_{r} : \mathcal{T}^{d}_{r} \to J^{d}_{r} \ch{L}$ on the open subset $J^{d}_{r} U_{i} \times \GL_{m}$ as the composition $q \circ \iota \circ \beta_{i}$, where $\iota : \GL_{m} \to \GL_{m}$ is the inversion. Note that the fact that the $\beta_{i}$'s are compatible with each other follows from the fact that the condition that $b^{k}$ give a flat basis is independent of the chosen frame. The condition that the map $\alpha^{d}_{r}$ be $\GL_{m}$-invariant follows from the fact that if $f = [f_{ij}]$ solves the differential system in (\ref{diffeq}), then so does $f A$, where $A \in GL_{m}$ is any matrix. We are now ready to prove:

\begin{proof}[Proof of \autoref{constrjetcorrespthm}]
It is evident from our construction that both the map $\alpha^{d}_{r} : \mathcal{T}^{d}_{r} \to J^{d}_{r} \ch{L}$ and the projection $p^{d}_{r} : \mathcal{T}^{d}_{r} \to J^{d}_{r} S$ is computable, so it suffices to show the required compatibility with local period maps. This is a local claim, so it suffices to assume that $S = U$ and we have a filtration compatible frame $v^{i} : U \to \mathcal{H}$ and the map $\alpha^{d}_{r}$ is defined as $q \circ \iota \circ \beta$, where $\beta$ is as in \autoref{jetevalprop}. Supposing that $f = [f_{ij}]$ where $f_{ij} : B \to \GL_{m}$ are analytic functions on an open set $B \subset U$ such that $b^{k} = \sum_{i = 1}^{m} f_{ik} v^{i}$ gives a flat frame for $\mathbb{V}_{\mathbb{C}}(B)$, then $\iota \circ f$ gives the change-of-basis matrix from the frame $v^{1}, \hdots, v^{m}$ to the frame $b^{1}, \hdots, b^{m}$. By definition (see \autoref{locpermapdef}), the map $\psi = q \circ \iota \circ f$ is a local period map. We then have that for $(j, M) \in J^{d}_{r} U \times \GL_{m}$ with $f(\pi(j)) = M$, where $\pi : J^{d}_{r} U \to U$ is the usual projection, that
\[ \alpha^{d}_{r}(j, M) = q \circ \iota \circ \beta(j, M) = \psi \circ j , \]
where we have used the defining property of $\beta$ in \autoref{jetevalprop}.
\end{proof}

\subsection{Germs defined by Period Maps}
\label{germsofpermaps}

For later use, let us also extract some consequences of the proof of \autoref{jetevalprop} that we will find useful later. We consider an open subset $U \subset S$ on which $\Omega^{1}_{S}$ is trivial, and let $dz_{1}, \hdots, dz_{n}$ be a trivialization over $U$. Then the coordinates $z_{1}, \hdots, z_{n}$ induce, for each point $s \in U(\mathbb{C})$, a natural basis for the $r$'th order neighbourhoods $\widehat{\mathcal{O}}_{S_{\mathbb{C}}, s}/\mathfrak{m}_{S_{\mathbb{C}},s}^{r+1}$, where the basis elements are the monic moniomials in the functions $(z_{1} - s_{1}), \hdots, (z_{n}-s_{n})$, and $s_{i} \in \mathbb{C}$ is the value of $z_{i}$ at $s$. The analogous fact is true for a point $M \in \GL_{m}(\mathbb{C})$ with respect to the natural coordinates on $\GL_{m}$.

\begin{prop}
\label{permapgermprop}
Fix a filtration compatible frame $v^{1}, \hdots, v^{m}$ for $\mathcal{H}$ on $U$. For each $s \in U(\mathbb{C})$ and $M \in \GL_{m}(\mathbb{C})$,   consider the map $\tau^{r}_{s,M} : \widehat{\mathcal{O}}_{\GL_{m,\mathbb{C}}, M}/\mathfrak{m}^{r+1}_{\GL_{m,\mathbb{C}}, M} \to \widehat{\mathcal{O}}_{S_{\mathbb{C}}, s}/\mathfrak{m}^{r+1}_{S_{\mathbb{C}}, s}$ induced by the analytic map $f : B \to \GL_{m}$ for which $b^{k} = \sum_{i = 1}^{m} f_{ik} v^{i}$ gives a flat frame with $f(s) = M$. Then with respect to the bases induced by the coordinates $z_{1}, \hdots, z_{n}$ and the natural coordinates on $\GL_{m}$, the entries of the matrix defining the linear transformations $\tau^{r}_{s,M}$ are $K$-algebraic functions on $U \times \GL_{m}$. 
\end{prop}

\begin{proof}
The map $\tau_{s,M} : \widehat{\mathcal{O}}_{\GL_{m,\mathbb{C}}, M} \to \widehat{\mathcal{O}}_{S_{\mathbb{C}}, s}$ induced by $f$ has, when written as a power series in the variates $(z_{1} - s_{1}), \hdots, (z_{n} - s_{n})$, coefficients given by constant scalar multiples of the derivatives $\partial_{i_{1}, \hdots, i_{q}} f_{jk}$, where the notation is as in \autoref{deriveq}. In particular, as in \autoref{jetevalprop} and \autoref{deriveq}, these derivatives are algebraic functions $\xi_{(i_{1}, \hdots, i_{q}),jk}$ on $U \times \GL_{m}$. Writing $M = [M_{jk}]$ and letting $a_{jk} - M_{jk}$ be the natural set of generators for $\mathfrak{m}_{\GL_{m,\mathbb{C}},M}/\mathfrak{m}^{2}_{\GL_{m,\mathbb{C}},M}$, the map $\tau_{s,M}$ sends $a_{jk} - M_{jk}$ to the power series expansion for $f_{jk} - M_{jk}$ at $s$. Thus the required functions are then obtained from the $\xi_{(i_{1}, \hdots, i_{q}),jk}$ and by applying the multiplicativity of the map $\tau_{s,M}$. 
\end{proof}

\begin{cor}
\label{pullbackfunclocper}
Continue with the setup of \autoref{permapgermprop}, and let $g$ be a $K$-algebraic algebraic function on the space of $m \times m$ matrices. Then the coordinates of $\tau^{r}_{s,M}(g)$ are $K$-algebraic functions on $U \times \GL_{m}$.
\end{cor}

\begin{proof}
The coordinates of the image of $g$ in the rings $\widehat{\mathcal{O}}_{\GL_{m,\mathbb{C}}, M}/\mathfrak{m}^{r+1}_{\GL_{m,\mathbb{C}}, M}$ are constant scalar multiples of the derivatives of $g$ at $M$, which are $K$-algebraic functions on $\GL_{m}$, and which can be extended to functions on $U \times \GL_{m}$ via the identity. Since the matrix representing $\tau^{r}_{s,M}$ has as entries $K$-algebraic functions on $U \times \GL_{m}$, the coordinates of the vector $\tau^{r}_{s,M}(g)$ are thus also of this form.
\end{proof}

Lastly, let us note that the proof of \autoref{jetevalprop} also gives an algorithm for computing the $K$-algebraic functions referenced in \autoref{permapgermprop} and \autoref{pullbackfunclocper}, as it essentially amounts to computing the polynomials $\xi_{i_{1},\hdots,i_{q},jk}$. 

\section{Moduli of Weakly Specials from Jets}
\label{modulisec}

In this section we fix a polarized variation of integral Hodge structure $\mathbb{V}$ on the smooth connected quasi-projective $K$-variety $S$, with $K \subset \mathbb{C}$ a subfield, and a projective compactification $S \subset \overline{S}$. We then describe the locus $\mathcal{W} \subset \textrm{Hilb}(\overline{S})$ inside the Hilbert scheme consisting of those points $[\overline{Z}]$ such that $Z = S \cap \overline{Z}$ is a weakly special subvariety for the variation $\mathbb{V}$. This will lead immediately to a proof of \autoref{mainthm2}. We postpone computational considerations until \autoref{compsec}.

\subsection{Preliminaries on Families of Jets}

In this section we develop some useful results about compatible families and sequences of jets, some of which were established in a previous paper by the author.

\begin{defn}
Suppose that $\mathcal{T}_{r} \subset (J^{d}_{r} Z)(\mathbb{C})$ for $r \geq 0$ are sets. We say that $\{ \mathcal{T}_{r} \}_{r \geq 0}$ is a \emph{compatible family} if the projections $\pi^{r}_{r-1} : J^{d}_{r} Z \to J^{d}_{r-1} Z$ restrict to give maps $\mathcal{T}_{r} \to \mathcal{T}_{r-1}$ for all $r \geq 1$.
\end{defn}

\begin{defn}
We say that $\{ j_{r} \}_{r \geq 0}$ with $j_{r} \in J^{d}_{r} Z$ is a \emph{compatible sequence} if we have $\pi^{r}_{r-1}(j_{r}) = j_{r-1}$ for all $r \geq 1$.
\end{defn}

\begin{notn}
Given an analytic space $X$ and a point $x \in X$, we denote by $(J^{d}_{r} X)_{x}$ the fibre of $J^{d}_{r} X$ above $x$.
\end{notn}

\begin{lem}
\label{compextlem}
Suppose we have a compatible family $\mathcal{T}_{r} \subset (J^{d}_{r} Z)(\mathbb{C})$, and each $\mathcal{T}_{r}$ is a constructible algebraic subset. Let $j_{r_{0}} \in \mathcal{T}_{r_{0}}$ be a jet for some fixed $r_{0}$. Then if the fibre of $\pi^{r'}_{r_{0}} : J^{d}_{r'} Z \to J^{d}_{r_{0}} Z$ above $j_{r_{0}}$ intersects $\mathcal{T}_{r'}$ for each $r' \geq r_{0}$, then $j_{r_{0}}$ belongs to a compatible sequence $\{ j_{r} \}_{r \geq 0}$ with $j_{r} \in \mathcal{T}_{r}$ for each $r \geq 0$.
\end{lem}

\begin{proof}
See \cite[Lem 5.3]{periodimages}.
\end{proof}

\begin{defn}
If $Z$ is a variety (either algebraic or analytic) and $j \in J^{d}_{r} Z$ is a jet, viewed as a map $j : \mathbb{D}^{d}_{r} \to Z$, we say that $j$ is \emph{non-degenerate} if $j$ induces an embedding on tangent spaces. The open subvariety of non-degenerate jets is denoted $J^{d}_{r,nd} Z \subset J^{d}_{r} Z$.
\end{defn}

\begin{lem}
\label{landinlem}
Let $f : (X, x) \to (Y, y)$ be a map of germs of analytic spaces, suppose that $(Z, y) \subset (Y, y)$ is an analytic subgerm, and that $(X, x)$ is smooth. Suppose we have a compatible family $j_{r} \in (J^{d}_{r,nd} X)_{x}$ of non-degenerate jets with $d = \dim (X, x)$, and that $f \circ j_{r}$ lies in $(J^{d}_{r} Z)_{y}$ for all $r$. Then $f$ factors through the inclusion $(Z, y) \subset (Y, y)$.
\end{lem}

\begin{proof}
As $(X, x)$ is smooth we may reduce to the case where $(X, x) = (\mathbb{C}^{d}, 0)$. We recall that the coordinate ring of the space $\mathbb{D}^{d}_{r}$ is $A^{d}_{r} = \mathbb{C}[t_{1}, \hdots, t_{d}]/\langle t_{1}, \hdots, t_{d} \rangle^{r+1}$. Letting $A^{d}_{\infty} = \varprojlim_{r} A^{d}_{r}$ be the associated formal power series ring, the maps $j_{r}$ glue to give a single map $\widehat{j}^{\sharp}_{\infty} : \widehat{\mathcal{O}}_{X, x} \to A^{d}_{\infty}$ of formal power series rings. The non-degeneracy of the $j_{r}$ and the formal inverse function theorem implies that $\widehat{j}^{\sharp}_{\infty}$ is an isomorphism, so the kernel of the map $\widehat{f}^{\sharp} : \widehat{\mathcal{O}}_{Y, y} \to \widehat{\mathcal{O}}_{X, x}$ agrees with the kernel of the map $\widehat{j}^{\sharp}_{\infty} \circ \widehat{f}^{\sharp}$. The fact that $f \circ j_{r}$ lies inside $(Z, y)$ for each $r$ translates to the fact that the ideal $\widehat{I} \subset \widehat{\mathcal{O}}_{Y, y}$ corresponding to $(Z, y)$ lies in the kernel of $\widehat{j}^{\sharp}_{\infty} \circ f^{\sharp}$, and hence in the kernel of $\widehat{f}^{\sharp}$. But this implies the existence of the desired factorization of the map $f^{\sharp} : \mathcal{O}_{Y, y} \to \mathcal{O}_{X, x}$ through the ring $\mathcal{O}_{Y, y} / I$ on the level of local rings, where $I$ is the ideal defining $(Z, y)$, hence of the desired factorization of $f$.
\end{proof}

\begin{lem}
\label{jetconvergence}
Suppose that $f : X \to S$ is a smooth map of $K$-schemes. Then for each $s \in S(\mathbb{C})$, we have a fibre $X_{s}$, and a jet space $J^{d}_{r} X_{s}$. Suppose that $s_{i} \to s$ is a sequence in $S(\mathbb{C})$, and $j \in (J^{d}_{r} X_{s})(\mathbb{C})$ is a jet. Then there exists $j_{i} \in (J^{d}_{r} X_{s_{i}})(\mathbb{C})$ such that $j_{i} \to j$, with the limit taken inside $(J^{d}_{r} X)(\mathbb{C})$. 
\end{lem}

\begin{proof}
From the family $f$, we obtain a family $J^{d}_{r} X \to S$, and for each point $s \in S(\mathbb{C})$ we have a subscheme $J^{d}_{r} X_{s} \subset J^{d}_{r} X$. We wish to understand $\bigcup_{s \in S(\mathbb{C})} (J^{d}_{r} X_{s})(\mathbb{C})$ as the complex points of a subscheme $X^{d}_{r} \subset J^{d}_{r} X$. We define $X^{d}_{r}$ as the fibre of $J^{d}_{r} X$ above the subscheme $c(S) \subset J^{d}_{r} S$ of constant jets. Suppose that $j \in X^{d}_{r}(\mathbb{C})$, and that $f \circ j \in (J^{d}_{r} S)(\mathbb{C})$ lies above the point $s \in S(\mathbb{C})$. Then $j$ may be viewed as a map $j : \mathbb{D}^{d}_{r, \mathbb{C}} \to X$ such that the composition $f \circ j$ factors as $f \circ j = s \circ \alpha$ where $\alpha : \mathbb{D}^{d}_{r, \mathbb{C}} \to \Spec \mathbb{C}$ is the structure map. This is all depicted in the following diagram:

\begin{center}
\begin{tikzcd}
\mathbb{D}^{d}_{r, \mathbb{C}}
\arrow[bend left]{drr}{\alpha}
\arrow[bend right,swap]{ddr}{j}
\arrow[dashed]{dr} & & \\
& X_{s} \arrow{r} \arrow[d, hook]
& \Spec{\mathbb{C}} \arrow{d}{s} \\
& X \arrow[swap]{r}{f}
& S
\end{tikzcd} .
\end{center}
The dashed arrow from the universal property of the fibre product then tells us that $j$ lies in $(J^{d}_{r} X_{s})(\mathbb{C})$.

Since the morphism $f : X \to S$ of $K$-schemes is smooth, so is the induced morphism of $\mathbb{D}^{d}_{r}$-schemes $X \times_{K} \mathbb{D}^{d}_{r} \to S \times_{K} \mathbb{D}^{d}_{r}$ obtained by base-change. The Weil-restriction functor\footnote{We note that the functor $J^{d}_{r}$ may be identified with a Weil-restriction functor; see \cite[\S2]{periodimages}.} sends smooth morphisms to smooth morphisms (see \cite[Prop. 2.2(v)]{ji1997weil}), hence so is the induced morphism $J^{d}_{r} X \to J^{d}_{r} S$. We thus obtain, by base changing the map $J^{d}_{r} X \to J^{d}_{r} S$, a smooth map $X^{d}_{r} \to c(S)$, which after using the natural identification between $c(S)$ and $S$ becomes a smooth map $f' : X^{d}_{r} \to S$. The map $f'$ is smooth and of finite presentation, hence it is in particular flat and of finite presentation. By \cite[\href{https://stacks.math.columbia.edu/tag/01UA}{Lemma 01UA}]{stacks-project} this means it is universally open. Koll\'ar has shown \cite[Thm. 31]{kollar2021fundamental} that if $g : M \to N$ is a universally open map of $\mathbb{C}$-schemes then the map $M(\mathbb{C}) \to N(\mathbb{C})$ of topological spaces satisfies the \emph{local path lifting} property: if $\gamma : [0,1] \to N(\mathbb{C})$ is a path and $m \in g^{-1}(\gamma(0))$ then there exists a path $\gamma' : [0, \ep] \to M(\mathbb{C})$ for some $\ep > 0$ such that $\gamma'(0) = m$ and $\gamma = g \circ \gamma'$. 

In our setting we have $M = X^{d}_{r}$ and $N = S$. We note that every algebraic variety can be triangulated (see \cite{Lojasiewicz1964}), so passing to a finite subsequence we may assume that the points $s_{i}$ for $i \geq 1$ and the point $s$ all lie in a subspace of $S(\mathbb{C})$ homeomorphic to a simplex. It is then clear we may construct a continuous path $\gamma : [0, 1] \to S(\mathbb{C})$ such that $\gamma(2^{-(i-1)}) = s_{i}$ and $\gamma(0) = s$. The fact that this path lifts to some path $\gamma'$ on the interval $[0, \ep]$ such that $\gamma'(0) = j$ then lets us define the $j_{i}$ for infinitely many $i$. As the claim that $j_{i} \to j$ only depends on the tail end of the sequence, this completes the proof.
\end{proof}

\subsection{Preliminaries on Weakly Specials}
\label{weakspprelimsec}

We now establish some basic facts about weakly special subvarieties that will be of use in the coming sections. We will write $V = \mathbb{Z}^{m}$, fix a polarization $Q_{V} : V \otimes V \to \mathbb{Z}$ compatible with a polarization $Q : \mathbb{V} \otimes \mathbb{V} \to \mathbb{Z}$ on $\mathbb{V}$, and let $\ch{D} \subset \ch{L}$ denote the $\mathbb{Q}$-algebraic subvariety of flags which satisfy the first Hodge-Riemann bilinear relation with respect to $Q_{V}$. There is a natural analytic open subset $D \subset \ch{D}$ consisting of polarized Hodge structures. If we let $G = \textrm{Aut}(V, Q_{V}) \subset \GL(V)$ denote the subgroup of polarization-preserving automorphisms, and let $\Gamma = G(\mathbb{Z})$, then we obtain from the variation $\mathbb{V}$ a canonical analytic period map $\varphi : S \to \Gamma \backslash D$ which sends the point $s$ to the isomorphism class of the integral polarized Hodge structure $(\mathbb{V}_{s}, Q_{s})$.

\begin{defn}
Given a Hodge structure $t \in D$, denote by 
\[ \textrm{Hg}(t) \subset \bigoplus_{m,n \geq 1}^{\infty} V^{\otimes m} \otimes (V^{*})^{\otimes n} \]
the subspace of all Hodge tensors associated to the Hodge structures tensorially generated by $t$. We denote by $\textrm{MT}(t) \subset G$ the Mumford-Tate group of $t$, which is the stabilizer of these tensors.
\end{defn}

\begin{prop}
\label{weakspredef}
A weakly special subvariety is exactly an irreducible component of a variety of the form $\varphi^{-1}(\pi(N(\mathbb{R}) \cdot t))$, where $N \subset \textrm{MT}(t)$ is a $\mathbb{Q}$-algebraic normal subgroup, $t \in D$ is a point, and $\pi : D \to \Gamma \backslash D$ is the natural projection.
\end{prop}

\begin{proof}
This is easily deduced from \cite[Def. 3.1, Cor. 3.14]{closurepositivelocus}. See also \cite[Prop. 3.8]{closurepositivelocus} for a proof that such varieties are algebraic.
\end{proof}

\begin{lem}
\label{monocont}
Suppose that $M_{i} \subset G$ for $i \in \{ 1, 2 \}$ are two Mumford-Tate groups, that $N_{1} \subset M_{1}$ and $N_{2} \subset M_{2}$ are two connected semisimple $\mathbb{Q}$-algebraic normal subgroups, and let $h \in D$ be a Hodge structure with $\textrm{MT}(h) \subset M_{1} \cap M_{2}$. Suppose moreover that each $\mathbb{Q}$-simple factor of $N_{1}$ does not have a compact set of real points. Then if $N_{1} \cdot h \subset N_{2} \cdot h$ we in fact have $N_{1} \subset N_{2}$.
\end{lem}

\begin{proof}
Letting $M_{3} = M_{1} \cap M_{2}$ and $N_{3} = N_{1} \cap N_{2}$, it suffices to show that $N_{1} = N_{3}$. Since elements of $M_{3}$ normalize both $N_{1}$ and $N_{2}$, the problem is unchanged by replacing $h$ with $m \cdot h$ for some $m \in M_{3}(\mathbb{R})$, so we may assume that $\textrm{MT}(h) = M_{3}$. We begin by showing that $N_{3} \cdot h = N_{1} \cdot h$; indeed, since these are group orbits, it suffices to check the equality of tangent spaces $T_{h} (N_{3} \cdot h) = T_{h} (N_{1} \cdot h)$ at $h$, and using using the openness of the real orbits of $h$ (c.f. \cite[pg. 56]{GGK} on the openness of real orbits of Mumford-Tate groups) we may reduce to checking the statement $T_{h} (N_{3}(\mathbb{R}) \cdot h) = T_{h} (N_{1}(\mathbb{R}) \cdot h)$. Analogously to \cite[pg.54]{GGK} we may model these tangent spaces in the following way: the Hodge structure $h : \mathbb{S} \to \GL(V)_{\mathbb{R}}$ induces, through the adjoint action on the Lie algebras $\mathfrak{n}_{i}$ of $N_{i}$ for $i = 1, 2, 3$, a weight zero Hodge structure. We then obtain a decomposition $\mathfrak{n}_{i,\mathbb{C}} = \mathfrak{s}_{i} \oplus \mathfrak{p}_{i}$ into summands of negative and non-negative weight, and such that the orbit map $N_{i}(\mathbb{C}) \to N_{i}(\mathbb{C}) \cdot h$ defined by $h$ identifies $\mathfrak{s}_{i}$ with $T_{h}(N_{i}(\mathbb{C}) \cdot h) = T_{h} (N_{i}(\mathbb{R}) \cdot h)$ and $\mathfrak{p}_{i}$ with the complex Lie algebra of the parabolic subgroup of $N_{i}$ stabilizing $h$. In particular, since $N_{1} \cdot h \subset N_{2} \cdot h$ we have $\mathfrak{s}_{1} \subset \mathfrak{s}_{2}$. Finally, these decompositions satisfy the property that $\mathfrak{s}_{3} = \mathfrak{s}_{1} \cap \mathfrak{s}_{2}$, from which it follows that $\mathfrak{s}_{1} = \mathfrak{s}_{3}$ and $T_{h} (N_{3}(\mathbb{R}) \cdot h) = T_{h}(N_{1}(\mathbb{R}) \cdot h)$.

The equality $\mathfrak{s}_{1} = \mathfrak{s}_{3}$ together with the Hodge symmetry of the (polarizable) Hodge structures on $\mathfrak{n}_{1}$ and $\mathfrak{n}_{3}$ imply that the containment $\mathfrak{n}_{3} \subset \mathfrak{n}_{1}$ of polarizable Hodge structures decomposes as $\mathfrak{n}_{1} = \mathfrak{n}_{3} \oplus \mathfrak{a}$, where $\mathfrak{a}$ is a $\mathbb{Q}$ sub-Hodge structure concentrated entirely in degree zero. Since $\mathfrak{a}$ is defined over $\mathbb{Q}$ and concentrated in degree zero, it is spanned entirely by endomorphisms of the Hodge structure $h$, so elements of $\mathfrak{a}$ are fixed by the conjugation action of $\textrm{MT}(h) = M_{3}$, and in particular by $N_{3}$. It follows that $[\mathfrak{n}_{3}, \mathfrak{a}] = 0$, so $\mathfrak{a}$ is an ideal of $\mathfrak{n}_{1}$, and hence a semisimple summand of $\mathfrak{n}_{1}$. If we let $A = \textrm{exp}(\mathfrak{a})$, then since $\mathfrak{a} \subset \mathfrak{p}_{1}$ we find that $A$ lies in the stabilizer $P_{1} \subset N_{1}$ of the Hodge flag corresponding to $h$. Applying \cite[IV.A.2]{GGK} the closed subgroup of $N_{1}(\mathbb{R})$ stabilizing $h$ is compact, so it follows that $A(\mathbb{R})$ is compact. We know by assumption that $N_{1}$ does not have $\mathbb{Q}$-simple factors with compact sets of real points, hence $\mathfrak{a} = 0$.
\end{proof}

\medskip

Another useful property of weakly special subvarieties will be the following consequence of the Ax-Schanuel Theorem for variations of Hodge structures proven by Bakker and Tsimerman \cite{AXSCHAN}. In what follows we will treat algebraic monodromy groups $\mathbf{H}_{Z}$ associated to subvarieties $Z \subset S$ as subgroups of $G$. Strictly speaking, the group $\mathbf{H}_{Z}$ only defines a canonical conjugacy class of subgroup of $G$, but what we mean is to consider $\mathbf{H}_{Z}$ as a subgroup of $G$ relative to a polarization-preserving identification $V \simeq \mathbb{V}(B)$ used to define a local lift $\psi : B \to D$ of the fixed period map $\varphi$. When we view algebraic monodromy groups $\mathbf{H}_{Z}$ as subgroups of $G$ in the coming sections we will always mean with respect to some (sometimes implicitly) fixed local lift.

\begin{lem}
\label{axshanres}
Suppose that $\psi : B \to \ch{D}$ is a local lift of $\varphi$, and let $\ch{E} \subset \ch{D}$ be an algebraic subvariety. Fix a point $s \in B$ such that $\psi(s) \in \ch{E}$. Then:
\begin{itemize}
\item[(i)] an analytic germ $(C, s) \subset (B, s)$ of an irreducible algebraic subvariety, maximal among such germs for the property that $\psi(C, s) \subset (\ch{E}, \psi(s))$, is the germ of a weakly special subvariety;
\item[(ii)] if $(C, s)$ is as in (i) an analytic germ of the algebraic subvariety $Y \subset S$, then the Zariski closure of $\psi(C, s)$ is the orbit $\mathbf{H}_{Y} \cdot \psi(s)$, where $\mathbf{H}_{Y}$ is the algebraic monodromy of $Y$.
\item[(iii)] if we have irreducible analytic germs $(C_{i}, s) \subset (B, s)$ of algebraic $Y_{i} \subset S$ for $i = 1, 2$ such that $\mathbf{H}_{Y_{1}} \cdot \psi(s) \subset \mathbf{H}_{Y_{2}} \cdot \psi(s)$ then we have $\mathbf{H}_{Y_{1}} \subset \mathbf{H}_{Y_{2}}$ (inside $\GL(\mathbb{V}_{s})$).
\end{itemize}
\end{lem}

\begin{proof}
We start with (iii), which follows from \autoref{monocont} above if we can show that algebraic monodromy groups do not contain $\mathbb{Q}$-simple factors with compact sets of real points. This is shown below in \autoref{nocomprealfac}. Turning to (ii), if $(C, s) \subset (B, s)$ is any analytic germ of an algebraic subvariety $Z \subset S$, it follows as in \cite[Lem. 4.7]{periodimages} from the Ax-Schanuel Theorem that the Zariski closure of $\psi(C, s)$ in $\ch{D}$ is the orbit $\mathbf{H}_{Z} \cdot \psi(s)$, where $\mathbf{H}_{Z}$ the algebraic monodromy group of $Z$. The result (i) then follows from our original definition of weakly special in \autoref{weakspdef} and (iii).
\end{proof}

\begin{lem}
\label{nocomprealfac}
If $\mathbb{V}$ is a variation of Hodge structure on an irreducible algebraic variety $Z$, then the algebraic monodromy group of $Z$ does not admit a $\mathbb{Q}$-simple factor with a compact set of real points.
\end{lem}

\begin{proof}
Noting that algebraic monodromy is an invariant unchanged by replacing $Z$ with an open subvariety, we may assume $Z$ is smooth. The reason for this is then as follows: according to \cite[III.A.1]{GGK}, the monodromy group $\Gamma_{Z}$ of a variation of Hodge structure on a variety $Z$ splits, after possibly replacing $Z$ with a finite \'etale covering, as $\Gamma_{Z} = \Gamma_{1} \times \cdots \times \Gamma_{k}$ where the Zariski closure $\overline{\Gamma_{i}}^{\textrm{Zar}} = N_{i}$, and $N_{1} \cdots N_{k}$ is the $\mathbb{Q}$-simple decomposition of $\mathbf{H}_{Z}$. Since algebraic monodromy groups are unchanged by passing to finite \'etale coverings, none of the $N_{i}$ can have a compact set of real points as then the equality $\overline{\Gamma_{i}}^{\textrm{Zar}} = N_{i}$ would not hold due to the fact that $\Gamma_{i}$ would have to be finite, and algebraic monodromy groups are assumed connected.
\end{proof}

\subsection{Differential Constraints and Types}
\label{typeintrosec}

We continue with the notation in the previous section. In this section we introduce an invariant of a geometrically irreducible subvariety $Z \subset S$ induced by the variation $\mathbb{V}$ called its \emph{type}, and use it to both construct certain so-called ``differential constraints'' satisfied by all varieties of a certain type and give an alternate description of the weakly special subvarieties.

Recall that we have a fixed a polarized lattice $(V, Q_{V})$ with $V = \mathbb{Z}^{m}$, and that we have an associated $\mathbb{Q}$-algebraic variety $\ch{L}$ consisting of all Hodge flags on $V$ on which $\GL(V) = \GL_{m}$ acts. The variety $\ch{D}$ consisting of all Hodge flags satisfying the first Hodge-Riemann bilinear relation we view as a $\mathbb{Q}$-algebraic subvariety $\ch{D} \subset \ch{L}$, and the open submanifold $D \subset \ch{D}$ of polarized Hodge flags we similarly view as embedded in $\ch{L}$. We also have a canonical period map $\varphi : S \to \Gamma \backslash D$, where $\Gamma = \textrm{Aut}(V, Q_{V})(\mathbb{Z})$, sending a polarized Hodge structure to its isomorphism class.

\begin{defn}
Given two closed subvarieties $V_{1}, V_{2} \subset \ch{L}$, we say that $V_{1} \sim_{\GL} V_{2}$ if there exists $g \in \GL_{m}(\mathbb{C})$ such that $V_{2} = g \cdot V_{1}$.
\end{defn}

\begin{defn}
Given a closed subvariety $V \subset \ch{L}$, we call the equivalence class of $V$ under $\sim_{\GL}$ a \emph{type}, and use the notation $\mathcal{C}(V)$ to denote this equivalence class. We will also say \emph{the type associated to }$V$ to refer to $\mathcal{C}(V)$.
\end{defn}

\begin{defn}
Given two types $\mathcal{C}_{1}$ and $\mathcal{C}_{2}$, we write $\mathcal{C}_{1} \leq \mathcal{C}_{2}$ if there exists $V_{1}, V_{2} \subset \ch{L}$ such that $V_{1} \subset V_{2}$ and $\mathcal{C}_{i} = \mathcal{C}(V_{i})$ for $i = 1, 2$.
\end{defn}

\begin{defn}
Let $\psi : B \to D$ be a local lift of $\varphi$, viewed as a map $B \to \ch{L}$ using the fixed embedding $D \subset \ch{L}$. Let $C \subset B \cap Z$ be an analytic component with $Z \subset S$ a geometrically irreducible algebraic subvariety. Then the \emph{type} of $Z$ is the equivalence class $\mathcal{C}(\overline{\psi(C)})$, where $\overline{\psi(C)}$ is the Zariski closure of $\psi(C)$. We use the notation $\mathcal{C}(Z)$ to denote the type of $Z$.
\end{defn}

\begin{lem}
The type of $Z \subset S$ is independent of the lift $\psi$ and the component $C$ chosen, and as $Z$ ranges over all geometrically irreducible algebraic subvarieties of $S$ only finitely many types arise.
\end{lem}

\begin{proof}
For the first statement we note that the different Zariski closures of $\psi(C)$, all of which are orbits of $\mathbf{H}_{Z}$ by \autoref{axshanres}(ii), will be related by elements of $G(\mathbb{Z}) \subset \GL_{m}(\mathbb{C})$ since a local lift of $\psi$ at a point is determined up to a monodromy group lying inside $G(\mathbb{Z})$. Moreover one can see that varying the neighbourhood $B$ intersecting $Z$ or the component $C$ will not change the Zariski closure using irreducibility and analytic continuation.

For the second claim we argue as follows. Recall that we have a fixed subvariety $\ch{D} \subset \ch{L}$ of Hodge flags satisfying the first Hodge-Riemann bilinear relation. Points $t \in \ch{D}$ may be identified with complex cocharacters $\mu : \mathbb{C}^{\times} \to \GL_{m,\mathbb{C}}$, with the action on $t$ by $\GL_{m}(\mathbb{C})$ corresponding to the action on $\mu$ by conjugation; a general version of this fact appears in \cite[VI.B.9]{GGK}, and it is also used in the proof of \cite[Thm. 4.14]{hodgelocivoisin}. Moreover, those points $t \in D \subset \ch{D}$ that correspond to polarized Hodge structures with Mumford-Tate group $M$ necessarily give rise by \cite[VI.B.9]{GGK} to cocharacters $\mu$ which factor through $M_{\mathbb{C}}$. It will be enough to show the stronger claim that if $M \subset G$ is a Mumford-Tate group, $t \in \ch{D}$ is a point whose cocharacter $\mu_{t} : \mathbb{C}^{\times} \to \GL_{m,\mathbb{C}}$ factors through $M_{\mathbb{C}}$, and $N \subset M_{\mathbb{C}}$ is a \emph{complex} semisimple normal algebraic subgroup, then there are finitely many possibilities for the orbit $N \cdot t \subset \ch{D}$ up to translation by $\GL_{m}(\mathbb{C})$. We use the fact, proven in \cite[Thm. 4.14]{hodgelocivoisin}, that there are finitely many $\GL_{m}(\mathbb{C})$-conjugacy classes of Mumford-Tate groups. We claim this implies the stronger fact that triples $(M, N, t)$ consiting of a Mumford-Tate group $M$, a complex semisimple normal subgroup $N \subset M_{\mathbb{C}}$, and a point $t$ whose cocharacter $\mu_{t}$ factors through $M_{\mathbb{C}}$, are finite up to $\GL_{m}(\mathbb{C})$-equivalence. (The action of $\GL_{m}(\mathbb{C})$ on $M$ and $N$ is by conjugacy and on $t$ by translation using the embedding $\ch{D} \subset \ch{L}$.) To constrain the possibilites for $N$ we observe that there are finitely many possibilities for the normal factor $N$ of $M$. To constrain the possibilities for $t$ it have that by \cite[VI.B.9]{GGK} the points $t$ whose cocharacters $\mu_{t}$ factor through $M_{\mathbb{C}}$ belong to finitely many $M(\mathbb{C})$-orbits (the components of $\ch{\textrm{NL}}_{M}$ in the notation of \cite{GGK}) in $\ch{D}$. Finally, if two triples $(M_{1}, N_{1}, t_{1})$ and $(M_{2}, N_{2}, t_{2})$ are in the same $\GL_{m}(\mathbb{C})$-orbit, we then have for some $g \in \GL_{m}(\mathbb{C})$ that
\[ N_{1} \cdot t_{1} = (g \cdot N_{2} \cdot g^{-1}) \cdot (g \cdot t_{2}) = g \cdot (N_{2} \cdot t_{2}) , \]
hence the result.
\end{proof}

\medskip

We now fix an enumeration $\mathcal{C}_{1}, \hdots, \mathcal{C}_{q}$ of the types of geometrically irreducible subvarieties $Z \subset S$ arising from our fixed variation $\mathbb{V}$. We fix representatives $(N_{1}, t_{1}), \hdots, (N_{q}, t_{q})$ of each type such that $\mathcal{C}_{i} = \mathcal{C}(N_{i} \cdot t_{i})$.

\begin{notn}
Given a variety $V \subset \ch{L}$ we write
\[ \mathcal{T}^{d}_{r}(V) = \mathcal{T}^{d}_{r}(\mathcal{C}(V)) = (\eta^{d}_{r})^{-1}(\GL_{m} \cdot J^{d}_{r} V) , \]
where $\eta^{d}_{r}$ is as in \autoref{constrjetcorrespthm}. For $Z \subset S$, we set $\mathcal{T}^{d}_{r}(Z) = \mathcal{T}^{d}_{r}(\mathcal{C}(Z))$, where $\mathcal{C}(Z) = \mathcal{C}(N_{i} \cdot t_{i})$ y $i$.
\end{notn}

\begin{defn}
Sets of the form $\mathcal{T}^{d}_{r}(\mathcal{C})$ for some type $\mathcal{C}$ we call \emph{differential constraints}. If for some $Y \subset S$ we have $J^{d}_{r} Y \subset \mathcal{T}^{d}_{r}(\mathcal{C})$, we say that $Y$ \emph{satisfies} the differential constraint $\mathcal{T}^{d}_{r}(\mathcal{C})$.
\end{defn}

We note that by the defining property of the correspondence in \autoref{constrjetcorrespthm} we have that $Z$ satisfies the differential constraint $\mathcal{T}^{d}_{r}(Z)$ for all $d$ and $r$. We now use this to give the following characterization of weakly special subvarieties.


\medskip

\begin{prop}
\label{weakspcharprop}
Let $S$ be a smooth connected algebraic variety admitting a polarized integral variation of Hodge structure $\mathbb{V}$, and suppose that $Z \subset S$ is a geometrically irreducible algebraic subvariety. Then the following conditions are equivalent:
\begin{itemize}
\item[(i)] $Z$ is weakly special;
\item[(ii)] $Z$ is a maximal geometrically irreducible algebraic subvariety of type $\mathcal{C}(Z)$;
\item[(iii)] $Z$ is maximal among geometrically irreducible algebraic subvarieties $Y \subset S$ which satisfy $\mathcal{T}^{d}_{r}(Z)$ for all $d$ and $r$.
\end{itemize}
\end{prop}

\medskip

The key input to \autoref{weakspcharprop} is the following:

\begin{lem}
\label{liesintrans}
Suppose that $Y \subset S$ and $V \subset \ch{L}$ are subvarieties, and that $Y$ satisfies $\mathcal{T}^{d}_{r} = \mathcal{T}^{d}_{r}(V)$ for all $d$ and $r$. Let $\psi : B \to \ch{L}$ be a local period map, and suppose that $C \subset Y \cap B$ is an irreducible analytic component. Then $\psi(C)$ lies in a $\GL_{m}(\mathbb{C})$-translate of the variety $V$.
\end{lem}

\begin{proof}
Choose a point $s \in C$ in the smooth locus of $C$, and let $j_{r} \in J^{d}_{r,nd} C$ for $r \geq 0$ be a sequence of compatible non-degenerate jets with $j_{0} = s$ and $d = \dim C$. By the defining property of the $\eta^{d}_{r}$ in \autoref{constrjetcorrespthm}, the jets $\psi \circ j_{r}$ all lie inside $\GL_{m} \cdot J^{d}_{r} V$. Let $\mathcal{G}_{r} \subset \GL_{m}$ be the constructable algebraic subset of $g \in \GL_{m}(\mathbb{C})$ for which $\psi \circ j_{r}$ lies inside $g \cdot J^{d}_{r} V$. This gives a descending chain of constructible sets, with $\mathcal{G}_{r+1} \subset \mathcal{G}_{r}$, non-empty at every finite stage; from constructibility it follows that the intersection $\mathcal{G}_{\infty} = \bigcap_{r} \mathcal{G}_{r}$ is non-empty. Taking $g \in \mathcal{G}_{\infty}$, we find that $\psi(C) \subset g \cdot V$ by applying \autoref{landinlem} with $(X, x) = (C, s)$, $(Y, y) = (\ch{L}, \psi(s))$, and $(Z, y) = (g \cdot V, \psi(s))$, and using the irreducibility of $C$. 
\end{proof}

\medskip

\begin{proof}[Proof of \autoref{weakspcharprop}:]
To see that $(i) \implies (ii)$, suppose $Z$ is weakly special and that $Z \subset Y$ is an inclusion of two geometrically irreducible subvarieties of the same type. Then fixing a local lift $\psi : B \to D$ and components $C_{Z} \subset C_{Y}$ of $Z \cap B$ and $Y \cap B$, we learn from \autoref{axshanres}(iii) and the equality $\overline{\psi(C_{Z})} = \overline{\psi(C_{Y})}$ that $\mathbf{H}_{Y} \subset \mathbf{H}_{Z}$, so necessarily $Y = Z$ as $Z$ is weakly special. 

To see that $(ii) \implies (iii)$ we use the fact that if $Y \subset S$ satisfies the differential constraints $\mathcal{T}^{d}_{r}(Z)$ for all $d$ and $r$, then by \autoref{liesintrans} we have $\mathcal{C}(Y) \leq \mathcal{C}(Z)$. Thus if $Z$ is a maximal variety of type $\mathcal{C}(Z)$ it is also a maximal variety satisfying the constraints $\mathcal{T}^{d}_{r}(Z)$. 

Finally, to see that $(iii) \implies (i)$, suppose that $Z \subset S$ is maximal among irreducible algebraic subvarieties that satisfy $\mathcal{T}^{d}_{r}(Z)$ for all $d$ and $r$, and suppose we have an irreducible algebraic $Y$ such that $Z \subset Y \subset S$ and $Y$ has the same algebraic monodromy group as $Z$. It then follows from the definitions and \autoref{axshanres}(ii) that $\mathcal{T}^{d}_{r}(Y) = \mathcal{T}^{d}_{r}(Z)$ for all $d$ and $r$, so since $Y$ satisfies $\mathcal{T}^{d}_{r}(Y)$ for all $d$ and $r$, it follows that that $Y = Z$.
\end{proof}

\subsection{Proving Maximality}
\label{provmaxsec}

We will soon see that the results of the preceding section will enable us to construct, for a fixed type $\mathcal{C}_{i}$ arising from the variation $\mathbb{V}$, families $g : \mathcal{Z} \to T$ of subvarieties of $S$ such that for every $t \in T(\mathbb{C})$ we have $\mathcal{C}(\mathcal{Z}_{t}) \leq \mathcal{C}_{i}$. To understand weakly specials, we will then wish to understand the condition that the fibres of $g$ be maximal algebraic subvarieties of $Z$ of type $\mathcal{C}_{i}$ in terms of constructible conditions on $T$. In this section we give a method for doing this.

Recall that weakly special subvarieties $Z$ of $S$ are defined in \autoref{weakspdef} using the monodromy of the normalization $Z^{\textrm{nor}} \to Z$. Since the natural map of fundamental groups induced by a Zariski open subvariety of a normal variety is surjective, it follows that $Z$ has the same algebraic monodromy as any open subvariety of $Z$. In particular, the condition that $Z$ be weakly special in $S$ can be checked after replacing $S$ with an open subvariety $U \subset S$ which intersects $Z$. In what follows we will therefore assume that $S$ is affine and that $\Omega^{1}_{S}$ and the Hodge bundles $F^{\bullet} \mathcal{H}$ over $S$ are all trivial; we note that because $S$ is smooth one can always find a finite open cover whose constituent open sets have these properties.

In what follows we call a subscheme $A$ of $S_{\mathbb{C}}$ with underlying reduced space a point $p \in S(\mathbb{C})$ and which is contained in the $r$'th order infinitesimal neighbourhood of $p$ an (infinitesimal) \emph{germ of order} $r$. We fix a collection of functions $z_{1}, \hdots, z_{n}$ in the coordinate ring of $S$ such that $dz_{1}, \hdots, dz_{n}$ spans $\Omega^{1}_{S}$, as well as a filtration-compatible frame $v^{1}, \hdots, v^{m}$ for $\mathcal{H}$. We recall from \autoref{torsorconssec} that for each point $s \in S(\mathbb{C})$ and $M \in \GL_{m}(\mathbb{C})$ we have a map $f : B \to \GL_{m}$ with $s \in B$ satisfying the property that $f(s) = M$ and $b^{k} = \sum_{i = 1}^{m} f_{ik} v^{i}$ gives a flat frame on $B$ whose germ at $s$ is uniquely determined by the choice of $M$. Moreover, the results of \autoref{germsofpermaps} show that with respect to the natural bases determined by $z_{1}, \hdots, z_{n}$ and the natural coordinates on $\GL_{m}$, the maps of complex vector spaces $\tau^{r}_{s,M} : \widehat{\mathcal{O}}_{\GL_{m,\mathbb{C}}, M}/\mathfrak{m}^{r+1}_{\GL_{m,\mathbb{C}}, M} \to \widehat{\mathcal{O}}_{S_{\mathbb{C}}, s}/\mathfrak{m}^{r+1}_{S_{\mathbb{C}}, s}$ induced by the analytic maps $f$ have their matrix entries given by $K$-algebraic functions on $S \times \GL_{m}$. We note here that the ``basis determined by $z_{1}, \hdots, z_{n}$'' inside $\widehat{\mathcal{O}}_{S_{\mathbb{C}}, s}/\mathfrak{m}^{r+1}_{S_{\mathbb{C}}, s}$ is the basis of monic monomials in $(z_{1} - s_{1}), \hdots, (z_{n} - s_{n})$, where $s_{i}$ is the value of $z_{i}$ at $s$.

\begin{prop}
\label{Agermloc}
Fix a variety $V \subset \ch{L}$ and a $K$-algebraic family $g : \mathcal{Z} \to T$ of subvarieties of $S$ with projection $p : \mathcal{Z} \to S$ restricting to an embedding on fibres of $g$. Consider the subset $\mathcal{A}(g, V, r) \subset \mathcal{Z} \times \GL_{m}$ given by
\[ \left\{ (z, M) : \begin{pmatrix} \textrm{ the order }r\textrm{ germ of } p(\mathcal{Z}_{g(z)}) \textrm{ at }p(z)\textrm{ lies in } \\
\psi^{-1}(V) \textrm{ where } \psi = q \circ \iota \circ f \textrm{ with }f(s) = M \end{pmatrix} \right\} . \]
Then $\mathcal{A}(g, V, r)$ is a constructible $K$-algebraic subset of $\mathcal{Z} \times \GL_{m}$. 
\end{prop}

\begin{proof}
By replacing $V$ with $W = (q \circ \iota)^{-1}(V)$ it suffices to consider
\[ \left\{ (z, M) : \begin{pmatrix} \textrm{ the order }r\textrm{ germ of } p(\mathcal{Z}_{g(z)}) \textrm{ at }p(z) \\
\textrm{ lies in } f^{-1}(W) \textrm{ where }f(s) = M \end{pmatrix} \right\} . \]
The variety $W$ may be assumed to be the restriction to $\GL_{m}$ of an affine subvariety of the space $\mathbb{M}$ of $m \times m$ matrices, and by realizing $W$ as an intersection of varieties of codimension one we may reduce to the situation where $W$ is defined by a single algebraic function $h$ in the natural coordinates $[a_{jk}]$ with $1 \leq j, k \leq m$ on $\mathbb{M}$. Applying \autoref{pullbackfunclocper} we find that the coordinates of $\tau^{r}_{s,M}(h)$ with respect to the basis induced by the coordinates $z_{1}, \hdots, z_{n}$ on $S$ are $K$-algebraic functions on $S \times \GL_{m}$.

Using the maps $p$ and $g$, we may view $\mathcal{Z}$ as a closed subvariety of the product $T \times S$, and reduce to realizing the above as a constructible locus inside $T \times S \times \GL_{m}$. We may reduce to the case where $T = \Spec B$ is affine, and since $g$ is a family of affine subvarieties of $S = \Spec R$, we may assume $\mathcal{Z}$ is defined inside $T \times S$ by functions $j_{1}, \hdots, j_{k} \in B \otimes_{K} R$; we will denote by $j_{i,t}$ their restrictions to $\{ t \} \times S$. The basis $dz_{1}, \hdots, dz_{n}$ for $\Omega^{1}_{S}$ admits a dual basis $\partial_{1}, \hdots, \partial_{n}$ for the tangent sheaf, which we may identify with the sheaf of algebraic derivations of functions on $S$. Moreover, we can extend the operators $\partial_{1}, \hdots, \partial_{n}$ to $T \times S$ by regarding functions on $T$ as constants. Then in the basis for $\widehat{\mathcal{O}}_{S_{\mathbb{C}}, s}/\mathfrak{m}^{r+1}_{S_{\mathbb{C}}, s}$ induced by the coordinates $z_{1}, \hdots, z_{n}$, the coefficients of the image of $j_{i,t}$ are given by constant scalar multiples of the values at $s$ of the partial derivatives $\partial_{i_{1} \cdots i_{q}} j_{i,t}$ up to order $r$. In particular, the coefficients of the image of $j_{i,t}$ in $\widehat{\mathcal{O}}_{S_{\mathbb{C}}, s}/\mathfrak{m}^{r+1}_{S_{\mathbb{C}}, s}$ are algebraic functions on $T \times S$.

The conditions on points $(t, s, M) \in T \times S \times \GL_{m}$ we wish to enforce are then the following. First, we must have $s \in \mathcal{Z}_{t}$, which simply comes from requiring that $(t, s)$ lie in the image of the embedding $\mathcal{Z} \hookrightarrow T \times S$. Secondly we require that the ideal generated by $j_{1,t}, \hdots, j_{k,t}$ inside $\widehat{\mathcal{O}}_{S_{\mathbb{C}}, s}/\mathfrak{m}^{r+1}_{S_{\mathbb{C}}, s}$ contain the function $h$. The second condition reduces to a linear algebraic condition on the coefficients of $j_{1,t}, \hdots, j_{k,t}$ and $\tau^{r}_{s,M}(h)$ written in the natural basis: in particular, one may obtain a finite spanning set of the ideal in $\widehat{\mathcal{O}}_{S_{\mathbb{C}}, s}/\mathfrak{m}^{r+1}_{S_{\mathbb{C}}, s}$ generated by $j_{1,t}, \hdots, j_{k,t}$ by multiplying by monomials in $(z_{1} - s_{1}), \hdots, (z_{n} - s_{n})$, where $s_{i}$ is the value of $z_{i}$ at $s$. The vectors in such a spanning set $\mathcal{S}(s,t)$ have coordinates which are algebraic functions on $T \times S$ (because this is true for $j_{1,t}, \hdots, j_{k,t}$). Moreover, the condition that $h$ lie in the ideal generated by $j_{1, t}, \hdots, j_{k,t}$ can described by saying that there exists a linear combination of the vectors in $\mathcal{S}(s,t)$ which is equal to $\tau^{r}_{s,M}(h)$. Because the coefficients of the elements of $\mathcal{S}(s, t) \cup \{ h \}$ are all $K$-algebraic functions on $T \times S \times \GL_{m}$, this is a constructible $K$-algebraic condition on $T \times S \times \GL_{m}$. 
\end{proof}

To give our criterion for maximality we will also need a way of understanding how the size of a germ, measured by the vector space dimension of its coordinate ring, relates to the dimension of the variety it comes from. For this we use the following:

\begin{defn}
We denote by $H(r,d)$ the vector space dimension 
\[ H(r,d) := \dim_{\mathbb{C}} \mathbb{C}[t_{1}, \hdots, t_{d}]/\langle t_{1}, \hdots, t_{d} \rangle^{r+1} . \]
\end{defn}

\begin{lem}
\label{angermlarge}
Let $A$ an analytic variety of pure dimension $d$, and let $a \in A$ be a point. Then we have $\dim_{\mathbb{C}} \widehat{\mathcal{O}}_{A,a} / \mathfrak{m}^{r+1}_{A,a} \geq H(r,d)$. 
\end{lem}

\begin{proof}
Using the density of the smooth locus of $A$, it suffices to show that the function $a \mapsto \dim_{\mathbb{C}} \widehat{\mathcal{O}}_{A,a} / \mathfrak{m}^{r+1}_{A,a}$ is upper semi-continuous on $A$. This question is local on $A$, so after replacing $A$ with an open neighbourhood we may assume that $A$ is a closed analytic subvariety of $\mathbb{C}^{n}$, defined as the vanishing locus of an ideal $I = \langle f_{1}, \hdots, f_{k} \rangle$ generated by functions holomorphic in the coordinates $x_{1}, \hdots, x_{n}$ on $\mathbb{C}^{n}$. We will write $a = (a_{1}, \hdots, a_{n}) \in \mathbb{C}^{n}$. The question then reduces to showing that $a \mapsto \dim_{\mathbb{C}} I_{a}$ is lower semi-continuous, where $I_{a} \subset \mathcal{O}_{\mathbb{C}^{n}, a} / \mathfrak{m}^{r+1}_{\mathbb{C}^{n},a}$ is the image of $I$. With respect to the canonical bases for the germ $\mathcal{O}_{\mathbb{C}^{n}, a}/\mathfrak{m}^{r+1}_{\mathbb{C}^{n},a}$ given by monomials in $(x_{1} - a_{1}), \hdots, (x_{n} - a_{n})$, the coordinates for the elements of $I_{a}$ are given by the values at $a$ of polynomial functions in the partial derivatives of the generators $f_{1}, \hdots, f_{k}$; in particular, the condition that $\dim I_{a} \geq k$ for some $k$ may be enforced using determinant conditions defining an open locus in $a$. The result follows.
\end{proof}



In what follows we assume that the varieties in the family $g : \mathcal{Z} \to T$ have a common dimension $d$. Choose $V_{i} \subset \ch{L}$ such that $\mathcal{C}_{i} = \mathcal{C}(V_{i})$, where $\mathcal{C}_{i}$ is one of the types arising from the variation $\mathbb{V}$. The locus $\mathcal{A}(g, V_{i}, r)$ has a sublocus $\mathcal{B}(g, V_{i}, r)$ consisting of points $(z, M)$ where for every $r' \leq r$ the coordinate ring $\widehat{\mathcal{O}}_{\psi^{-1}(V), p(z)} / \mathfrak{m}^{r'+1}_{\psi^{-1}(V), p(z)}$ has vector space dimension $\geq H(r',d+1)$; this locus may be seen to be $K$-algebraically constructible, since the generators of the ideal defining the order-$r'$ infinitesimal neighbourhood of $\psi^{-1}(V)$ at $p(z)$ have coordinates which are $K$-algebraic functions on $S \times \GL_{m}$, and this is a linear-algebraic condition on the vector space dimension of this ideal. The complement of the projection of $\mathcal{B}(g, V_{i}, r)$ to $T$ we denote by $\mathcal{U}(g, V_{i}, r)$. Informally, this is the locus of $t \in T(\mathbb{C})$ such that no germ of $p(\mathcal{Z}_{t})$ is contained in an order $r$ infinitesimal germ $A$ of $\psi^{-1}(V_{i})$ for a local period map $\psi$ such that the coordinate ring of $A$ is large enough that the germ of $\psi^{-1}(V_{i})$ at $s$ could have dimension $d+1$.

\begin{prop}
\label{maxcondloc}
Suppose that $t \in T(\mathbb{C})$ is a point such that $\mathcal{C}(p(\mathcal{Z}_{t})) = \mathcal{C}_{i}$. Then $p(\mathcal{Z}_{t})$ is weakly special of type $\mathcal{C}_{i}$ if and only if the point $t$ lies in the union $\bigcup_{r \geq 0} \mathcal{U}(g, V_{i}, r)$.
\end{prop}

\begin{proof}
Let us first suppose that $Z = p(\mathcal{Z}_{t})$ is not weakly special of type $\mathcal{C}_{i} = \mathcal{C}(V_{i})$. Choose a point $s \in Z(\mathbb{C})$ and a local lift $\psi : B \to D$ of the canonical period map $\varphi : S \to \Gamma \backslash D$ sending the point $s$ to the isomorphism class of the polarized Hodge structre at $s$. Applying \autoref{axshanres}(ii) we have that $\overline{\psi(B \cap Z)}^{\textrm{Zar}} = \mathbf{H}_{Z} \cdot \psi(s)$, and applying \autoref{weakspredef} the irreducible component $Y$ of $\varphi^{-1}(\pi(\mathbf{H}_{Z} \cdot \psi(s)))$ containing $Z$ is a weakly special subvariety of type $\mathcal{C}_{i}$; here $\pi : D \to \Gamma \backslash D$ is the projection. Necessarily $Y$ is of larger dimension, and so by the compatibility of $\psi$ with $\varphi$ the germ of $\psi^{-1}(\mathbf{H}_{Z} \cdot \psi(s))$ at $s$ has dimension at least $d+1$ and contains the corresponding germ of $Z$. The containment moreover continues to hold if we consider the $r$'th order neighbourhoods of these germs, and moreover by \autoref{angermlarge} the dimension of the coordinate ring of the $r$th order neighbourhood of the germ of $\psi^{-1}(\mathbf{H}_{Z} \cdot \psi(s))$ at $s$ is at least $H(r,d+1)$. The fact that $\mathcal{C}_{i} = \mathcal{C}(Z)$ means that there exists $g \in \GL_{m}(\mathbb{C})$ such that $g \cdot (\mathbf{H}_{Z} \cdot \psi(s)) = V_{i}$, hence $\psi^{-1}(\mathbf{H}_{Z} \cdot \psi(s)) = (g \cdot \psi)^{-1}(V_{i})$. It follows that $t$ does not lie in $\mathcal{U}(g, V_{i}, r)$ for every $r$, and so does not lie in the union.

Conversely, let us suppose that $Z = p(\mathcal{Z}_{t})$ has type $\mathcal{C}_{i}$, but that $t$ does not lie in the specified union; in particular, we assume that $t \notin \mathcal{U}(g, V_{i}, r)$ for every $r$. Then by definition of $\mathcal{U}(g, V_{i}, r)$, we find that the constructible algebraic subset $\mathcal{G}_{r} \subset \mathcal{B}(g, V_{i}, r)$ consisting of points $(t, s, M)$ such that $s$ lies in $p(\mathcal{Z}_{t})(\mathbb{C})$ is non-empty for every $r$. We observe that $\mathcal{G}_{r+1} \subset \mathcal{G}_{r}$ for each $r$, hence by constructibility the countable intersection $\bigcap_{r \geq 0} \mathcal{G}_{r}$ of non-empty constructible sets is non-empty containing an element $(t, s_{0}, M_{0})$. If we then let $\psi$ be the local period map defined by $M_{0}$, then we claim that the germ of $\psi^{-1}(V_{i})$ at $s_{0}$ has dimension at least $d+1$: indeed, by \cite[\href{https://stacks.math.columbia.edu/tag/00KQ}{Proposition 00KQ}]{stacks-project} the map $r' \mapsto \dim_{\mathbb{C}} \widehat{\mathcal{O}}_{\psi^{-1}(V), s} / \mathfrak{m}^{r'+1}_{\psi^{-1}(V), s}$ is a polynomial in $r'$ of degree equal to the dimension of $\psi^{-1}(V)$, and $H(r',d+1)$ is a polynomial of degree $d+1$ in $r'$. By construction, the germ of $Z$ at $s$ lands inside $V_{i}$, hence $\overline{\psi(B \cap Z)}^{\textrm{Zar}} \subset V_{i}$. Since $\mathcal{C}(Z) = \mathcal{C}_{i}$, we in fact must have the equality $\overline{\psi(B \cap Z)}^{\textrm{Zar}} = V_{i}$ as the two irreducible algebraic loci have the same dimension. Choosing $g \in \GL_{m}(\mathbb{C})$ such that $g \cdot \psi$ is a local lift of the period map $\varphi : S \to \Gamma \backslash D$ sending a point $s$ to its associated polarized Hodge structure, we find by \autoref{axshanres}(ii) that $g \cdot V_{i}$ is identified with an orbit $\mathbf{H}_{Z} \cdot (g \cdot \psi)(s_{0})$ of the algebraic monodromy group of $Z$. We thus find that the weakly special subvariety $Y$ which is a component of $\varphi^{-1}(\pi(g \cdot V_{i}))$ and contains $Z$ has dimension at least $d+1$, hence the result.
\end{proof}

Let us now return to the case of a general base $S$ (i.e., we do not assume that $\Omega^{1}_{S}$, the Hodge bundles $F^{\bullet} \mathcal{H}$, etc., are trivial). The main result of this section is then the following:

\begin{prop}
\label{countunionprop}
Let $g : \mathcal{Z} \to T$ be a $K$-algebraic family of irreducible $d$-dimensional subvarieties of $S$, all of type $\mathcal{C}_{i}$. Then the locus in $T$ of maximal algebraic varieties of type $\mathcal{C}_{i}$ (i.e., the weakly special locus) is a countable union of constructible $K$-algebraic sets.
\end{prop}

\begin{proof}
Given a $K$-algebraic open subset $U \subset S$, we may reduce to showing the same statement for just those fibres of $g$ which intersect $U$, noting that an irreducible $d$-dimensional complex variety $Z$ of $S$ which intersects $U$ is weakly special if and only if its intersection with $U$ is weakly special for the variation $\restr{\mathbb{V}}{U}$. We may then find a finite cover of $S$ by open sets $U$ where $\Omega^{1}_{U}, F^{\bullet} \mathcal{H}$, etc., are all trivial. The result then follows from \autoref{maxcondloc}.
\end{proof}

\subsection{Constructible Loci Associated To Families}

In preparation for our main arguments, we will need to know that two additional sorts of conditions associated to families of varieties are described by constructible algebraic conditions on their moduli. We recall we have fixed a projective compactification $S \subset \overline{S}$, and that $\textrm{Hilb}(\overline{S})$ is the assocaited Hilbert scheme.

\begin{lem}
\label{constgivesconstrcond}
Let $\mathcal{T}^{d}_{r} \subset J^{d}_{r} S$ be a constructible $K$-algebraic set. Then the condition that $J^{d}_{r} Z$ (or $J^{d}_{r} Z_{\textrm{sm}}$) lie inside $\mathcal{T}^{d}_{r}$, where $Z = \overline{Z} \cap S$, gives a constructible $K$-algebraic condition on $[\overline{Z}] \in \textrm{Hilb}(\overline{S})$ within each component of $\textrm{Hilb}(\overline{S})$.
\end{lem}

\begin{proof}
We will handle the case where $Z$ satisfies the condition $\mathcal{T}^{d}_{r}$, the case for the smooth locus $Z_{\textrm{sm}}$ being entirely analogous.

Fix a component $C$ of $\textrm{Hilb}(\overline{S})$, and denote by $\overline{u} : \overline{\mathcal{Y}} \to C$ the universal family over $C$; we note that we may view $\overline{\mathcal{Y}}$ as a closed subscheme of $\overline{S} \times C$, and let $p : \overline{\mathcal{Y}} \to \overline{S}$ be the natural projection. We let $\mathcal{Y} = p^{-1}(S)$, and let $u : \mathcal{Y} \to C$ be the natural restriction. We may furthermore construct a family $u^{d}_{r} : \mathcal{Y}^{d}_{r} \to C$ such that the fibre above $[\overline{Z}]$ may be identified with $J^{d}_{r} Z$, where $Z = \overline{Z} \cap S$; we do this by considering the subscheme of $J^{d}_{r} \mathcal{Y}$ whose jets are constant on $C$. Then the projection $p$ admits a natural extension $p^{d}_{r} : \mathcal{Y}^{d}_{r} \to J^{d}_{r} S$, the set $p^{d,-1}_{r}(\mathcal{T}^{d}_{r})$ is $K$-algebraically constructible, and the desired locus in $C$ is the image of the complement $\mathcal{Y}^{d}_{r} \setminus p^{d,-1}_{r}(\mathcal{T}^{d}_{r})$ under $u^{d}_{r}$. This is a constructible $K$-algebraic set by Chevalley's Theorem.
\end{proof}

\begin{prop}
\label{liesinotherfam}
Suppose that $g_{i} : \mathcal{Y}_{i} \to T_{i}$ for $i = 1,2$ are two $K$-algebraic families of geometrically irreducible subvarieties of $S$, and suppose that all the fibres of $g_{1}$ have dimension $d$. Then the locus of $t \in T_{1}$ where the fibre $\mathcal{Y}_{1,t}$ lies inside a fibre of $g_{2}$ is a constructible $K$-algebraic locus.
\end{prop}

\begin{proof}
The locus in question may be constructed by considering the fibre product $\mathcal{Y}_{12} = \mathcal{Y}_{1} \times_{S} \mathcal{Y}_{2}$ with respect to the natural projections $p_{i} : \mathcal{Y}_{i} \to S$ for $i = 1, 2$, and considering the family $g_{12} : \mathcal{Y}_{12} \to T_{1} \times T_{2}$. The fibre above $(t_{1}, t_{2})$ of $g_{12}$ may be identified with the intersection $\mathcal{Y}_{1,t_{1}} \cap \mathcal{Y}_{1,t_{2}}$ of fibres. By geometric irreducibility, it suffices to compute the $K$-constructible locus in $T_{1} \times_{K} T_{2}$ where the fibres have dimension $d$ and project this to $T_{1}$, hence the result.
\end{proof}

\subsection{Digression on Definability}

In the next section we will need some tools from o-minimal geometry; in this short section we review the necessary background. For details see for instance \cite{defpermap}.

By a structure $\mathcal{S}$ we mean, for each $k$, a certain collection of subsets $S_{k} \subset \mathcal{P}(\mathbb{R}^{k})$ of the power set of $\mathbb{R}^{k}$, closed under finite unions, intersections, products, and projections; and such that $\mathbb{R} \in \mathcal{S}_{1}$, the diagonal in $\mathbb{R}^{2}$ is in $\mathcal{S}_{2}$, and the graph of addition function $+ : \mathbb{R}^{2} \to \mathbb{R}$ is in $\mathcal{S}_{3}$. We call the sets which are elements of some $\mathcal{S}_{k}$ definable sets. Given a function $f : A \to B$ between two definable subsets, we say that $f$ is definable if its graph is a definable subset. The structure of interest in Hodge theory is $\mathcal{S} = \mathbb{R}_{\textrm{an}, \textrm{exp}}$. It includes as definable any real analytic function on a compact interval in $\mathbb{R}$, as well as the graph of the real exponential. Using the identification $\mathbb{C} = \mathbb{R}^{2}$ one may also speak of definable complex analytic functions. This lets us define a definable analytic variety $X$ as an analytic variety admitting a finite cover $\{ U_{i} \}_{i = 1}^{n}$ by definable subsets $U_{i} \subset \mathbb{C}^{n}$, such that $X \cap U_{i}$ is cut out by definable analytic functions, and such that the transition maps are likewise definable analytic. Algebraic varieties have a canonical definable analytic structure, as one can see by compactifying. Using the fact that derivatives of definable analytic functions are definable analytic, one easily checks (for instance using the explicit models for jet spaces described in \cite[\S2]{periodimages}) that $J^{d}_{r}$ is naturally a functor on the category of definable analytic spaces.

We endow $D \subset \ch{D}$ with the definable structure induced from the algebraic variety $\ch{D}$. The main application of definable analytic spaces we will need is then the following:

\begin{prop}
\label{deflifts}
For any algebraic subvariety $Y \subset S$, there exists a finite cover $B_{1}, \hdots, B_{k}$ of $Y$ by definable open sets and definable local lifts $\psi_{i} : B_{i} \to D$ of the period map $Y \to \Gamma \backslash D$ obtained by restricting $\varphi : S \to \Gamma \backslash D$.
\end{prop}

\begin{proof}
As we may always restrict local lifts, it suffices to show this in the case where $S = Y$. We may reduce to the same statement with $S$ replaced by a finite \'etale cover, so we may in particular assume that the monodromy $\Gamma_{S}$ lies in a neat arithmetic subgroup $\Gamma' \subset G(\mathbb{Z})$, and is unipotent around a normal crossing divisor $E = X \setminus S$ corresponding to some log-smooth compactification $S \subset X$. Construct an analytic period map $\varphi' : S \to \Gamma' \backslash D$. Arguing as in \cite[\S3]{AXSCHAN}, we may cover $S$ by finitely many definable neighbourhoods of the form $\Delta^{k} \times (\Delta^{*})^{\ell}$ for some $k$ and $\ell$, where $\Delta$ is the open unit disc in $\mathbb{C}$ and $\Delta^{*}$ is the disk punctured at zero. Moreover, as in \cite[\S3]{AXSCHAN} each such neighbourhood fits into a definable analytic diagram
\begin{center}
\begin{tikzcd}
\Delta^{k} \times \Sigma^{\ell} \arrow[d, "\textrm{exp}"] \arrow[r, "\widetilde{\varphi'}"] & D \\
\Delta^{k} \times (\Delta^{*})^{\ell}
\end{tikzcd}
\end{center}
where the vertical map is a restriction of the universal covering $\Delta^{k} \times \mathbb{H}^{\ell} \to \Delta^{k} \times (\Delta^{*})^{\ell}$ and the horizontal map is a lift of $\varphi'$. It now suffices to take a finite covering of each neighbourhood $\Delta^{k} \times (\Delta^{*})^{\ell}$ by simply-connected definable open sets and locally invert $\exp$.
\end{proof}

Lastly, we recall the Peterzil-Starchenko theorem, proven in \cite{peterzil2008complex}, which tells us that a definable analytic subset of an algebraic variety is in fact algebraic.

\subsection{Rigidity and Families of Weakly Specials}

Recall we have fixed a projective compactification $S \subset \overline{S}$. Let $\textrm{Hilb}(\overline{S})$ be the associated Hilbert scheme, and denote by $\textrm{Var}(S)$ the locus consisting of those geometrically integral subschemes which intersect $S$; this is a $K$-algebraic constructible locus. Using the properness and flatness of the universal family (see [EGA IV, 12.2.1]), one checks this locus is in fact open, so we may endow it with the induced reduced scheme structure. To understand the locus $\mathcal{W} \subset \textrm{Var}(S)$ of weakly special subvarieties, we will need some preliminary results about families of weakly special subvarieties in $S$ and the conditions they induce on $\textrm{Var}(S)$.

The first notion we will need is the following, which is referred to as a \emph{weakly non-factor} weakly special subvariety in \cite{fieldsofdef}:

\begin{defn}
Given a weakly special subvariety $Z \subset S$, we say that $Z$ is \emph{rigid} if there does not exist a weakly special variety $Y$ with $Z \subsetneq Y \subset S$ such that $\mathbf{H}_{Z}$ is normal in $\mathbf{H}_{Y}$.
\end{defn}

\noindent As explained in \cite[\S2]{fieldsofdef}, the set of rigid weakly special subvarieties is countable, and corresponds to those weakly special subvarieties which cannot be ``Hodge-theoretically deformed'' inside $S$.

\begin{lem}
\label{existsbiggerrigid}
Suppose that $Z \subset S$ is a weakly special. Then there exists a rigid weakly special $Y \subset S$ containing $Z$ such that $\mathbf{H}_{Z}$ is normal in $\mathbf{H}_{Y}$.
\end{lem}

\begin{proof}
From the definition it follows that there exists a finite chain $Y_{0} \subsetneq Y_{1} \subsetneq \cdots \subsetneq Y_{k}$ where each $Y_{i}$ is weakly special, $Z = Y_{0}$, $\mathbf{H}_{Y_{i}}$ is normal in $\mathbf{H}_{Y_{i+1}}$, and $Y_{k}$ is rigid. Since each algebraic monodromy group is semisimple, it follows that $\mathbf{H}_{Z} = \mathbf{H}_{Y_{0}}$ is a product of simple factors of $\mathbf{H}_{Y_{k}}$.
\end{proof}

Fix a geometrically irreducible variety $Y \subset S$. Then for each point $s \in Y$, the algebraic monodromy group is well-defined inside the automorphism group $\mathbf{H}_{Y, s} \subset \GL(\mathbb{V}_{s})$ of each fibre. Given a path $\gamma$ in $Y$ from $s$ to $s'$, we obtain isomorphisms $\alpha_{\gamma} : \mathbb{V}_{s} \xrightarrow{\sim} \mathbb{V}_{s'}$. By writing the path $\gamma$ as a composition $\gamma = \gamma_{k} \circ \cdots \circ \gamma_{1}$, where each $\gamma_{i}$ can be lifted to the normalization $Y^{\textrm{nor}} \to Y$, one may check that $\mathbf{H}_{Y, s'} = \alpha_{\gamma} \circ \mathbf{H}_{Y,s} \circ \alpha_{\gamma}^{-1}$. However, given a $\mathbb{Q}$-algebraic normal subgroup $N \subset \mathbf{H}_{Y,s}$, it is not clear that $N' = \alpha_{\gamma} \circ N \circ \alpha_{\gamma}^{-1}$ is independent of $\gamma$. Different choices of $\gamma$, however, will lead to normal $\mathbb{Q}$-algebraic groups $N' \subset \mathbf{H}_{Y,s'}$ conjugate under an element $\sigma \in \Gamma_{Y} \subset \GL(\mathbb{V}_{s'})$, where $\Gamma_{Y}$ is the monodromy on $Y$. The element $\sigma$ acts through the outer automorphism group of $\mathbf{H}_{Y,s'}$, giving finitely many possibilities $N'_{1}, \hdots, N'_{\ell}$ for the normal $\mathbb{Q}$-algebraic subgroup $N'$. Moreover, if $\psi : B \to D$ is a local lift of the period map on $Y$ and $s' \in B$, the orbits $N'_{i} \cdot \psi(s')$ for $i = 1, \hdots, \ell$ all define the same type, and this type is independent of $\psi$. To see this, note that if $\psi$ is replaced by $\sigma \cdot \psi$ and $N'$ by $\sigma \cdot N' \cdot \sigma^{-1}$ then the resulting orbit is $\sigma \cdot (N' \cdot \psi(s'))$. The independence of the neighbourhood $B$ follows by analytic continuation. We denote this type by $\mathcal{C}(Y, N)$.

\begin{prop}
\label{closedfamlocusprop}
Let $Y \subset S$ be a geometrically irreducible variety, and let $N \subset \mathbf{H}_{Y}$ be a connected normal $\mathbb{Q}$-algebraic subgroup. Consider the locus $\mathcal{D}(Y, N) \subset \textrm{Var}(S)$ consisting of $[Z]$ such that $Z \subset Y$ and such that $\mathbf{H}_{Z}$ lies inside one of the parallel translates $N_{1}, \hdots, N_{\ell}$ of $N$; the reasoning in the preceding paragraph shows this locus is well-defined.
Then 
\begin{enumerate}
\item[(i)] $\mathcal{D}(Y, N)$ is a closed algebraic locus in $\textrm{Var}(S)$; and
\item[(ii)] for every $[Z] \in \mathcal{D}(Y, N)$ we have $\mathcal{C}(Z) \leq \mathcal{C}(Y, N)$.
\end{enumerate}
\end{prop}

\begin{proof}
By \autoref{deflifts} we may choose finitely many definable analytic local lifts $\psi_{i} : B_{i} \to D$ for $1 \leq i \leq n$ of the restriction of $\varphi$ to $Y$, and after translating we may assume that if $y \in B_{i} \cap B_{i'}$ then $\psi_{i}$ and $\psi_{i'}$ differ on $B_{i} \cap B_{i'}$ by an element of the monodromy group $\Gamma_{Y} \subset \GL(\mathbb{V}_{y})$. More precisely, we may assume that the isomorphisms $\mathbb{V}_{y} \simeq V$ which determine the $\psi_{i}$ differ on overlaps by an element of $\Gamma_{Y}$; here we recall the notation in \autoref{weakspprelimsec}. This implies that the maps $\{ \psi_{i} \}_{i = 1}^{n}$ determine a common embedding $\mathbf{H}_{Y} \subset G$ (recall the discussion preceding \autoref{axshanres}), and hence applying \autoref{axshanres}(ii) the images $\psi_{i}(B_{i})$ for $1 \leq i \leq n$ lie in a common orbit $\ch{D}_{Y} = \mathbf{H}_{Y} \cdot t \subset \ch{D}$. For each of the realizations $N_{1}, \hdots, N_{\ell}$ of the group $N$, we obtain a semisimple splitting $\mathbf{H}_{Y} = N_{j} \cdot N'_{j}$, and corresponding factorizations $\ch{D}_{Y} = \ch{D}_{j} \times \ch{D}'_{j}$ by orbits. We define $\mathcal{K}^{d}_{r,j} \subset J^{d}_{r} \ch{D}_{Y}$ to be the locus of jets constant on $\ch{D}'_{j}$, and define $\mathcal{J}^{d}_{r} \subset J^{d}_{r} Y$ by $\mathcal{J}^{d}_{r} \cap J^{d}_{r} B_{i} = (J^{d}_{r} \psi_{i})^{-1}(\mathcal{K}^{d}_{r,1} \cup \cdots \cup \mathcal{K}^{d}_{r,\ell})$. 

From the fact that $\Gamma_{Y}$ preserves the union $\mathcal{K}^{d}_{r,1} \cup \cdots \cup \mathcal{K}^{d}_{r,\ell}$, it is clear that the resulting locus is well-defined and closed analytic in $J^{d}_{r} Y$. As the maps $J^{d}_{r} \psi_{i}$ are definable analytic, this locus is also definable, so by the Peterzil-Starchenko Theorem \cite{peterzil2008complex}, algebraic. Thus by \autoref{constgivesconstrcond} above the locus $\mathcal{N}^{d}_{r} \subset \textrm{Var}(S)$ of points $[Z]$ which satisfy the differential conditions $\mathcal{J}^{d}_{r}$ is algebraically constructible in each component of $\textrm{Var}(S)$. It is clear that since $\mathcal{J}^{d}_{r} \subset J^{d}_{r} Y$ for every $d$ and $r$, any such $[Z]$ which satisfies $\mathcal{J}^{d}_{r}$ must lie in $Y$. Let $C \subset B_{i} \cap Z$ be an analytic component for such a $Z$ and pick a smooth point $z \in C$. Then the irreducibility of the fibre $(J^{d}_{r} C)_{z}$ implies there is a $j$ such that $(J^{d}_{r} \psi_{i})((J^{d}_{r} C)_{z}) \subset \mathcal{K}^{d}_{r,j}$. The constancy of the jets on the factor $\ch{D}'_{j}$ then implies that $\psi(C) \subset \ch{D}_{j} \times \{ t_{2} \}$ for some $t_{2} \in \ch{D}'_{j}$, and hence $\psi(C)$ lies in an $N_{j}$ orbit; it follows that $\mathcal{C}(Z) \leq \mathcal{C}(Y, N)$. Using \autoref{axshanres}(ii), \autoref{monocont} and \autoref{nocomprealfac} we moreover learn that $\mathbf{H}_{Z} \subset N_{j}$. Finally, for any $Z$ such that $\mathbf{H}_{Z} \subset N_{j}$ for some $j$, it is clear from the definitions and \autoref{axshanres}(ii) that $Z$ satisfies $\mathcal{J}^{d}_{r}$ for all $d$ and $r$. We may thus identify $\mathcal{D}(Y,N)$ with the intersection $\bigcap_{d, r} \mathcal{N}^{d}_{r}$.

It now suffices to show the countable intersection $\mathcal{D}(Y,N)$ of algebraically constructible sets is closed. Consider a sequence of points $[Z_{j}] \in \mathcal{D}(Y, N)$ converging to a limit $[Z] \in \textrm{Var}(S)$. It suffices to argue that $Z$ satisfies the constraints $\mathcal{J}^{d}_{r}$ for all $d$ and $r$. Reasoning as in the preceding paragraph, it suffices to check an irreducible component $C$ of an analytic germ of $Z$ satisfies $\mathcal{J}^{d}_{r}$ for all $d$ and $r$. This is equivalent to showing that after shrinking $C$ so that it lies in some $B_{i}$, the image $\psi_{i}(C)$ lies in an orbit of one of the groups $N_{1}, \hdots, N_{\ell}$. Applying \autoref{landinlem} it will suffice to take $d = \dim Z$ and show the locus $\mathcal{J}^{d}_{r}$ contains a single compatible family $\sigma_{r} \in J^{d}_{r,nd} C$ of non-degenerate jets with $\sigma_{0}$ lying in the smooth locus of $Z$; let us fix such a sequence. After passing to a subsequence we may assume that the points $[Z_{j}]$ lie in a single component $T$ of $\textrm{Var}(S)$, and let $g : \mathcal{Y} \to T$ be the universal family over $T$ (i.e., the restriction of the universal family over $\textrm{Hilb}(\overline{S})$ obtained by replacing the fibres with their intersections with $S$). After replacing $\mathcal{Y}$ with the open locus $\mathcal{Y}_{\textrm{sm}}$ which parametrizes the smooth loci of the fibres of $\mathcal{Y}$, the restriction $\mathcal{Y}_{\textrm{sm}} \to T$ is flat, of finite presentation, and has smooth fibres: it is a smooth morphism. It then follows from \autoref{jetconvergence} that for each $r$ we can lift the sequence $[Z_{j}]$ to a sequence of jets $\sigma_{j,r} \in J^{d}_{r} Z_{j}$ which converge to $\sigma_{r}$. We conclude that $\sigma_{r} \in \mathcal{J}^{d}_{r}$ for all $r$ from the fact that $\mathcal{J}^{d}_{r}$ is closed.
\end{proof}

\subsection{The Locus $\mathcal{W}$}
\label{weaksplocsec}

We are now ready to study the locus $\mathcal{W} \subset \textrm{Var}(S)$ of weakly special subvarieties. To aid in the discussion, let us denote for any type $\mathcal{C}$ the locus $\mathcal{W}(\mathcal{C}) \subset \textrm{Var}(S)$ of varieties that satisfy the differential constraints $\mathcal{T}^{d}_{r}(\mathcal{C})$ for all $d$ and $r$. From \autoref{liesintrans} it is clear that this is the same as the locus of $[Z] \in \textrm{Var}(S)$ such that $\mathcal{C}(Z) \leq \mathcal{C}$. We write $\mathcal{C} = \mathcal{C}(V)$ for some variety $V$ over the field $L$.



\begin{prop}
\label{typeinequalityloci}
The locus $\mathcal{W}(\mathcal{C}) \subset \textrm{Var}(S)$ is, within each component of $\textrm{Var}(S)$, a closed $L$-algebraic subvariety. Moreover, each geometric component $C$ of $\mathcal{W}(\mathcal{C})$ agrees with a component of $\mathcal{D}(Y, N)$ for some choice of rigid weakly special $Y$ and algebraic monodromy group $N = \mathbf{H}_{Z_{0}}$, where $[Z_{0}] \in C$ is some point.
\end{prop}

\begin{proof}
We fix an $L$-component $T \subset \textrm{Var}(S)$ and argue about the locus $\mathcal{W}(\mathcal{C}) \cap T$. The condition on $[Z] \in \textrm{Var}(S)$ that $Z$ satisfies the constraint $\mathcal{T}^{d}_{r}(\mathcal{C})$ is an $L$-algebraically constructable condition on $T$ by \autoref{constgivesconstrcond}, so the condition that we satisfy all the differential constraints exhibits $\mathcal{W}(\mathcal{C}) \cap T$ as a countable intersection of $L$-algebraically constructable sets. On the other hand, we claim that $\mathcal{W}(\mathcal{C})$ is also a countable union of closed algebraic sets $\mathcal{D}(Y_{j}, N_{j})$ for choices of rigid $Y_{j}$ and subgroups $N_{j} \subset \mathbf{H}_{Y}$ coming from the algebraic monodromy of points of $\mathcal{W}(\mathcal{C})$, which will suffice to show the result. (See also \autoref{unionintlem} below.)

To see this, consider a point $[Z_{0}] \in \mathcal{W}(\mathcal{C})$. By \autoref{existsbiggerrigid}, we can find a rigid weakly special $Y$ such that $\mathbf{H}_{Z_{0}}$ is a $\mathbb{Q}$-algebraic normal subgroup of $\mathbf{H}_{Y}$. By \autoref{closedfamlocusprop} we have a closed algebraic locus $\mathcal{D}(Y,\mathbf{H}_{Z_{0}})$, and the points $[Z]$ in this locus satisfy the condition $\mathcal{C}(Z) \leq \mathcal{C}(Y,\mathbf{H}_{Z_{0}}) = \mathcal{C}(Z_{0})$. Since $[Z_{0}] \in \mathcal{W}(\mathcal{C})$ implies $\mathcal{C}(Z_{0}) \leq \mathcal{C}$, it follows that $\mathcal{D}(Y, \mathbf{H}_{Z_{0}}) \subset \mathcal{W}(\mathcal{C})$. On the other hand, as we vary the point $[Z_{0}]$, only countably many rigid weakly specials $Y$ and $\mathbb{Q}$-algebraic groups $\mathbf{H}_{Z_{0}}$ can arise, hence $\mathcal{W}(\mathcal{C}) = \bigcup_{j = 1}^{\infty} \mathcal{D}(Y_{j}, N_{j})$ for some countable collection $\{ (Y_{j}, N_{j}) \}_{j = 1}^{\infty}$ as described.
\end{proof}

\begin{lem}
\label{unionintlem}
Let $X$ be a $L$-algebraic variety, and suppose that a set $\mathcal{S} \subset X(\mathbb{C})$ is both a countable intersection $\mathcal{S} = \bigcap_{i = 1}^{\infty} D_{i}$ of constructible $L$-algebraic sets $D_{i}$ and a countable union $\mathcal{S} = \bigcup_{i = 1}^{\infty} E_{i}$ of constructible $L$-algebraic sets $E_{i}$. Then $\mathcal{S}$ is $L$-constructible, and there exists $m$ and $n$ such that
\[ \mathcal{S} = \bigcap_{i = 1}^{m} D_{i} = \bigcup_{i = 1}^{n} E_{i} . \]
\end{lem}

\begin{proof}
By replacing $D_{i}$ with $D_{1} \cap \cdots \cap D_{i}$ and $E_{i}$ with $E_{1} \cup \cdots \cup E_{i}$ we may assume that $D_{i} \supset D_{i+1}$ and $E_{i} \subset E_{i+1}$ for all $i$. We recall that any constructible $L$-algebraic subset of $X$ is of the form $(V_{1} \cap U_{1}) \cup \cdots \cup (V_{k} \cap U_{k})$ for closed $L$-algebraic sets $V_{1}, \hdots, V_{k}$ and open $L$-algebraic sets $U_{1}, \hdots, U_{k}$. The sequence of closures $\overline{D_{1}}, \hdots, \overline{D_{i}}, \hdots$ therefore eventually stabilizes to some $L$-algebraic union $V_{1} \cup \cdots \cup V_{k}$ of closed $L$-algebraic sets, and replacing the sequence $\{ D_{i} \}_{i = 1}^{\infty}$ by a subsequence we may assume that each $D_{i}$ is of the form $(V_{1} \cap U_{1,i}) \cup \cdots \cup (V_{k} \cap U_{k,i})$.

Using the fact that a countable intersection of Zariski open subsets of an irreducible closed set remains dense, the closure of the set $\mathcal{S}$ is equal to $V_{1} \cup \cdots \cup V_{k}$. It follows from the fact that closed complex algebraic sets form meagre subsets of algebraic sets of larger dimension that for large enough $i$ the sets $E_{i}$ are also of the form $(V_{1} \cap U'_{1,i}) \cup \cdots \cup (V_{k} \cap U_{k,i})$. The complements $V_{j} \setminus (V_{j} \cap U'_{1,i})$ are a sequence of closed subsets of $V_{j}$, decreasing in the index $i$, so must stabilize for large $i$. The result follows.
\end{proof}

\begin{thm}
\label{Wisconstrthm}
Suppose that $\mathbb{V}$ is a polarized integral variation of Hodge structure on the smooth $K$-algebraic variety $S$, and that the associated algebraic data admits a model $(\mathcal{H}, F^{\bullet}, \nabla)$ over $K$. Then the locus $\mathcal{W} \subset \textrm{Var}(S)$ of weakly special varieties is a constructible $K$-algebraic subset within each component of $\textrm{Var}(S)$.
\end{thm}

\begin{proof}
It suffices to show that $\mathcal{W}$ is a constructible complex algebraic set in each component; the $K$-algebraicity follows from \cite[Prop 3.1(b)]{fieldsofdef}, which shows that an $\textrm{Aut}(\mathbb{C}/K)$-conjugate of a weakly special variety is again weakly special.

Recall we have a fixed list of types $\mathcal{C}_{1}, \hdots, \mathcal{C}_{q}$ arising from the variation $\mathbb{V}$. It suffices to show this statement for the locus $\mathcal{W}_{j} \subset \mathcal{W}(\mathcal{C}_{j})$ of varieties which are weakly special (i.e., maximal) for the differential constraints $\mathcal{T}^{d}_{r}(\mathcal{C}_{j})$; we write $\mathcal{C}_{j} = \mathcal{C}(V_{j})$ for some $V_{j} \subset \ch{L}$. Thus we need to compute the sublocus of $\mathcal{W}(\mathcal{C}_{j})$ consisting of points $[Z] \in \mathcal{W}(\mathcal{C}_{j})$ such that for any $[Y] \in \mathcal{W}(\mathcal{C}_{j})$ with $Z \subset Y$ we have $Z = Y$. Let us denote by $\{ C_{i} \}_{i = 1}^{\infty}$ the components of $\mathcal{W}(\mathcal{C}_{j})$, and fix some $i_{0}$. We denote by $C_{i_{0}}^{\circ} \subset C_{i_{0}}$ the subset consiting of $[Z]$ such that $\mathcal{C}(Z) = \mathcal{C}_{j}$; this is the open locus obtained by removing the intersections of $C_{i_{0}}$ with $\mathcal{W}(\mathcal{C}_{j'})$ for those $j' \neq j$ such that $\mathcal{C}_{j'} \leq \mathcal{C}_{j}$. We let $d$ be the dimension of the subvarieties parametrized by $C_{i_{0}}^{\circ}$. We may compute the maximal points of $C_{i_{0}}^{\circ}$ in two ways:
\begin{itemize}
\item[(i)] on one hand a point $[Z] \in C_{i_{0}}^{\circ}$ is maximal if it lies outside the loci
\[ E_{i} = \{ [Z] \in C_{i_{0}}^{\circ} : \textrm{there exists } [Y] \in C_{i} \textrm{ such that } Z \subset Y \}, \]
for every $i \neq i_{0}$;
\item[(ii)] on the other hand for each $r$ we have from \autoref{countunionprop} a locus $\mathcal{U}(C_{i_{0}}^{\circ}, V_{j}, r) \subset C_{i_{0}}^{\circ}$ such that $[Z] \in C_{i_{0}}^{\circ}$ is maximal if and only if it lies inside $\mathcal{U}(C_{i_{0}}^{\circ}, V_{j}, r)$ for some $r$.
\end{itemize}
Both the conditions described in (i) and (ii) are constructible; for (i) this is true by \autoref{liesinotherfam}, and for (ii) this is true by \autoref{countunionprop}. It follows that the locus $\mathcal{W}_{j} \cap C_{i_{0}}^{\circ} = \mathcal{W}_{j} \cap C_{i_{0}}$ is both a countable union and a countable intersection of constructible sets, hence constructible by \autoref{unionintlem}.
\end{proof}

\section{The Computation of $\mathcal{W}$}
\label{compsec}

Thus far in \autoref{modulisec} we have ignored computational considerations. We now explain why the components of the moduli scheme $\textrm{Var}(S)$ as well as the constructible conditions defining $\mathcal{W}$ in each component are computable. We fix a closed embedding $\overline{S} \subset \mathbb{P}^{n}$ with respect to which we compute degrees of subvarieties of $\overline{S}$. We will assume throughout that we have an open $K$-algebraic affine cover $S = \bigcup_{i = 1}^{n} S_{i}$ and models $(\mathcal{H}_{i}, F^{\bullet}_{i}, \nabla_{i})$ for the $K$-algebraic data associated to $\restr{\mathbb{V}}{S_{i}}$, where $\mathbb{V}$ is a polarizable integral variation of Hodge structure $\mathbb{V}$, and $K \subset \mathbb{C}$ is a computable subfield.

\subsection{Loci Defined by Families}

\begin{prop}
\label{familyalgos}
Let $f : X \to T$ be a projective $K$-algebraic morphism over the Noetherian base $T$. Then there exists algorithms to compute:
\begin{itemize}
\item[(i)] a \emph{flattening stratification} of $T$, i.e. a partition $T = T_{1} \cup \cdots \cup T_{k}$ into locally closed $K$-algebraic sets such that the base-change maps $f : X_{T_{i}} \to T_{i}$ are flat for $i = 1, \hdots, k$;
\item[(ii)] a partition $T = T_{1} \cup \cdots \cup T_{k}$ of $T$ into constructible sets such that the geometrically irreducible components of the fibres are given uniformly (see \cite[Def. 1]{ayad2010computing}) in each set (in particular, some subset of the sets $\{ T_{i} \}_{i = 1}^{k}$ parametrize exactly the geometrically irreducible fibres);
\item[(iii)] the constructible locus $T_{\textrm{red}} \subset T$ where the fibres of $f$ are reduced;
\end{itemize}
\end{prop}

\begin{proof}
~
\begin{itemize}
\item[(i)] The theory of parametric Gr\"obner bases solves this problem, as is explained by Swinarski in \cite[\S1.4]{swinarski2016equations}. This is also discussed in the paper of \cite{wibmer2007grobner}; see in particular \cite[Cor. 2]{wibmer2007grobner}.
\item[(ii)] This problem is solved by Ayad in \cite{ayad2010computing}. Note that Ayad works over $\mathbb{Q}$, but explains how to extend the methods to a general $K$ in his PhD thesis \cite{ayad2006complexite}.
\item[(iii)] By (i) we may reduce to the situation where $f$ is flat, and by taking an irreducible decomposition of $T$, where $T$ is irreducible and reduced. In this situation the locus $T_{\textrm{red}} \subset T$ is an open subset, either empty or dense. It will suffice to decide whether $T_{\textrm{red}}$ is non-empty, and if it is, compute a non-empty open subset $V \subset T_{\textrm{red}}$ and reduce to the same problem for the family $X_{T \setminus V} \to T \setminus V$, which we can resolve by induction on the dimension of the base.

It suffices to make the argument in \cite[\href{https://stacks.math.columbia.edu/tag/0578}{Lemma 0578}]{stacks-project} constructive. More specifically, \cite[\href{https://stacks.math.columbia.edu/tag/0578}{Lemma 0578}]{stacks-project} in combination with the openness of the reduced locus tells us $T_{\textrm{red}}$ is non-empty if and only if $X_{K(T)}$ is reduced, where $K(T)$ is the function field of $T$. We can decide if $X_{K(T)}$ is reduced using any of the algorithms for computing the reduction of a finite-type scheme over a field (e.g., using \cite{krick1991algorithm}), so it suffices to compute $V$ in the case where $T_{\textrm{red}}$ is non-empty. 

For this we can follow the method of proof in \cite[\href{https://stacks.math.columbia.edu/tag/0578}{Lemma 0578}]{stacks-project} and  \cite[\href{https://stacks.math.columbia.edu/tag/0573}{Lemma 0573}]{stacks-project}, which computes the locus $X_{\textrm{sm}} \subset X$ consisting of $x \in X$ such that the germ $(X_{f(x)}, x)$ is smooth, and then argues we can take for $V$ an open locus consisting of $t \in T$ where $X_{\textrm{sm}, t}$ is scheme-theoretically dense in $X_{t}$. The proof in \cite[\href{https://stacks.math.columbia.edu/tag/0573}{Lemma 0573}]{stacks-project} which establishes the existence of a non-empty $V$ is constructive, except that one needs a constructive version of the ``generic flatness'' result that if $A$ is an integral finite type algebra and $M$ is a finitely-generated $A$-module then $M$ is flat on an open neighbourhood of $\Spec A$. For this one can use the stronger property that $M$ is generically free, for which a constructive proof appears in \cite{blechschmidt2018elementary}.
\end{itemize}
\end{proof}

\begin{lem}
\label{compliesinfam}
There exists an algorithm to compute the locus described by \autoref{liesinotherfam} given the families $g_{1}$ and $g_{2}$ as input.
\end{lem}

\begin{proof}
As in the proof of \autoref{liesinotherfam}, we may construct a family $g_{12} : X_{12} \to T_{1} \times_{K} T_{2}$ whose fibre above $(t_{1}, t_{2})$ may be identified with the intersection $X_{1,t_{1}} \cap X_{2,t_{2}}$. We then use the algorithm \autoref{familyalgos}(ii) to find the locus where the fibres of $g_{12}$ have dimension $d$, and project this locus to $T_{1}$. 
\end{proof}

\subsection{Computing with Type Data}

To compute the differential constraints $\mathcal{T}^{d}_{r}(\mathcal{C}_{i})$ for the finitely many types $\mathcal{C}_{1}, \hdots, \mathcal{C}_{q}$ associated to the variation $\mathbb{V}$, we must first explain how to compute certain data associated to those types. We will continue with the notation established in \autoref{modulisec}; in particular, we have a fixed polarized lattice $V$, and fixed embeddings $D \subset \ch{D} \subset \ch{L}$.

\begin{defn}
We say a type $\mathcal{C}$ is \emph{Hodge-theoretic} if $\mathcal{C} = \mathcal{C}(V)$ where $V = N \cdot t$ for a point $t \in D$ and a connected $\mathbb{Q}$-algebraic normal subgroup $N \subset \textrm{MT}(t)$.
\end{defn}

\begin{prop}
\label{computetypes}
There exists an algorithm that computes a finite list $V_{1}, \hdots, V_{\ell} \subset \ch{L}$ of subvarieties such that the set of Hodge-theoretic types is a subset of $\{ \mathcal{C}(V_{1}), \hdots, \mathcal{C}(V_{\ell}) \}$.
\end{prop}

\begin{proof}
We first note that the task of computing \emph{Hodge-theoretic} types reduces to the following problem: compute pairs $(h, N)$ up to $\GL(V)(\mathbb{C})$-equivalence, with $h \in D$ and $N$ a $\mathbb{Q}$-algebraic normal subgroup of $\textrm{MT}(h)$. Here we regard $\GL(V)(\mathbb{C})$ acting on $h$ through the embedding $D \subset \ch{L}$ and on $N$ by conjugation. Given that $(N', h') = g \cdot (N, h)$ for $g \in \GL(V)(\mathbb{C})$, we have as orbits in $\ch{L}$
\[ N' \cdot h' = (g \cdot N \cdot g^{-1}) \cdot (g \cdot h) = g \cdot (N \cdot h) , \]
so two equivalent pairs determine the same type.

We will denote by $\{ h^{p,q} \}_{p+q=n}$ the Hodge numbers associated to the period domain $D$. Given a point $h \in D$ thought of as a morphism $h : \mathbb{S} \to \GL(V)_{\mathbb{R}}$, the Mumford-Tate group is the $\mathbb{Q}$-Zariski closure of $h(\mathbb{S})$. Using the construction of $\mathbb{S}$ in \cite[Ch. 1]{GGK}, there is an isomorphism $\mathbb{S}_{\mathbb{C}} \simeq \mathbb{G}_{m,\mathbb{C}} \times \mathbb{G}_{m,\mathbb{C}}$, and pulling back the morphism $h_{\mathbb{C}}$ along the inclusion of the second factor we obtain a cocharacter $\mu_{h} : \mathbb{G}_{m,\mathbb{C}} \to \GL(V)_{\mathbb{C}}$. The weights of $\mu_{h}$ determine the Hodge flag associated to $h$, in the sense that $\mu_{h}$ acts through the character $z^{p}$ with multiplicity $h^{p,q}$, and the weight space decomposition on $V_{\mathbb{C}}$ is the Hodge decomposition associated to $h$.

The preceding paragraph shows that pairs $(h, N)$ are represented up to $\GL_{m}(\mathbb{C})$-equivalence by pairs $(\mu, N')$ with the following properties:
\begin{itemize}
\item[(i)] the weight space decomposition of $V_{\mathbb{C}}$ with respect to $\mu$ has the weight $z^{p}$ occur with multiplicity $h^{p,q}$;
\item[(ii)] the group $N'$ is a (possibly complex) connected semisimple algebraic subgroup of $\GL(V)_{\mathbb{C}}$;
\item[(iii)] the character $\mu$ factors through the connected component of the normalizer of $N'$.
\end{itemize}
To see that pairs $(\mu, N')$ satisfying (i), (ii) and (iii) are finite up to $\GL_{m}(\mathbb{C})$-conjugacy, we may reason as follows. First, we use the known fact that semisimple subgroups of $\GL(V)_{\mathbb{C}}$ are finite up to conjugacy, hence this is also true for the identity components of their normalizers. It then suffices to show that cocharacters $\mu : \mathbb{G}_{m, \mathbb{C}} \to H$ of a complex algebraic group $H$ with a fixed set of weights and multiplicities lie inside finitely many $H$-conjugacy classes. This reduces to the case where $H$ is a torus, where it is immediate. (For a similar argument, see the proof in \cite{hodgelocivoisin} that there are finitely many Mumford-Tate groups up to $\GL(V)(\mathbb{C})$-conjugacy.)

The proof is now completed as follows. In \cite{de2011constructing} de Graaf gives an algorithm that classifies all semisimple subalgebras of $\mathfrak{gl}(V)_{\mathbb{C}}$ up to linear equivalence, hence using \cite{de2009constructing} we may obtain representatives $N'_{1}, \hdots, N'_{\ell}$ for all connected semisimple subgroups of $\GL(V)_{\mathbb{C}}$ up to conjugacy. Computing the identity component of the normalizers $H'_{1}, \hdots, H'_{\ell}$ and fixing maximal tori $T'_{i} \subset H'_{i}$ for each $i$, we may then compute appropriate cocharacters $\mu_{1i}, \hdots, \mu_{k_{i} i}$ of $T'_{i}$. Then the desired set of types can be taken to be $\mathcal{C}(N'_{i} \cdot \mu_{ji})$ for indices $1 \leq i \leq \ell$ and $1 \leq j \leq k_{i}$.

\end{proof}

The types $\mathcal{C}(V_{i})$ computed by \autoref{computetypes} need not all come from the variation $\mathbb{V}$, in the sense that they may not be among the types $\mathcal{C}_{1}, \hdots, \mathcal{C}_{q}$ associated to subvarieties $Z \subset S$. As a consequence, to implement the ideas present in the previous sections we will require a way of testing when a certain differential constraint is the ``right one'' for a family of subvarieties, for which the following two lemmas will be useful:

\begin{lem}
\label{identifyflattensors}
Suppose $Z \subset S$ is a geometrically irreducible subvariety, and $U \subset Z$ is a Zariski open subset of the smooth locus. Let $\mathbb{V}^{m,n}_{U}$ be the restriction to $U$ of $\mathbb{V}^{\otimes m} \otimes (\mathbb{V}^{*})^{\otimes n}$, and let $\mathcal{H}^{m,n}_{U}$ be the associated algebraic Hodge bundle. Then the natural map $\mathbb{V}^{m,n}_{U, \mathbb{C}} \to \mathcal{H}^{m,n}_{U,\textrm{an}}$ identifies global monodromy invariant sections and algebraic sections in the kernel of $\nabla$. 
\end{lem}

\begin{proof}
The required result follows from the Riemann-Hilbert correspondence: global sections of $\mathbb{V}^{m,n}_{U, \mathbb{C}}$ may be identified with morphisms $\mathbb{C} \to \mathbb{V}^{m,n}_{U, \mathbb{C}}$ from the trivial local system, which correspond exactly to connection-compatible morphisms $\mathcal{O}_{U} \to \mathcal{H}^{m,n}_{U}$ of (algebraic, by Riemann-Hilbert) vector bundles.
\end{proof}

\begin{defn}
Given a type $\mathcal{C}(V)$ with $V \subset \ch{L}$ a subvariety, we refer to the dimension $\dim V$ as the dimension of $\mathcal{C}(V)$.
\end{defn}

\begin{defn}
Suppose that $g : X \to T$ is a subfamily of the universal family over $\textrm{Var}(S)$, with $T \subset \textrm{Var}(S)$ a geometrically irreducible closed subscheme. By the type $\mathcal{G}(T)$ of the generic fibre of $g$, we mean the type of a fibre $X_{t}$ above a point $t$ lying in the complement of those loci $\mathcal{W}(\mathcal{C}_{i})$ for $i = 1, \hdots, q$ which intersect $T$ properly. That this complement is a non-empty open subset follows from the closedness of the loci $\mathcal{W}(\mathcal{C}_{i})$ established in \autoref{typeinequalityloci}.
\end{defn}

\begin{lem}
\label{computepreimages}
Suppose we are given a $K$-algebraic family $g : X \to T$ over a geometrically irreducible reduced Noetherian base $T$ which is a subfamily of the universal family over $\textrm{Var}(S)$, along with the natural projection $p : X \to S$. Then there exists an algorithm to compute $\dim \mathcal{G}(T)$.
\end{lem}

\begin{proof}
Fix $(m,n)$ large enough so that any semisimple subgroup $N \subset \GL(V)_{\mathbb{C}}$ is determined by its tensor invariants of bidegree at most $(m,n)$; since \cite{de2011constructing} gives an algorithm for determining semisimple subalgebras of $\mathfrak{gl}(V)_{\mathbb{C}}$ up to linear equivalence, such an $(m,n)$ can be determined algorithmically. Given a vector bundle $\mathcal{V}$ on $S$, we denote by $\mathcal{V}_{g}$ its pullback under $p$ to $X$; this we can regard as a family of vector bundles on the fibres of $g$. We may compute the various tensor powers $\mathcal{H}^{j,h}_{g} = \left(\mathcal{H}^{\otimes j} \otimes (\mathcal{H}^{*})^{h} \right)_{g}$ for $0 \leq j \leq m$ and $0 \leq h \leq n$, as well as the subsheaf $\mathcal{I}^{j,h}_{g} \subset \mathcal{H}^{j,h}_{g}$ of flat tensors. We may also compute the locus $X_{\textrm{sm}} \subset X$ where the morphism $g$ is smooth. 

For each point $t \in T(\mathbb{C})$, let us explain the relationship between the bundles $\restr{\mathcal{I}^{j,h}_{g}}{X_{\textrm{sm}, t}}$ and the algebraic monodromy group of the fibre $X_{t}$; we will denote the bundles $\restr{\mathcal{I}^{j,h}_{g}}{X_{\textrm{sm}, t}}$ and $\restr{\mathcal{H}^{j,h}_{g}}{X_{\textrm{sm}, t}}$ simply by $\mathcal{I}^{j,h}$ and $\mathcal{H}^{j,h}$ in what follows. By \autoref{identifyflattensors}, the canonical isomorphism $\sigma : \mathcal{H}^{\textrm{an}} \xrightarrow{\sim} \mathbb{V} \otimes \mathcal{O}_{X_{\textrm{sm}, t}^{\textrm{an}}}$ identifies the bundles $\mathcal{I}^{j,h}$ with the subbundles generated by monodromy invariant tensors. Thus, it identifies the algebraic monodromy groups $\mathbf{H}_{X_{t}}$ with appropriate stabilizer groups of the tensors lying in the bundles $\mathcal{I}^{j,h}$ for $0 \leq j \leq m$ and $0 \leq h \leq n$. Moreover, if we regard $\mathbf{H}_{X_{t}}$ as a subgroup of $\GL(\mathbb{V}_{s})$ for some $s \in X_{t}$, the dimension $\dim \mathcal{C}(X_{t})$ is equal to $\dim \mathbf{H}_{X_{t}} - \dim (\textrm{Stab}(F^{\bullet}_{s}) \cap \mathbf{H}_{X_{t}})$, where $\textrm{Stab}(F^{\bullet}_{s})$ is the stabilizer of the Hodge flag at $s$.

Compute a generic point $t \in T(\mathbb{C})$; this we can do by computing $\mathcal{I}^{j,h}$ over the generic point $\eta \in T$ (i.e., over the function field $K(T)$) and then specializing. Choose a point $s \in X_{t,\textrm{sm}}(\mathbb{C})$. We may then compute the subgroup $N \subset \GL(\mathcal{H}_{s})$ which stabilizes $\mathcal{I}^{j,h}_{s}$ for all $0 \leq j \leq m$ and $0 \leq h \leq n$, and the subgroup $P \subset N$ which stabilizes $F^{\bullet}_{s}$. Applying \autoref{identifyflattensors} and the fact that $\sigma$ preserves filtrations, we get
\[ \dim N - \dim P = \dim \mathbf{H}_{X_{t}} - \dim (\textrm{Stab}(F^{\bullet}_{s}) \cap \mathbf{H}_{X_{t}}) = \dim \mathcal{C}(X_{t}) = \dim \mathcal{G}(T) . \]
\end{proof}

\subsection{Computing Hilbert Schemes}

\begin{lem}
\label{computehilbpolys}
There exists finitely many Hilbert polynomials $P_{1}, \hdots, P_{m}$ associated to pure dimensional reduced subschemes of $\overline{S}$ of degree at most $d$. Moreover, there exists an algorithm to compute the polynomials $P_{1}, \hdots, P_{m}$.
\end{lem}

\begin{proof}
The finiteness claim is proven by Starr in \cite{264428}. To see that the relevant polynomials are computable, we may inspect Starr's proof: for varieties of pure dimension $k$ and for each $m$ such that $n-k \leq m \leq {{n+d} \choose n}$ he considers the natural projective family over the Grassmannian of $m$-dimensional subspaces of $H^{0}(\mathbb{P}^{n}, \mathcal{O}(d))^{*}$, and argues that all of the polynomials $P_{1}, \hdots, P_{m}$ occur among the fibres of these families. Constructing these families for each $k$ and each $m$, we may use \autoref{familyalgos}(i) above to compute flattening stratifications, and then obtain all the polynomials $P_{1}, \hdots, P_{m}$ by computing the reduced, geometrically irreducible locus of each family using \autoref{familyalgos}(ii) and \autoref{familyalgos}(iii), and if this locus is non-empty, choosing a fibre and computing the Hilbert polynomial using \cite{bayer1992computation}.
\end{proof}

\begin{lem}
\label{univfamofvars}
There exists an algorithm to compute the subscheme $\textrm{Var}(S)_{d} \subset \textrm{Var}(S)$ of $\textrm{Var}(S) \subset \textrm{Hilb}(\overline{S})$ parametrizing varieties with closure in $\overline{S}$ of degree at most $d$, as well as the universal family $u_{d} : \mathcal{Y}_{d} \to \textrm{Var}(S)_{d}$ and the projection $\mathcal{Y}_{d} \to S$.
\end{lem}

\begin{proof}
By \autoref{computehilbpolys}, we may enumerate the Hilbert polynomials $P_{1}, \hdots, P_{m}$ corresponding to pure dimensional geometrically reduced subschemes of degree at most $d$. For each polynomial $P_{i}(t)$, one can construct the associated component $\textrm{Hilb}(\overline{S})_{i}$ of the Hilbert scheme using subvarieties of appropriate Grassmanians described explicitly in the Pl\"uker coordinates; see for instance \cite[\S1]{lella2012computable}. Note that such methods naturally also compute the associated family $\overline{\mathcal{Y}}_{i} \to \textrm{Hilb}(\overline{S})_{i}$ and the associated projection, as they work with explicit families of generators for homogeneous ideals in $\mathbb{P}^{n}$ and so produce equations for the universal family inside $\mathbb{P}^{n} \times \textrm{Hilb}(\overline{S})_{i}$. For the geometric irreducibility and reducedness we use \autoref{familyalgos}(ii) and \autoref{familyalgos}(iii) above, and to impose the condition that a scheme $Y \subset \overline{S}$ intersects $S$ it suffices to require that $Y$ does not lie inside $\overline{S} \setminus S$, which we can do using \autoref{compliesinfam}.
\end{proof}

\subsection{Computing Loci Satisfying Differential Constraints}

Let us suppose that we have computed the varieties $V_{1}, \hdots, V_{k} \subset \ch{L}$ given by \autoref{computetypes}, and denote by $\mathcal{C}(V_{1}), \hdots, \mathcal{C}(V_{k})$ the associated types. We know that the types $\mathcal{C}_{1}, \hdots, \mathcal{C}_{q}$ associated to subvarieties of $S$ are a subset of the types $\mathcal{C}(V_{1}), \hdots, \mathcal{C}(V_{k})$, but we don't know which subset this is. To remedy this we will work with he following loci:

\begin{defn}
For any locus $\mathcal{W}(\mathcal{C})$ (recall \autoref{weaksplocsec}) associated to a type $\mathcal{C}$, define the locus $\mathcal{W}(\mathcal{C})_{\textrm{opt}}$ to be the sublocus of $\mathcal{W}(\mathcal{C})$ consisting of just those components $C \subset \mathcal{W}(\mathcal{C})$ such that we have an equality of types $\mathcal{G}(C) = \mathcal{C}$. We refer to such components $C$ as \emph{optimal} components; this notion is distinct from other notions of optimality in theories of exceptional intersections. We note that if $\mathcal{C} = \mathcal{C}(V)$ for $V$ irreducible, this is the same condition as $\dim \mathcal{G}(C) = \dim \mathcal{C}(V)$. 
\end{defn}

\begin{lem}
\label{optnonemptylem}
The locus $\mathcal{W}(V_{i})_{\textrm{opt}}$ is non-empty if and only if the type $\mathcal{C}(V_{i})$ is among the types $\mathcal{C}_{1}, \hdots, \mathcal{C}_{q}$ associated to subvarieties $Z \subset S$.
\end{lem}

\begin{proof}
Given any subfamily $g : X \to C$ of the universal family over $\textrm{Var}(S)$, where $C \subset \textrm{Var}(S)$ is closed and irreducible, the generic type $\mathcal{G}(C)$ satisfies the property that for each $c \in C$ we have $\mathcal{C}(X_{c}) \leq \mathcal{G}(C)$. (This is because $C \subset \mathcal{W}(\mathcal{G}(C))$ since the latter set is closed.) Thus if we consider a component $C$ of $\mathcal{W}(\mathcal{C}(Z))$ which contains $[Z]$, then necessarily $\mathcal{C}(Z) \leq \mathcal{G}(C) \leq \mathcal{C}(Z)$, and so $\mathcal{W}(\mathcal{C}(Z))_{\textrm{opt}}$ contains $C$, and is non-empty. 

Conversely, supposing that $\mathcal{W}(V_{i})_{\textrm{opt}}$ is non-empty containing the component $C$, we learn that for a generic point $c \in C$ we have $\mathcal{C}(X_{c}) = \mathcal{C}(V_{i})$, hence the result.
\end{proof}



\vspace{1em}

\noindent We now explain how to compute the loci $\mathcal{W}(\mathcal{C}(V_{i}))_{\textrm{opt}}$.


\begin{lem}
\label{algtocomputesatlocus}
Given a $K$-algebraic family $g : X \to T$ which is a subfamily of the universal family over $\textrm{Var}(S)$, and given a type $\mathcal{C}(V)$, there exists an algorithm to compute the constructible locus of $[Z] \in T$ such that $Z$ satisfies the differential constraint $\mathcal{T}^{d}_{r}(\mathcal{C}(V))$.
\end{lem}

\begin{proof}
Applying \autoref{compdiffconst}, we may compute the differential constraints $\mathcal{T}^{d}_{r}(\mathcal{C}(V)) \subset J^{d}_{r} S$. We then observe that the proof of \autoref{constgivesconstrcond} is constructive, so it suffices to construct the appropriate families and judiciously apply \autoref{algbasicsprop}(ii).
\end{proof}

\begin{prop}
\label{computeoptprop}
Given a component $T$ of $\textrm{Var}(S)$ and a geometrically irreducible variety $V \subset \ch{L}$, there exists an algorithm to compute $T \cap \mathcal{W}(\mathcal{C}(V))_{\textrm{opt}}$.
\end{prop}

Before proceeding with the proof of \autoref{computeoptprop}, we note that \autoref{univfamofvars} above guarantees that we may assume we have computed the universal family $u : \mathcal{Y} \to T$ over $T$ along with the projection $p : \mathcal{Y} \to S$. We will fix $d$ to be the common dimension of the varieties parametrized by $T$ throughout. We denote by $\mathcal{W}^{d}_{r,i} \subset T$ the constructible locus of $[Z] \in T$ such that $Z$ satisfies the differential constraint $\mathcal{T}^{d}_{r}(\mathcal{C}(V_{i}))$. Note that it follows from \autoref{liesintrans} that $T \cap \mathcal{W}(\mathcal{C}(V_{i})) = \bigcap_{r} \mathcal{W}^{d}_{r,i}$.

\begin{proof}[Proof of \ref{computeoptprop}]
We set $V_{0} = V$. We will show the stronger claim that we may compute all the loci $T \cap \mathcal{W}(\mathcal{C}(V_{0}))_{\textrm{opt}}, T \cap \mathcal{W}(\mathcal{C}(V_{1}))_{\textrm{opt}}, \hdots, T \cap \mathcal{W}(\mathcal{C}(V_{k}))_{\textrm{opt}}$; in particular, we consider those indicies $i_{1}, \hdots, i_{p}$ with $0 \leq i_{j} \leq k$ such that $\dim V_{i_{j}} = \ell$, and show that we may compute the loci $T \cap \mathcal{W}(\mathcal{C}(V_{i_{1}}))_{\textrm{opt}}, \hdots, T \cap \mathcal{W}(\mathcal{C}(V_{i_{p}}))_{\textrm{opt}}$. 

We compute the loci $T \cap \mathcal{W}(\mathcal{C}(V_{i_{p}}))_{\textrm{opt}}$ in parallel. More specifically, we compute for each $i_{j}$ the locus $\mathcal{W}^{d}_{r, i_{j}}$. We observe that there must be some $r_{0}$ at which:
\begin{itemize}
\item[(i)] the loci $\mathcal{W}^{d}_{r_{0}, i_{j}}$ are all closed;
\item[(ii)] if $j \neq j'$ and $C \subset \mathcal{W}^{d}_{r_{0},i_{j}}$ and $C' \subset \mathcal{W}^{d}_{r_{0},i_{j'}}$ are two components such that $\dim \mathcal{G}(C) = \dim \mathcal{C}(V_{i_{j}})$ and $\dim \mathcal{G}(C') = \dim \mathcal{C}(V_{i_{j'}})$, then $C$ and $C'$ intersect properly. 
\end{itemize}
Indeed, the existence of such an $r_{0}$ follows because this is true in the case where $\mathcal{W}^{d}_{r_{0},i_{j}} = T \cap \mathcal{W}(\mathcal{C}(V_{i_{j}}))$ for all $i_{j}$, which occurs for some $r_{0}$ large enough as the loci $\mathcal{W}(\mathcal{C}(V_{i_{j}}))$ are all closed. (Note that the condition $\dim \mathcal{G}(C) = \dim \mathcal{C}(V_{i_{j}})$ in fact implies $\mathcal{G}(C) = \mathcal{C}(V_{i_{j}})$ for components $C$ of $\mathcal{W}(\mathcal{C}(V_{i_{j}}))$, hence optimal components of $\mathcal{W}(\mathcal{C}(V_{i_{j}}))$ and $\mathcal{W}(\mathcal{C}(V_{i_{j}}))$ must intersect properly for $j \neq j'$ otherwise we would have $\mathcal{C}(V_{i_{j}}) = \mathcal{C}(V_{i_{j'}})$.) 

We claim that at this stage we must have that the optimal components of $T \cap \mathcal{W}(\mathcal{C}(V_{i_{j}}))$ are exactly those components $C$ of $\mathcal{W}^{d}_{r_{0}, i_{j}}$ whose generic type has dimension $\dim V_{i_{j}}$. Indeed, we know that for each such $C$, it must lie inside an optimal component of $T \cap \mathcal{W}(\mathcal{C}(V_{i_{j'}}))$ for some $j'$. But the choice of $r_{0}$ guarantees that $C$ lies inside $\mathcal{W}(\mathcal{C}(V_{i_{j'}}))$ only when $j' = j$, so the result follows.

We conclude by observing that we can compute all the loci $\mathcal{W}^{d}_{r,i_{j}}$ and the optimal components $C$ by computing the differential conditions $\mathcal{T}^{d}_{r}(\mathcal{C}(V_{i}))$ and using \autoref{algtocomputesatlocus}.
\end{proof}

\subsection{Main Result}

\begin{thm}
\label{mainalgres}
There exists an algorithm to compute the intersection of $\mathcal{W}$ with $\textrm{Var}(S)_{d}$.
\end{thm}

\begin{proof}
By \autoref{computeoptprop} and \autoref{univfamofvars}, we may compute the loci $\mathcal{W}(\mathcal{C}(V_{i}))_{\textrm{opt}} \cap \textrm{Var}(S)_{d}$ for each of the types $\mathcal{C}(V_{1}), \hdots, \mathcal{C}(V_{k})$. By \autoref{optnonemptylem} these loci will be non-empty only for types belonging to the subset $\mathcal{C}_{1}, \hdots, \mathcal{C}_{q}$, so we may assume that all the types we work with belong to this set. It suffices to prove the claim of the theorem for the locus $\mathcal{W}_{j}$ of weakly special subvarieties of type $\mathcal{C}_{j}$. Since $\mathcal{W}_{j} \subset \mathcal{W}(\mathcal{C}_{j})_{\textrm{opt}}$, we may therefore fix a component $C \subset \mathcal{W}(\mathcal{C}_{j})_{\textrm{opt}} \cap \textrm{Var}(S)_{d}$ parametrizing varieties of dimension $k$ and show we can compute $\mathcal{W}_{j} \cap C$. Removing points $[Z] \in C$ which lie inside $\mathcal{W}(\mathcal{C}_{j'})_{\textrm{opt}}$ for some $j' \neq j$ such that the intersection $\mathcal{W}(\mathcal{C}_{j'}) \cap \mathcal{W}(\mathcal{C}_{j})$ is proper, we obtain the locus $C^{\circ} \subset C$ which appears in the proof of \autoref{Wisconstrthm}. We then observe that:
\begin{itemize}
\item[(i)] By \autoref{computeoptprop} above, we have an algorithm which computes any component $T_{i}$ of $\mathcal{W}(\mathcal{C}_{j})_{\textrm{opt}}$, and in particular those components parametrizing $[Z] \in \textrm{Var}(S)$ with $\dim Z > k$. By applying \autoref{compliesinfam} to the universal family over $C^{\circ}$ as well as the universal family over $T_{i}$, we may therefore compute constructible loci $E_{i} \subset C^{\circ}$, where $E_{i}$ is as in the proof of \autoref{Wisconstrthm}.
\item[(ii)] We may compute the loci $\mathcal{U}(C^{\circ}, V_{j}, r)$ that appear in the proof of \autoref{Wisconstrthm} by observing that the arguments that appear in \autoref{germsofpermaps} and \autoref{provmaxsec} are constructive. In particular, the constructions that appear only require linear algebra and the tools present in \autoref{algbasicsprop}.
\end{itemize}
By constructibility of $\mathcal{W}_{j} \cap C^{\circ}$ and \autoref{unionintlem}, we know that there exists $i_{0}$ and $r_{0}$ such that
\[ \mathcal{W}_{j} \cap C = \bigcap_{i = 1}^{i_{0}} (C^{\circ} \setminus E_{i}) = \bigcup_{r = 1}^{r_{0}} \mathcal{U}(C^{\circ}, V_{j}, r) . \]
Moreover, the sets $E_{i}$ and $\mathcal{U}(C^{\circ}, V_{j}, r)$ are disjoint for each choice of $i$ and $r$. It follows that we may compute $\mathcal{W}_{j} \cap C^{\circ}$ by computing the decreasing intersections $\bigcap_{i = 1}^{i'} (C^{\circ} \setminus E_{i})$ and increasing unions $\bigcup_{r = 1}^{r'} \mathcal{U}(C^{\circ}, V_{j}, r)$, and terminating when we achieve equality.
\end{proof}

\section{Applications}

We now prove a stronger form of the conjecture that appears in \cite{daw2018applications}.

\begin{thm}
\label{conjugacythm}
Let $\mathbb{V}$ be a polarizable integral variation of Hodge structure on a smooth algebraic variety $S$. Fix a projective compactification $S \subset \overline{S}$ and an ample line bundle $\mathcal{L}$ on $\overline{S}$. Then as $Z \subset S$ ranges over (geometrically) irreducible subvarieties of $S$ whose closures $\overline{Z}$ in $\overline{S}$ have degree at most $d$ with respect to $\mathcal{L}$, the following holds:
\begin{itemize}
\item[(i)] up to conjugacy by the monodromy group $\Gamma_{S}$ of $\mathbb{V}$, only finitely many algebraic monodromy groups $\mathbf{H}_{Z}$ arise;
\item[(ii)] of those $Z$ that are weakly special, only finitely many are rigid (weakly non-factor in the language of \cite{fieldsofdef}), and each corresponds to an isolated point $[Z]$ in the weakly special locus $\mathcal{W} \subset \textrm{Var}(S)$;
\item[(iii)] the non-rigid weakly specials $Z$ belong to finitely many constructible families, with each family consisting of varieties lying inside a common rigid (weakly) special.
\end{itemize}
\end{thm}

\begin{proof}
By replacing $\mathcal{L}$ with a tensor power we may embed $\overline{S}$ in $\mathbb{P}^{n}$, and assume the degrees of the closures $\overline{Z}$ correspond with their degrees as subvarieties of $\mathbb{P}^{n}$. We are then in the setting of \autoref{modulisec}, from which the desired results may be obtained in the following way:
\begin{itemize}
\item[(i)] As there are finitely many types, it suffices to show the statement for those points $[Z] \in \textrm{Var}(S)_{d}$ of type $\mathcal{C}_{i}$, and hence those $[Z] \in \mathcal{W}(\mathcal{C}_{i}) \cap \textrm{Var}(S)_{d}$ that do not also lie in $\mathcal{W}(\mathcal{C}_{j})$ for some type $\mathcal{C}_{j} \leq \mathcal{C}_{i}$ with $j \neq i$. In this case the component $C$ of $\mathcal{W}(\mathcal{C}_{i})$ containing $[Z]$ will satisfy $\mathcal{G}(C) = \mathcal{C}(Z)$, and applying \autoref{typeinequalityloci} we may identify $C$ with a component of $\mathcal{D}(Y,\mathbf{H}_{Z}) \cap \textrm{Var}(S)_{d}$, where $Y$ is a rigid weakly special containing $Z$ such that $\mathbf{H}_{Z}$ is normal in $\mathbf{H}_{Y}$. By \autoref{closedfamlocusprop} a generic point $[Z'] \in C$ (outside of $\mathcal{W}(\mathcal{C}_{j})$ for $\mathcal{C}_{j} \leq \mathcal{C}_{i}$ and $j \neq i$) has algebraic monodromy group $\mathbf{H}_{Z'}$ conjugate to $\mathbf{H}_{Z}$ under the monodromy $\Gamma_{Y}$ on $Y$, and hence by the monodromy $\Gamma_{S}$ on $S$. Only finitely many such components occur for degrees $\leq d$, hence the result. 

\item[(ii)] Suppose that $[Z] \in \mathcal{W}(\mathcal{C}(Z))$ is a rigid weakly special, and let $C$ be a component of $\mathcal{W}(\mathcal{C}(Z))$ containing $Z$. By \autoref{typeinequalityloci} we may identify $C$ with a component of $\mathcal{D}(Y,N)$ for $Y$ a rigid weakly special and  and where $N = \mathbf{H}_{Z_{0}}$ is normal in $\mathbf{H}_{Y}$, with $[Z_{0}] \in C$ a point. On one hand, from \autoref{closedfamlocusprop} we know that $\mathcal{C}(Z) \leq \mathcal{C}(Y, N) = \mathcal{C}(Z_{0})$. On the other hand since $C \subset \mathcal{W}(\mathcal{C}(Z))$ we have $\mathcal{C}(Z_{0}) \leq \mathcal{C}(Z)$, so it follows that $\mathcal{C}(Z) = \mathcal{C}(Z_{0}) = \mathcal{C}(Y,N)$. 

Once again applying \autoref{closedfamlocusprop} we learn that $\mathbf{H}_{Z}$ lies in one of the $\mathbb{Q}$-algebraic subgroups $N_{1}, \hdots, N_{\ell}$ of $\mathbf{H}_{Y}$ which are conjugate to $\mathbf{H}_{Z_{0}}$ under $\Gamma_{Y}$; say that $\mathbf{H}_{Z} \subset N_{i}$. If we fix a local period map $\psi : B \to D$ such that $B$ contains a point $s \in Z(\mathbb{C})$, then we have that $\mathbf{H}_{Z} \cdot \psi(s) \subset N_{i} \cdot \psi(s)$. Using the fact that $\dim \mathcal{C}(Z) = \dim \mathcal{C}(Y,N)$ we learn that in fact $\mathbf{H}_{Z} \cdot \psi(s) = N_{i} \cdot \psi(s)$. Since $N_{i}$ is conjugate under $\Gamma_{Y}$ to $\mathbf{H}_{Z_{0}}$, it does not contain a $\mathbb{Q}$-simple factor with a compact set of real points by \autoref{nocomprealfac}. Applying \autoref{monocont} twice we learn that $\mathbf{H}_{Z} = N_{i}$, and hence $\mathbf{H}_{Z}$ is normal in $\mathbf{H}_{Y}$. By rigidity, it follows that $Y = Z$. Since all elements of $C$ lie in $Y$ and have the same dimension as $Z$, we learn that $C$ is a point.

We have learned that for each type $\mathcal{C}_{i}$ arising from the variation $\mathbb{V}$, every rigid weakly special $Z$ with $\mathcal{C}(Z) = \mathcal{C}_{i}$ is an isolated point in $\mathcal{W}(\mathcal{C}_{i})$. Since there can only be finitely many isolated points in the finite-type locus $\mathcal{W}(\mathcal{C}_{i}) \cap \textrm{Var}(S)_{d}$ and there are only finitely many types $\mathcal{C}_{1}, \hdots, \mathcal{C}_{q}$ arising from $\mathbb{V}$, the result follows.

\item[(iii)] This is the statement that the locus $\mathcal{W}_{i} \cap \textrm{Var}(S)_{d}$ is covered by finitely many of the loci $\mathcal{D}(Y, N)$ of \autoref{closedfamlocusprop} and \autoref{typeinequalityloci}.
\end{itemize}
\end{proof}

As a corollary, we obtain the following conjecture which appears in \cite{daw2018applications}

\begin{cor}
\label{dawrencor}
Suppose that $S = \Gamma \backslash X$ is a Shimura variety associated to the Shimura datum $(G, X)$. Then there exists a finite set $\Omega$ of semisimple $\mathbb{Q}$-algebraic subgroups of $G$ such that if $Z \subset S$ is a special subvariety of $S$ with $\deg(Z) \leq d$ (relative to the Bailey-Borel line bundle) and defined by a Shimura subdatum $(H, X_{H}) \subset (G, X)$, then
\[ H^{\textrm{der}} = \gamma F \gamma^{-1} , \]
for some $\gamma \in \Gamma$ and $F \in \Omega$.
\end{cor}

\begin{proof}
Without loss of generality, passing to a finite covering if necessary, we can assume $\Gamma \subset G(\mathbb{Q})$ is a neat arithmetic subgroup and so $S$ is a smooth quasi-projective algebraic variety. If we choose a faithful integral representation $\rho : G \to \GL(V)$, we obtain a variation $\mathbb{V}$ of polarizable integral Hodge structures on $S$ with monodromy $\Gamma_{S} = \Gamma$. For each special subvariety $Z$ defined by $(H, X_{H})$, the group $H$ may be identified with the Mumford-Tate group of $Z$, and the group $H^{\textrm{der}}$ with its algebraic monodromy group. We thus may take the set $\Omega$ to be the set of algebraic monodromy groups associated to subvarieties $Z \subset S$ of degree at most $d$, and the result follows from \autoref{conjugacythm}(i) above.
\end{proof}

As explained in \cite[\S10]{daw2018applications} and \cite[Rem. 3.8]{daw2020effective}, \autoref{dawrencor} also resolves \cite[Conj. 10.3]{daw2018applications} of Daw-Ren and the conjecture in \cite{daw2020effective} of Daw-Javanpekkar-K\"uhne that there should be finitely many non-facteur special subvarieties of bounded degree.

Lastly, we resolve \autoref{mainthmspcor}.

\begin{cor}
In the setting of \autoref{mainalgres}, there exists an algorithm to compute a finite subset $\mathcal{S} \subset \mathcal{W} \cap \textrm{Var}(S)_{d}$ containing the rigid weakly special subvarieties.
\end{cor}

\begin{proof}
By \autoref{conjugacythm}(ii) it suffices to compute the finite set of isolated points of $\mathcal{W} \cap \textrm{Var}(S)_{d}$, which can be done by taking the Zariski closure and computing an irreducible decomposition.
\end{proof}

\bibliography{hodge_theory}
\bibliographystyle{alpha}

\end{document}